\documentclass[12pt]{article}
\newcommand{\nl}{\mbox{} \\}

\newcommand{\IR}{\hbox{\rm I\kern-.200em R}}             

\newcommand{\ep}{\varepsilon}

\newtheorem{lemma}{Lemma}[section]
\newtheorem{theorem}{Theorem}[section]
\newtheorem{assumption}{Assumption}[section]
\newtheorem{definition}{Definition}[section]
\newcommand{\ba}{\begin{array}}
\newcommand{\ea}{\end{array}}
\newcommand{\be}{\begin{equation}}
\newcommand{\ee}{\end{equation}}
\newcommand{\bea}{\begin{eqnarray}}
\newcommand{\eea}{\end{eqnarray}}
\newcommand{\beaa}{\begin{eqnarray*}}
\newcommand{\eeaa}{\end{eqnarray*}}
\newcommand{\bi}{\begin{itemize}}
\newcommand{\ei}{\end{itemize}}
\newcommand{\bn}{\begin{enumerate}}
\newcommand{\en}{\end{enumerate}}
\newcommand{\bt}{\begin{theorem}}
\newcommand{\et}{\end{theorem}}
\newcommand{\bl}{\begin{lemma}}
\newcommand{\el}{\end{lemma}}

\newcommand{\f}{\frac}

\def\p2{2\pi}

\def\gp2{{1\over \p2}}

\def\'{^\prime}

\def\eps{\varepsilon}

\def\R    {\hbox{\rm I\kern-.200em R}}           

\def\eps {\varepsilon}                                   

\def\lam  {\lambda}
\def\al {\alpha}
\def\ni {\|_{\infty}}

\def\o0 {\omega^{(0)}}

\def\cD  {{\cal D}}

\def\qand{\quad\mbox{and}\quad}
\def\qwith {\quad\mbox{with}\quad}

\def\qif{\quad\mbox{if}\quad}

\def\qat{\quad\mbox{at}\quad}

\def\qas{\quad\mbox{as}\quad}
\def\qfor{\quad\mbox{for}\quad}

\def\bt {\begin{theorem}}
\def\bl {\begin{lemma}}
\def\bd {\begin{definition}}
\def\ba{\begin{assumption}}
\def\bp {\begin{proof}}
\def\et {\end{theorem}}
\def\el {\end{lemma}}
\def\ed {\end{definition}}
\def\ea {\end{assumption}}
\def\ep {\end{proof}}
\def\f {\frac}

\def\Box {\diamond}

\begin{document}

\mbox{} \vspace{-1.750cm} \\

\begin{center}
\mbox{} \hspace{-0.900cm}
\mbox{\Large \bf
The Navier-Stokes Equations for Incompressible Flows:} \\
\mbox{} \vspace{-0.200cm} \\
{\Large \bf
Solution Properties at Potential Blow-up Times} \\
\nl
\mbox{} \vspace{-0.100cm} \\
{\large \sc Jens Lorenz} \\
\mbox{} \vspace{-0.475cm} \\
{\normalsize Department of Mathematics and Statistics} \\
\mbox{} \vspace{-0.590cm} \\
{\normalsize University of New Mexico} \\
\mbox{} \vspace{-0.600cm} \\
{\normalsize Albuquerque, NM 87131, USA} \\
\nl
{\large \sc Paulo R. Zingano} \\
\mbox{} \vspace{-0.475cm} \\
{\normalsize Departamento de Matematica Pura e Aplicada} \\
\mbox{} \vspace{-0.590cm} \\
{\normalsize Universidade Federal Do Rio Grande Do Sul} \\
\mbox{} \vspace{-0.600cm} \\
{\normalsize Porto Alegre, RS 91509, Brasil} \\
\end{center}

\mbox{} \vspace{-1.150cm} \\

\begin{abstract}
In this paper we consider the Cauchy problem for
the 3D Navier--Stokes equations for incompressible flows.
The initial data are assumed to be smooth and
rapidly decaying at infinity.
A famous open problem is whether
classical solutions can develop singularities in finite time.
Assuming the maximal interval of existence to be finite,
we give a unified discussion of various known
solution properties \linebreak
as time approaches the blow-up time.
\end{abstract}

\mbox{} \vspace{-0.750cm} \\

\noindent
{\bf Key words:}
Navier--Stokes equations, incompressible flows, blow--up

\noindent
{\bf AMS subject classification}: 35G20, 35Q30, 76D03, 76D05

\noindent
{\bf Abbreviated Title:} Incompressible Navier--Stokes

%
%

\mbox{} \vspace{-0.750cm} \\

%
%
\begin{section}{Introduction}
\label{s1}
\mbox{} \vspace{-1.050cm} \\

In this paper we consider the Cauchy problem for the
3D Navier--Stokes equations,
\be
\label{eq.1.1}
\label{Navier.Stokes.eqs}
u_t + u\cdot \nabla u + \nabla p = \Delta u,
\quad
\nabla \!\:\!\cdot u = 0,
\quad
u(x,0) = f(x) \qfor x\in \R^3 \!.
\ee
The current mathematical theory of the problem (\ref{eq.1.1})
remains fundamentally incomplete:
it is known that a weak solution exists for all time
$t\geq 0$
if
$ f \in L^{2}(\R^3) $,
$ \nabla \!\cdot\! f = 0 $,
but it is not known if weak solutions are unique.
This is recognized as a major open problem
since the fundamental paper of Leray \cite{Leray1934}
(for a brief account of Leray's works,
see \cite{BorelHenkinLax2000}).
On the other hand, if $f$ is more regular,
then a unique classical solution exists
in some maximal interval $\;\!0\leq t < T_{\!\;\!f}$,
but it is not known
if $T_{\!\;\!f}$
can be finite or is always infinite.
In other words,
it is not known if classical solutions
can break down in finite time.

The Navier--Stokes equations are of fundamental importance
in continuum mechanics.
When one derives the equations
from the balance laws of mass and momentum
and from principle assumptions relating the stresses to velocity gradients,
then smoothness of the solution is assumed.
To make the model of the Navier--Stokes equations self--consistent,
one would like to prove that singularities in the solution
do not develop in finite time,
from smooth initial data with finite energy.
Thus far, however, this aim has not been achieved.
It remains one of the fundamental open problems in nonlinear analysis,
being included in the Millenium Prize Problems by the
Clay Mathematics Institute \cite{Fefferman2006}.
In fact, \linebreak
it has invariably appeared in all major recent lists
of the most important problems of Modern Mathematics,
see e.g.$\;$\cite{Browder2000, Constantin2001, Jackson2000, Smale1998, Wiegner1999}.

In this paper we consider only classical solutions of the problem (\ref{eq.1.1}),
and under our assumptions on $f$ these will be  $C^\infty$ functions.
If one normalizes the pressure so that $ p(x,t) \to 0 $ as $ |\,x\,| \to \infty $,
then the solution $\bigl( u(x,t), p(x,t) \bigr)$ is unique.
Its maximal interval \linebreak
of existence is denoted by $ 0 \leq t < T_{\!\;\!f} $.
Assuming $ T_{\!\;\!f} $ to be finite,
certain norms of $u(\cdot,t)$ \linebreak
will tend to infinity as $t \to T_{\!\;\!f} \!\;\!$,
while other norms remain bounded.
This issue is well \linebreak
studied in the literature,
see e.g.
\cite{Benameur2010, FoiasTemam1979, Galdi2000, %
Giga1986, Heywood1990, Kato1984, Leray1934, %
RobinsonSadowski2012, RobinsonSadowskiSilva2012, %
Serrin1963},
but the results are somewhat scattered.
We will review some results
that we consider to be very important
and will also derive lower bounds for some blow--up rates.
The results considered in this review
all fit in our unified discussion
that requires only a small selected set of
relatively basic ideas.
Thus, in spite of their undisputable importance,
many fundamental results such as
\cite{BeiraoDaVeiga2000, CaffarelliKohnNirenberg1982, %
ConstantinFefferman1993, Prodi1959, Serrin1962, Struwe1988}
and their recent developments
will be left out.\footnote{%
%
%
For a nice discussion of the celebrated
Caffarelli--Kohn--Nirenberg's regularity result,
see \cite{FridPerepelitsa2001}.
}
%
%

Another issue is to compare two functionals of $u(\cdot,t)$
that blow up as $ t \to T_{\!\;\!f} $.
Which one will blow up faster?
We believe that a better understanding of this issue is important
for further progress on the blow--up question and
recall some simple results in Section \ref{s5}.

The intent of this paper is to give rather complete proofs of
some solution properties for (\ref{eq.1.1}) that must hold,
as $t$ approaches $T_{\!\;\!f} $,
{\it if $ \,T_{\!\;\!f} $ is finite}.
These results and their proofs may be helpful
if one wants to construct a solution that actually does blow up.
They may also be helpful to show that blow--up is ultimately impossible.

For a treatment of many recent developments
regarding the Navier--Stokes equations
using methods of modern analysis,
the reader is referred to \cite{Cannone2004, Lemarie2002, Sohr2001}.
For similar blow--up questions
concerning the related Euler equations,
see e.g. \cite{BealeKatoMajda1984, Chae2008, %
LopesNussenzveigZheng1999, Majda1982, MarchioroPulvirenti1994}.

For simplicity of presentation
(and to avoid unessential complications near $t = 0$),
we put strong smoothness and spatial decay assumptions
on the initial state $ f $
and require \linebreak
(as in \cite{Fefferman2006})
that $ \!\;\!f $ is a divergence--free $C^\infty\!\;\!$ function
with all of its derivatives in $L^2(\R^3)$,
i.e.,
we assume that
$ f \in H^{n}(\R^3) $ for all $n$,
with $ \nabla \!\cdot\! f = 0 $. \\
\mbox{} \vspace{-0.940cm} \\
%

%
%

As to our {\bf notation},
we will be using the following standard definitions: \\
\mbox{} \vspace{-1.100cm} \\
\beaa
|\,u\;\!|^2 &=& u_1^2+u_2^2+u_3^2 \qfor u = (u_1, u_2, u_3) \in \R^3 \!, \\
|\;\!\alpha\;\!|\;&=& \alpha_1 + \alpha_2 + \alpha_3
\quad \mbox{for~a~multi--index} \quad \alpha =(\al_1,\al_2,\al_3), \\
D^\alpha \; &=& D_1^{\alpha_1} D_2^{\alpha_2} D_3^{\alpha_3},
\quad D_j =\, \partial/\partial x_j,
\quad \alpha =(\al_1,\al_2,\al_3), \\
\|\,u(\cdot,t)\,\|_{\mbox{}_{\scriptstyle L^q}} &=&\!\!
\biggl(\;\! \sum_{i\,=\,1}^{3}
\int_{\mbox{}_{\mbox{\scriptsize $\!\;\!\R^3$}}} \!\!\!\!\;\!
|\,u_{i}(x,t)\,|^q \, dx\:\!\biggr)^{\!\!\:\!1/q}\!\!\;\!,
\quad 1 \leq q < \infty,
\quad u = (u_1, u_2, u_3), \\
\|\,u(\cdot,t)\,\| \;\;\; &=&
\|\,u(\cdot,t)\,\|_{\mbox{}_{\scriptstyle L^2}}
\!\;\!,\\
\|\,u(\cdot,t)\,\|_\infty \;\!&=& \sup
\:\bigl\{\, |\,u_i(x,t)\,| \!:\;\!
x \in \R^3\!\;\!, \, \mbox{\small $1 \leq i \leq 3$}\,\bigr\}, \\
\|\,\cD^n u(\cdot, t)\,\|_{\mbox{}_{\scriptstyle L^q}} \!\!&=&\!\!
\biggl(\;\! \sum_{i\,=\,1}^{3} \sum_{\;\!j_{\mbox{}_{1}}=\,1}^{3} \!\cdots\!
\sum_{j_{\mbox{}_{n}}=\,1}^{3}
\int_{\mbox{}_{\mbox{\scriptsize $\!\;\!\R^3$}}} \!\!\!\!\;\!
|\;\!D_{\scriptstyle \!\;\!j_{\mbox{}_{1}}}
\!\!\:\!\cdots
D_{\scriptstyle \!\;\!j_{\mbox{}_{n}}}
u_{i}(x,t)\:|^q\, dx\:\!\biggr)^{\!\!\:\!1/q}\!\!\;\!,
\quad 1 \leq q < \infty, \\
\|\,\cD^n u(\cdot,t)\,\| \;\,&=&
\|\,\cD^n u(\cdot,t) \,\|_{\mbox{}_{\scriptstyle L^2}}
\!\;\!, \\
\|\,\cD^n u(\cdot,t)\,\|_\infty \!\!\:\!&=&
\sup \:\bigl\{\, |\;\! D^\alpha u_{i}(x,t)\,| \!\:\!:\;\!
x \in \R^3 \!\;\!, \,\mbox{\small $1 \leq i \leq 3$}, \,
|\,\mbox{\small $\alpha$}\,| = \mbox{\small $n$} \;\!\bigr\}, \\
(u,v) \;\;&=& \!\sum_{j\,=\,1}^3
\int_{\mbox{}_{\mbox{\scriptsize $\!\;\!\R^3$}}}
\!\!\!\!\;\!
u_j(x) \, v_j(x)\, dx\:\!.
\eeaa
\mbox{} \vspace{-0.300cm} \\
($\,\!$Note that
$ {\displaystyle
\|\, u(\cdot,t) \,\|_{\mbox{}_{\scriptstyle L^{q}}} \!\to\!\:\!
\|\, u(\cdot,t) \,\ni \!\;\!,\,
\|\,\cD^n u(\cdot,t) \,\|_{\mbox{}_{\scriptstyle L^{q}}} \!\to\!\:\!
\|\,\cD^n u(\cdot,t) \,\ni
} $
as $ \mbox{\small $q$} \to \!\;\!\mbox{\small $\infty$} $,
for all~\mbox{\small $n$}.$\mbox{}\!$\footnote{%
%
%
Moreover,
it is worth noticing that,
under the definitions above,
if an inequality of Gagliardo type
$ {\displaystyle
\|\,{\sf u}\,\|_{\mbox{}_{\scriptstyle L^{\!q}}}
\!\leq\!\;\!K\;\!
\|\,{\sf u}\,\|_{\mbox{}_{\scriptstyle L^{\!r_{\mbox{}_{\!\:\!1}}}}}^{1 - \theta}
\|\;\!\nabla {\sf u}\,
\|_{\mbox{}_{\scriptstyle L^{\!r_{\mbox{}_{\!\:\!2}}}}}^{\:\!\theta}
\!
} $,
$ 0 \leq \theta \leq 1 $,
holds for {\em scalar\/} functions $ {\sf u} $
(and some appropriate constant $K$),
then it will also be valid for
{\em vector\/} functions $ u $
with the {\em same\/} constant $K\!\:\!$
as in the scalar case.
Similarly,
one has
$ {\displaystyle
\;\!
\|\,\cD^n u(\cdot,t)\,\|_{\mbox{}_{\scriptstyle L^{\!q}}}
\!\:\!\leq\:\!
\|\,\cD^n u(\cdot,t)\,
\|_{\mbox{}_{\scriptstyle L^{\!q_{\mbox{}_{\!\:\!1}}}}}^{1 - \theta}
\|\,\cD^n u(\cdot,t)\,
\|_{\mbox{}_{\scriptstyle L^{\!q_{\mbox{}_{\!\:\!2}}}}}^{\:\!\theta}
\!\!\;\!
} $
if $ 1/q = (1 - \theta)/q_{\mbox{}_{1}} \!+ \theta/q_{\mbox{}_{2}} $,
$ 0 \leq \theta \leq 1 $,
and so on.
}$\!\;\!$) \linebreak
%
%
We may occasionally write $ \|\;u\;\|_{\mbox{}_{\scriptstyle L^q}}\!\;\!$
instead of $\|\,u(\cdot,t)\,\|_{\mbox{}_{\scriptstyle L^q}}\!\;\!$,
for simplicity.
Constants will be usually denoted by the letters
$ C \!\;\!$, $\!\;\!c$, $\!\;\!K$;
we write $ \:\!C_{\lambda}\!\;\! $
to indicate a constant whose value may depend on
a given parameter $\lambda$,
etc.
Also, for economy,
we often keep the same symbol for constants
in spite of possible changes in
their numerical values
(so, for example, we write $ C^2 \!\;\!$ again as $C\!\;\!$,
and so forth). \\
\mbox{} \vspace{-0.900cm} \\

An outline of the paper is as follows.
In Section \ref{s2},
we show that
a bound on the maximum norm
$\|\,u(\cdot,t)\,\ni$ in some interval
$0\leq t <T$ implies bounds for {\em all\/}
derivatives of $u(\cdot,t)$
in the same interval.
This is a well known result
that dates back to Leray \cite{Leray1934}
(see also \cite{KreissLorenz1989}),
but since it is the basis for all our blow--up results
we will prove it here.
An important implication is the following:
if
$\|\,u(\cdot,t)\,\ni$ is bounded in $0 \leq t < T$
for some finite $T$,
then $u(\cdot,t)$
can be continued as a $C^\infty$ solution beyond $T$.
This follows from well--known local constructions of solutions.
(See, for example,
\cite{KreissLorenz1989}
for an elementary proof.
See also \cite{Kato1984, KochTataru2001} for the
development of a local theory under much weaker assumptions
on the initial state $f$.)
In other words,
we can state the following first blow--up result:

%
%
%
%
\begin{theorem}
\label{t.1.max}
\label{theorem.u.Linfty}
If $ \,T_{\!\;\!f} < \infty $,
then \\
\mbox{} \vspace{-0.950cm} \\
\be
\label{sup.u.Linfty}
\sup_{0 \,\leq\, t \,<\, T_{\!\;\!f}}
\! \|\,u(\cdot,t)\,\ni
\;\!=\, \infty\:\!.\footnote{%
%
%
Actually,
$\;\!$if
$ \;\!T_{\!\;\!f} < \infty $
$\;\!$then one has
$ {\displaystyle
\lim_{t \nearrow \;\!T_{\!\:\!f}}
\!\;\!
\|\,u(\cdot,t)\,\|_{\mbox{}_{\infty}}
\!\:\!=\:\! \infty
} $,
%
%
$\;\!$cf.$\;$Theorem 1.3 below.
}
%
%
\ee
\end{theorem}
\mbox{} \vspace{-0.700cm} \\

In Sections \ref{s3} and \ref{s4} below,
we show that boundedness of
$\;\!\|\,u(\cdot,t) \,\|_{\mbox{}_{\scriptstyle L^q}}\!$
in some interval $0 \leq t < T$
for some $ q > 3 $,
or that of
$\;\!\|\,\cD u(\cdot,t) \,\|_{\mbox{}_{\scriptstyle L^q}}\!$
if $ \;\!3/2 < q \leq \infty $,
implies boundedness of $\;\!\|\,u(\cdot,t)\,\ni$
in the same time interval.
In particular, Theorem \ref{t.1.max}
yields the first part of the following result:

%
%
%
%
\begin{theorem}
\label{thm.1.2}
\label{theorem.Du.Lq}
$(i)$
Let $\,\f{3}{2} < q \leq \infty $.
If $ \,T_{\!\;\!f} < \infty $,
then \\
\mbox{} \vspace{-0.820cm} \\
\be
\label{sup.Du.Lq}
\sup_{0 \,\leq\, t \,<\, T_{\!\;\!f}} \!
\|\,\cD u(\cdot,t)\,\|_{\mbox{}_{\scriptstyle L^q}}
\,\!=\, \infty\:\!.
\ee
\mbox{} \vspace{-0.200cm} \\
$(ii)$
For each $ \,\f{3}{2} \leq q < 3 $,
there exists an absolute constant $\;\!c_{q}\!\,\! > 0 $,
independent of $\;\!t\!\;\!$ and $\!\;\!f\!\,\!$, \linebreak
with the following property$\;\!:$
if $ \,T_{\!\;\!f} < \infty $,
then \\
\mbox{} \vspace{-0.970cm} \\
\be
\label{lowerbound.Du.Lq}
\|\,\cD u(\cdot,t)\,\|_{\mbox{}_{\scriptstyle L^{q}}}
\:\!\geq\,
c_{q} \!\:\! \cdot (\:\!T_{\!\;\!f} -\;\! t\:\!
)^{\mbox{}^{\scriptstyle \!\!-\, \frac{\scriptstyle \;\!q \:-\: 3/2\,}
                                      {\scriptstyle q} }}
\quad \;\;
\forall \;\,
0 \leq t < T_{\!\;\!f}.
\ee
\end{theorem}
%
%
\mbox{} \vspace{-0.310cm} \\
Regarding (\ref{sup.Du.Lq}),
it will be shown
in Section \ref{s4} that
one actually has,
for each
$\;\! \f{3}{2} < q \leq \infty $,
the stronger property \\
\mbox{} \vspace{-0.950cm} \\
\be
\label{limit.Du.Lq}
\lim_{t \:\!\nearrow\;\! T_{\!\;\!f}}
\|\,\cD u(\cdot,t)\,\|_{\mbox{}_{\scriptstyle L^q}}
\:\!=\, \infty
\qquad
\mbox{($\,\!$if $\,T_{\!\;\!f\!\;\!} < \infty$)}.
\ee
\mbox{} \vspace{-0.850cm} \\

For $ q = 2 $,
Theorem \ref{thm.1.2}
was originally proved in \cite{Leray1934}.
%
%
The estimate (\ref{lowerbound.Du.Lq})
is only one of {\em many\/} similar
lower bound results for blow-up rates
of solution--size quantities $Q(u(\cdot,t)) $
that have been obtained since Leray \cite{Leray1934}.
In general,
these bounds result from some form of
local control on $Q(u(\cdot,t)) $
that is typically obtained
in one of the following basic ways: \\
\mbox{} \vspace{-0.200cm} \\
($\,\!$\mbox{\small I}$\,\!$) lower bounds for
the maximum existence time $T\!\,\!$ of solution $ u(\cdot,t) $
given in terms of~$ Q $; \\
\mbox{} \vspace{-0.300cm} \\
(\mbox{\small II}) differential or integral inequalities
satisfied by $ Q(u(\cdot,t)) $ while the solution exists; \\
\mbox{} \vspace{-0.300cm} \\
(\mbox{\small III}) relationships with
other functionals $ \tilde{Q}(u) $,
like Sobolev or interpolation inequalities. \\
\mbox{} \vspace{-0.100cm} \\
As to method (I),
we observe that
lower bound estimates for $T$
are a usual, important by--product of construction schemes
in existence theory,
so that blow--up estimates like (\ref{lowerbound.Du.Lq})
are actually very natural
and widespread
in the literature.
Thus,
for example,
for the \linebreak
solutions
$ u(\cdot,t) $ of (1.1),
it can be shown
(see e.g.$\;$\cite{Galdi2000, Heywood1980, Leray1934})
that \\
\mbox{} \vspace{-0.700cm} \\
\be
\label{existence.time.T}
\label{existence.time.T.Df}
T \;\!\geq\,
C \: \|\, \cD f \,\|^{- 4}
\ee
\mbox{} \vspace{-0.350cm} \\
for some absolute constant $ C > 0 $
independent of $f$.
It follows that,
given $ 0 \leq t_0 \!\:\!<\!\:\! T_{\!\;\!f} $
arbitrary,
we must have
$ \;\!T_{\!\;\!f} \!\;\!-\;\! t_0 \:\!>\;\!
C \: \|\, \cD u(\cdot,t_0) \,\|^{-4} \!\;\!$,
which is (\ref{lowerbound.Du.Lq}) in case $ q = 2 $.
As noted \linebreak
in \cite{RobinsonSadowskiSilva2012},
this is a two--way route:
had we had the estimate
(\ref{lowerbound.Du.Lq}) for $ q = 2 $
in the first place,
then we could have gotten
(\ref{existence.time.T})
from it
just as easily
by merely taking $ \,\!t = 0 $
there. \\
\mbox{} \vspace{-1.000cm} \\

In this review,
as we will be bypassing existence theory,$\!\;\!$\footnote{%
%
%
Some standard references for existence results are
\cite{ConstantinFoias1988, Galdi2000, %
Heywood1980, Kato1984, KreissLorenz1989, Ladyzhenskaya1969, %
Leray1934, Masuda1984, Serrin1963, Sohr2001, Temam1977, Weissler1981}.
}
%
%
our approach to obtaining
the blow--up estimates
considered here
will be based on the methods
(\mbox{\small II})
and (\mbox{\small III}) alone.\footnote{%
%
%
The way to obtain lower bound estimates for blow--up rates
out of nonlinear differential or integral inequalities
is very simple and is shown
in Lemmas \ref{How.differential.inequalities.are.used}
and \ref{How.integral.inequalities.are.used} below,
respectively.
}
%
%
%
In particular,
our derivation of (\ref{lowerbound.Du.Lq})
for $ q = 2 $
in Section~3
uses a well known differential inequality
satisfied by the function
$ \|\,\cD u(\cdot,t)\,\|^2\!$,
which will give us that
$ \;\!c_{\mbox{}_{2}}^{\,2} \!\;\!\geq\:\! 2 \;\!\pi\:\!\sqrt{\:\!2\;\!} $
\cite{Heywood1994}.
%
%
One should note that the blow--up rate of
$\;\!\|\,\cD u(\cdot,t)\,\|$,
as expressed by the estimate~(\ref{lowerbound.Du.Lq}),
is clearly consistent with the
fundamental upper bound \\
\mbox{} \vspace{-0.650cm} \\
\be
\label{eq.1.4a}
\label{integral.Du.L2}
\int_0^{T_{\!\;\!f}} \!\!\;\!
\|\;\!\cD u(\cdot,t)\,\|^2\, dt
\,\leq\:
\mbox{\small $ {\displaystyle \f{1}{2} }$} \: \|\,f\,\|^2
\ee
\mbox{} \vspace{-0.200cm} \\
that follows from
the energy equality satisfied by $ u(\cdot,t) $,
which we recall in Theorem \ref{t.2.energy}.
%
It then follows from
(\ref{lowerbound.Du.Lq}) and (\ref{eq.1.4a})
that,
if $ T_{\!\;\!f}\!\;\! < \infty $, \\
\mbox{} \vspace{-0.570cm} \\
$$ 2 \,c_{\mbox{}_{2}}^{\;\!2} \, T_{\!\;\!f}^{\;\!1/2}
\:\!=\;
c_{\mbox{}_{2}}^{\;\!2} \!\:\! \int_0^{T_{\!\;\!f}} \!
(\;\!T_{\!\;\!f} -\;\! t\:\!)^{-1/2}
\,dt
\:\leq
\int_0^{T_{\!\;\!f}} \!\!\;\!
\|\;\!\cD u(\cdot,t)\,\|^2\, dt
\,\leq\:
\mbox{\small $ {\displaystyle \f{1}{2} }$} \: \|\,f\,\|^2
\!\;\!,
$$ \\
\mbox{} \vspace{-0.800cm} \\
so that
we must have \\
\mbox{} \vspace{-0.620cm} \\
\be
\label{necessary.condition.1a}
c_{\mbox{}_{2}}^{\;\!4} \: \|\, \cD f \,\|^{-\;\!4}
\,\leq\:
T_{\!\;\!f} \,\leq\;
\mbox{\small $ {\displaystyle \frac{1}{\:\!16\,c_{\mbox{}_{2}}^{\;\!4} } }$} \: \|\,f\,\|^4
\qquad
\mbox{($\,\!$if $\,T_{\!\;\!f\!\;\!} < \infty$)},
\ee
\mbox{} \vspace{-0.220cm} \\
where, again, the first inequality
follows
from (\ref{lowerbound.Du.Lq})
by taking $ t = 0 $.
$\!$Therefore,
recalling that
$ c_{\mbox{}_{2}}^{\,2} \geq\:\! 2 \;\!\pi\:\!\sqrt{\:\!2\,}$,
finite--time
blow--up at
$ T_{\!\;\!f} $
is only possible if
we have \\
\mbox{} \vspace{-0.700cm} \\
\be
\label{necessary.condition.1}
\|\, u(\cdot,t) \,\| \;\:\!
\|\, \cD u(\cdot,t) \,\|
\:\geq\:
4\;\!\pi \:\!\sqrt{\:\!2\,}
\qquad
\forall \;\, 0 \leq t < T_{\!\;\!f}
\qquad
\mbox{($\,\!$if $\,T_{\!\;\!f\!\;\!} < \infty$)}.
\ee
\mbox{} \vspace{-0.240cm} \\
This is only one
of many {\em necessary\/} conditions
for finite--time blow--up
that have been found
since the fundamental paper
of Leray \cite{Leray1934}
(see (\ref{necessary.condition.Lq}) and (\ref{eq.17.s4})
for other similar examples),
the most celebrated of them
being the Beale-Kato-Majda condition
\cite{BealeKatoMajda1984} \\
\mbox{} \vspace{-0.570cm} \\
\be
\label{BealeKatoMajda1984}
\int_0^{T_{\!\;\!f}} \!\!\;\!
\|\, \nabla \!\times\!\;\! u(\cdot,t) \,
\|_{\mbox{}_{\scriptstyle L^{\infty}}} \,dt
\;=\; \infty
\qquad
\mbox{($\,\!$if $\,T_{\!\;\!f\!\;\!} < \infty$)},
\ee
\mbox{} \vspace{-0.220cm} \\
originally obtained for the
Euler equations,
but which also holds for the Navier--Stokes equations.
Actually,
the derivation of (\ref{BealeKatoMajda1984})
for (1.1)
is much easier,
as shown in Section~6.

%
\mbox{} \vspace{-0.975cm} \\

The proof of Theorem 1.2
is completed by method (\mbox{\small III})
after we consider,
in Section~\ref{s4},
the important solution norms \\
\mbox{} \vspace{-0.900cm} \\
$$\|\,u(\cdot,t)\,\|_{\mbox{}_{\scriptstyle L^q}}\!\;\!,
\qquad
3 < q \leq \infty. $$
\mbox{} \vspace{-0.420cm} \\
If any of these norms stays bounded
in some interval $0 \leq t < T \!\;\!$,
we show that
$ \|\, \cD u(\cdot,t)\,\| $
will also be bounded in that interval.
As this implies boundedness for
$ \|\, u(\cdot,t) \,\ni $
as well,
it follows that \\
\mbox{} \vspace{-1.100cm} \\
\be
\label{sup.Lq.norm}
\label{sup.u.Lq}
\sup_{0 \,\leq\, t \,<\, T_{\!\;\!f}} \!
\|\,u(\cdot,t) \,\|_{\mbox{}_{\scriptstyle L^q}}
=\, \infty,
\quad \;\;
3 < q \leq \infty
\ee
\mbox{} \vspace{-0.230cm} \\
if $ \;\!T_{\!\;\!f} < \infty $.
In fact,
one can again derive an algebraic lower bound
for the blow--up rate:

%
%
%
%
\begin{theorem}
\label{thm.1.3}
\label{theorem.rates.u.Lq}
For each $\;\!3 \leq q \leq \infty$,
there is a constant $c_q > 0$,
independent of $\;\!t$ and $f\!\;\!$,
such that
the following holds:
if $ \,T_{\!\;\!f} < \infty$,
then \\
\mbox{} \vspace{-0.170cm} \\
\mbox{} \hspace{+3.500cm}
$ {\displaystyle
\|\,u(\cdot,t)\,\|_{\mbox{}_{\scriptstyle L^q}}
\:\!\geq\,
c_q \!\:\!\cdot
(\;\!T_{\!\;\!f} - \;\!t\:\!)^{-\kappa}
\quad \;\;\,
\forall \;\,
0 \leq t < T_{\!\;\!f},
} $
\mbox{} \hfill $(1.12a)$ \\
\mbox{} \vspace{-0.450cm} \\
with \\
\mbox{} \vspace{-0.500cm} \\
\mbox{} \hspace{+3.000cm}
$ {\displaystyle
\kappa \,=\: \f{\:\!q - 3\:\!}{2\:\!q}
\quad \mbox{{\em if\/} $\,3 \leq q < \infty $},
\quad \;\;\,
\kappa \,=\: \mbox{\small $ {\displaystyle \f{1}{\:\!2\:\!} }$}
\quad \mbox{{\em if\/} $\,q \;\!=\;\! \infty $}.
} $
\mbox{} \hfill $(1.12b)$ \\
\setcounter{equation}{12}
\end{theorem}
%
%
\mbox{} \vspace{-0.700cm} \\
In particular,
from (1.12),
we have
$ {\displaystyle
\lim_{t \nearrow\;\!T_{\!\;\!f}}
\!\:\!
\|\, u(\cdot,t) \,\|_{\mbox{}_{\scriptstyle L^{q}}}
\!\,\!=\!\;\! \infty \:\!
} $
if
$ T_{\!\;\!f}\!\;\!< \infty $,
for each
$ 3 < q \leq \infty $. \\
\mbox{} \vspace{-0.150cm} \\
%
%
%
%
\noindent
{\sc Remark:}
The property (\ref{sup.Lq.norm})
is also valid for the limit case $ q = 3 $,
as shown in \cite{EsauriazaSereginSverk2003}.
The proof, however, is very involved
and will not be covered here.
More recently,
it has been shown by G. Seregin \cite{Seregin2012}
the stronger result \\
\mbox{} \vspace{-0.750cm} \\
\be
\label{limit.u.L3}
\lim_{t \nearrow\;\!T_{\!\;\!f}}
\!\;\!
\|\,u(\cdot,t)\,\|_{\mbox{}_{\scriptstyle L^3}}
\!\;\!=\,\infty
\qquad
\mbox{($\:\!$if $\,T_{\!\;\!f\!\;\!} < \infty $)}.
\ee
\mbox{} \vspace{-0.220cm} \\
%
%
It then follows from
(\ref{limit.u.L3})
and the 3D Sobolev inequality
$ {\displaystyle
\|\: {\sf u} \:\|_{\mbox{}_{\scriptstyle L^{3}}}
\!\;\!\leq \:\!
K \:
\|\: \nabla {\sf u} \:\|_{\mbox{}_{\scriptstyle L^{3/2}}}
\!\;\!
} $
that the properties
(\ref{sup.Du.Lq}) and (\ref{limit.Du.Lq}) above
are both valid for $ q = 3/2 $ as well.

%
%

\mbox{} \vspace{-0.750cm} \\
The estimates (1.12)
were originally given in \cite{Leray1934}
and reobtained in a more general setting
using semigroup ideas in \cite{Giga1986}.
They immediately imply
lower bounds for the
existence time of $ u(\cdot,t) $
of the form \\
\mbox{} \vspace{-0.930cm} \\
\be
\label{existence.time.Lq}
T \;\!\geq\: C_{\!\;\!q} \,
\|\, f \,\|_{\mbox{}_{\scriptstyle L^{q}}}
           ^{\mbox{}^{\scriptstyle \!- \,
\frac{\scriptstyle 2\:\!q}{\scriptstyle \:\!q \:\!-\:\! 3\:\!} }}
\!,
\quad \;\;\;\;
3 < q \leq \infty,
\ee
\mbox{} \vspace{-0.520cm} \\
where
$ {\displaystyle
\:\!C_{\!\;\!q} =\:\! c_q^{\mbox{}^{\scriptstyle
\frac{\scriptstyle 2\:\!q}{\scriptstyle \:\!q - 3\:\!}}}
\!\!\;\!
} $.
The estimates
(\ref{existence.time.Lq})
are obtained directly
from existence analysis
in \cite{Leray1934} (for $q = \infty$)
and in \cite{RobinsonSadowski2012}
(for $ 3 < q < \infty $),
thus providing another proof for
(1.12).
%
%
%
Again,
our derivation of (1.12),
which is carried out in Section \ref{s4},
follows the method (\mbox{\small II}),
firstly along the lines of \cite{Giga1986}
using some well established integral inequalities
satisfied by the quantities
$ \|\, u(\cdot,t) \,\|_{\mbox{}_{\scriptstyle L^{q}}} \!\,\!$
to obtain the result,
and then
by deriving some less known differential inequalities
that can be used for this purpose
just as easily.
%
%
We also obtain
from the latter analysis
another proof of the following
result
(see e.g.$\;$\cite{Galdi2000, Giga1986,%
Leray1934, RobinsonSadowski2012})
on the global existence of smooth solutions
for the Navier--Stokes problem (1.1).

%
%
%
%
\begin{theorem}
\label{thm.1.4}
\label{theorem.global.solvability}
For each $\:\! 3 \leq q \leq \infty $,
there exists a number $ \:\!\eta_{q} \!\:\!> 0 $,
depending only on~$\,\!q$,
such that
%
\mbox{} \vspace{-0.620cm} \\
\be
\label{global.existence.Lq}
\|\: f \:\|_{\mbox{}_{\scriptstyle \!\:\!L^{2}}}
           ^{\mbox{}^{\scriptstyle \frac{\scriptstyle \;\!2 q \;\!-\;\!6\;\!}
                                        {\scriptstyle \;\!3 q \;\!-\;\!6\;\!} }}
\:\!
\|\: f \:\|_{\mbox{}_{\mbox{}_{\scriptstyle \!\:\!L^{q}}}}
           ^{\mbox{}^{\scriptstyle \frac{\scriptstyle q}
                                        {\scriptstyle \;\!3 q \;\!-\;\!6\;\!} }}
<\: \eta_{q}
\;\;\;
\Longrightarrow
\;\;\,
T_{\!\;\!f} =\, \infty.
\ee
\end{theorem}
%
%
\mbox{} \vspace{-0.270cm} \\
In particular,
finite--time blow--up
of a smooth solution $ u(\cdot,t) $
can only occur if
we have \\
\mbox{} \vspace{-0.720cm} \\
\be
\label{necessary.condition.Lq}
\|\, u(\cdot,t) \,\|_{\mbox{}_{\scriptstyle \!\:\!L^{2}}}
                    ^{\mbox{}^{\scriptstyle
                        \frac{\scriptstyle \;\!2 q \;\!-\;\!6\;\!}
                             {\scriptstyle \;\!3 q \;\!-\;\!6\;\!} }}
\:\!
\|\, u(\cdot,t) \,\|_{\mbox{}_{\mbox{}_{\scriptstyle \!\:\!L^{q}}}}
                    ^{\mbox{}^{\scriptstyle
                        \frac{\scriptstyle q}
                             {\scriptstyle \;\!3 q \;\!-\;\!6\;\!} }}
\geq\: \eta_{q}
\quad \;\;\;
\forall \;\, 0 \leq t < T_{\!\;\!f}
\quad
\;\;\;(\:\!\mbox{if $\;\!T_{\!\;\!f}\!\:\!< \infty $}) \!\!\!\!\!
\ee
\mbox{} \vspace{-0.240cm} \\
for every $ \:\!3 \leq q \leq \infty $,
where $ \eta_{q} > 0 $
is the value given in
(\ref{global.existence.Lq})
above.
If this is the case,
using the 3D Sobolev inequality \\
\mbox{} \vspace{-0.720cm} \\
\be
\label{SNG.Du.Lq}
\|\: \mbox{\sf u} \:\|_{\mbox{}_{\scriptstyle L^{r(q)}}}
\,\!\leq\,
K_{\!\:\!q} \,
\|\, \nabla \!\;\!\mbox{\sf u} \:\|_{\mbox{}_{\scriptstyle L^{q}}},
\quad
\;\;
r(q) \;\!=\;\! \frac{3 \:\!q}{\;\!3 - q\,},
\quad
\;
\mbox{\normalsize $ {\displaystyle \frac{\;\!3\;\!}{2} }$} \leq q < 3,
\ee
\mbox{} \vspace{-0.270cm} \\
we obtain,
from (1.12),
the blow--up estimate
(\ref{lowerbound.Du.Lq})
in Theorem 1.2 above.
This illustrates
the use of method (\mbox{\small III})
to derive these results.
Other examples are found
in Sections 4, 5 and 6 below,
including the following general blow--up property
for arbitrary high order derivatives
of smooth solutions of (1.1).

%
%
%
%
\begin{theorem}
\label{thm.1.5}
\label{theorem.blow.up.Dnu}
Let $\,n \geq 2 \;\!$
be an integer.
If $\,T_{\!\;\!f} \!\;\!< \infty $,
we have \\
\mbox{} \vspace{-0.750cm} \\
\be
\label{limit.Dn.u.Lq}
\lim_{t \nearrow \;\!T_{\!\;\!f}} \;\!
\|\, \cD^{n} u(\cdot,t) \,\|_{\mbox{}_{\scriptstyle L^{q}}}
=\, \infty
\ee
\mbox{} \vspace{-0.250cm} \\
for every $\;\! 1 \leq q \leq \infty $.
\end{theorem}
%
%

It is also worth noticing here that,
from (1.12),
we clearly have \\
\mbox{} \vspace{-0.570cm} \\
$$
\int_0^{T_{\!\;\!f}} \!\!\;\!
\|\, u(\cdot,t) \,\|_{\mbox{}_{\scriptstyle L^{q}}}
                    ^{\mbox{}^{\scriptstyle \:\!r}}
\;\!dt
\;=\; \infty
\qquad
\mbox{($\:\!$if $\,T_{\!\;\!f\!\;\!} < \infty $)}
$$
\mbox{} \vspace{-0.200cm} \\
for any $ r \geq 2\:\!q/(q - 3) $,
or, equivalently,
for any $ r \geq 2 $
satisfying
$ \:\!2/r \:\!+\:\! 3/q \leq 1 $.
It is therefore natural
to expect that the so-called
{\em Prodi-Serrin condition}, \\
\mbox{} \vspace{-0.570cm} \\
\be
\label{Prodi.Serrin.Struwe}
\int_0^{T} \!\!\;\!
\|\, u(\cdot,t) \,\|_{\mbox{}_{\scriptstyle L^{q}}}
                    ^{\mbox{}^{\scriptstyle \:\!r}}
\;\!dt
\;<\; \infty,
\qquad
\frac{2}{\:\!r\:\!} \,+\, \frac{3}{\:\!q\:\!} \,\leq\, 1,
\ee
\mbox{} \vspace{-0.200cm} \\
for some
$ 2 \leq r < \infty $, $ 3 < q \leq \infty $ (arbitrary),
imposed on less regular {\em weak\/} solutions,
may be sufficient to guarantee strong regularity
and uniqueness properties.
This is indeed the case,
as shown
in \cite{Masuda1984, Prodi1959, Serrin1962, Struwe1988}
(see also \cite{Galdi2000, Giga1986, Kato1984, %
KozonoSohr1996, MontgomerySmith2005, Sohr2001}),
but it requires a more advanced analysis
and will not be discussed here.
Similar observations apply to
the other blow-up quantities
considered in (\ref{limit.Du.Lq}),
(\ref{BealeKatoMajda1984}) or (\ref{limit.u.L3}),
see e.g.$\;$\cite{EsauriazaSereginSverk2003, %
Galdi2000, Leray1934, RobinsonSadowski2012, Seregin2012}. \\
\mbox{} \vspace{-0.950cm} \\

If there is blow--up, one can ask:
do certain norms blow--up faster than others?
The answer is yes.
For example, if $\;\!3 \leq q < r \leq \infty $,
we show in Section \ref{s5}
that the $L^r\!\;\!$ norm blows up faster
than the $L^q\!\;\!$ norm,
with \\
\mbox{} \vspace{-0.600cm} \\
\be
\label{faster.blow.up.u.Lq.Lr}
\frac{\,\|\,u(\cdot,t)\,\|_{\mbox{}_{\scriptstyle L^{r}\!}}}
     {\,\|\,u(\cdot,t)\,\|_{\mbox{}_{\scriptstyle L^{q}\!}}}
\:\geq\;
c(q,r) \cdot
\|\, f \:\|^{\mbox{}^{\scriptstyle \!\;\!- \;\!\lambda}}
\hspace{-0.400cm} \cdot \hspace{+0.070cm}
(\:\!T_{\!\;\!f} -\:\!t\:\!)^{\mbox{}^{\scriptstyle \!-\;\!\gamma}} \!\!,
\quad
\;\;\,
\gamma \,=\, \frac{\;\!r - 3\;\!}{\;\!r - 2\;\!} \,
\biggl(\:\! \frac{1}{\;\!q\;\!} \;\!-\;\! \frac{1}{\;\!r\;\!} \:\!\biggr)
\ee
\mbox{} \vspace{-0.035cm} \\
for all $ \;\!0 \leq t < T_{\!\;\!f}\!\;\! $,
where $ c(q,r) > 0 $ depends only on $ q, \:\!r $,
and $ \lambda = 2 \;\!(r/q - 1)/(r - 2) $. \linebreak
These relations are typically obtained
by the approach \mbox{\small (III)}.
A further result of this kind,
which is related to (\ref{limit.u.L3}) above,
is also included in Section \ref{s5},
and given a direct proof that is
independent of (\ref{limit.u.L3}). \\
\mbox{} \vspace{-0.950cm} \\

In Section 6,
we briefly examine some related properties
for the flow vorticity
$ \omega(\cdot,t) = \nabla \!\:\!\times \!\;\! u(\cdot,t) $.
Our main goal in this Section
is to provide a short and very simple proof
for the Beale-Kato-Majda blow--up condition
(\ref{BealeKatoMajda1984})
for smooth solutions of the
Navier--Stokes equations.
This particular proof is not valid for
the inviscid Euler equations. \\
\mbox{} \vspace{-0.950cm} \\

Besides the famous major problems,
there are still many other open questions
related to our discussion
and some are indicated in the text.
For additional lower bound estimates and results
concerning other blow--up quantities,
the reader is referred to
\cite{Bae2012, Benameur2010, ConstantinFoias1988, %
FoiasGuillopeTemam1981, RobinsonSadowskiSilva2012}.

\mbox{} \vspace{-0.050cm} \\
{\bf Acknowledgments.}
This paper is based in part
on lectures given by the first author
during a recent visit to the
Universidade Federal do Rio Grande do Sul,
in Porto Alegre.
We express our gratitude
to CNPq for the financial support
provided for that visit.

\end{section}

\mbox{} \vspace{-1.750cm} \\
%
%
%
\begin{section}{A bound for $\|\,u(\cdot,t)\,\|_\infty$
implies bounds for all \\derivatives}
\label{s2}
\mbox{} \vspace{-0.900cm} \\
\setcounter{equation}{0}

Under our assumptions on $f$
stated in the introduction,
the Cauchy problem (\ref{eq.1.1})
has a unique $C^\infty$ solution
$\bigl( u(x,t), p(x,t) \bigr)$,
defined in some interval
$0 \leq t < T_{\!\;\!f}\!\:\!$,
with
$D^{\mbox{}^{\scriptstyle \alpha}} \!\:\!u(\cdot,t) \in L^2(\R^3)$
for all multi-index $\:\!\alpha\:\!$
and all $\;\! 0 \leq t < T_{\!\;\!f} \!\:\!$.
(Also, recall that we always require
$\:\!p(x,t) \to 0\;\!$ as $\;\!|\,x\,| \to \infty \:\!$
to make $\:\!p(\cdot,t)$ unique.)

As in \cite{Leray1934},
we set \\
%
%
\mbox{} \vspace{-0.900cm} \\
$$
J_n^2(t) \,=\!\!\!\;\!
\sum_{\;\,|\,\al\,|\,=\,n} \!\!\!\!\;\!
\|\;\!D^\alpha u(\cdot,t)\,\|^2
\!\;\!,
\quad \;\;\;
n \;\!=\;\! 0, 1, 2, \ldots
$$
\mbox{} \vspace{-0.150cm} \\
The most basic estimate for the solution of the Navier--Stokes
equations is the following well known energy estimate.
The proof follows from multiplying the
equation for $ u_{i}(\cdot,t) $
by $ u_{i}(\cdot,t) $
and integrating by parts
(see e.g.$\;$\cite{Galdi2000, KreissLorenz1989, %
Ladyzhenskaya1969, Serrin1963, Sohr2001, Temam1977}).

%
%
%
%
\begin{theorem}
\label{t.2.energy}
\label{theorem.energy.equation}
We have \\
\mbox{} \vspace{-0.225cm} \\
\mbox{} \hspace{+4.000cm}
$ {\displaystyle
\f{1}{2}\;\,\f{d}{dt}\,J_0^2(t)
\,=\, -\,J_1^2(t),
\qquad
\forall \;\,0 \leq t < T_{\!\;\!f}\!\;\!,
} $
\hfill $(2.1a)$ \\
\mbox{} \vspace{-0.270cm} \\
so that,
in particular, \\
\mbox{} \vspace{-0.200cm} \\
\mbox{} \hspace{+1.750cm}
$ {\displaystyle
J_0(t) \;\!\leq\, \|\,f\,\| \qfor 0 \leq t < T_{\!\;\!f}
\qand
\int_0^{\mbox{\footnotesize $t$}} \!\!\;\! J_1^2(s)\, ds
\;\leq\; \f{1}{2}\; \|\,f\,\|^2.
} $
\hfill $(2.1b)$ \\
\end{theorem}
\setcounter{equation}{1}
%
%
%
%
%
%
\mbox{} \vspace{-1.250cm} \\

Note that the integral bound in
$(2.1b)$ proves (\ref{eq.1.4a}).
To prove the next result,
we will use the 3D Sobolev inequality \\
\mbox{} \vspace{-0.800cm} \\
\be
\label{Sobolev.embedding}
\|\;\mbox{v}\;\ni
\;\!\leq\;\!
C \,\|\;\mbox{v}\;\|_{\mbox{}_{\scriptstyle H^2}}
\qquad
\forall \;\,
\mbox{v} \in H^{2}(\R^3),
\ee
which implies \\
\mbox{} \vspace{-0.800cm} \\
\be
\label{eq.3.s2}
\|\;\!D^j u(\cdot,t)\,\ni \;\!\leq\,
C_{j} \!\;\!\cdot\!\;\!
\bigl(\,\!J_{m}(t) + J_0(t)\:\!\bigr)
\qfor m \geq j + 2\:\!.
\ee
\mbox{} \vspace{-0.270cm} \\
(The bound $J_n^2(t) \:\!\leq\;\! C_{n}\!\!\;\!\cdot\!\;\!
\bigl(J_m^2(t) + J_0^2(t)\:\!\bigr)\;\!$
for $\;\!m > n\;\!$
follows by Fourier transformation.)

%
%
%
%
\begin{theorem}
\label{theorem.supnorm.Leray}
Assume that \\
\mbox{} \vspace{-1.000cm} \\
\be
\label{eq.2.M}
\label{bound.u.Linfty}
\sup_{0 \,\leq\, t \,<\, T} \!
\|\,u(\cdot,t)\,\ni
=:M < \infty.
\ee
\mbox{} \vspace{-0.250cm} \\
Then each function $J_n(t)$, $ n = 1, 2, \ldots$,
is bounded in $\;\!0 \leq t < T$
by some quantity $\:\!K_{\!\;\!n}\!> 0 $
depending only
on $\;\!n$, $M$, $T\!\;\!$ and
$\;\!\|\, f \,\|_{\mbox{}_{\scriptstyle H^{n}}}\!\:\!$,
that is,
$ K_{\!\;\!n} =
K\!\;\!(n, M, T\!\;\!, \|\,f\,\|_{\mbox{}_{\scriptstyle H^{n}}}) $.
\end{theorem}
%
%
{\bf Proof:}
Using (\ref{eq.2.M}),
we will first prove that $J_{1}(t), \;\! J_2(t)$
are bounded
in $\:\!0 \leq t < T \!\;\!$,
and then make an induction argument in $n$.
%
We have, for any multi--index $\al$, \\
\mbox{} \vspace{-0.670cm} \\
$$D^\al u_t \:\!+\;\! D^\al(u\cdot \nabla u)
\;\!+\:\! \nabla \;\!D^\al p
\,=\, \Delta \;\! D^\al u $$
\mbox{} \vspace{-0.420cm} \\
and,
because $ \nabla \!\cdot u = 0 $, \\
\mbox{} \vspace{-0.970cm} \\
\beaa
\f{1}{2}\;\f{d}{dt}\,J_n^2(t) &=&
\!\!\sum_{\;|\,\al\,|\,=\,n} \!\!(D^\al u, D^\al u_t) \\
 &\leq& \!\!- \!\!
     \sum_{\;|\,\al\,|\,=\,n} \!\!(D^\al u, D^\al (u\cdot \nabla u)\:\!)
         \,-\, J_{n+1}^2(t) \\
 &=:& S_n(t) \;\!-\;\! J_{n+1}^2(t).
\eeaa
It is convenient to use the short notation \\
\mbox{} \vspace{-0.750cm} \\
$$(D^i u, D^j u \;\!D^k u) $$
for any integral \\
\mbox{} \vspace{-1.000cm} \\
$$\int_{\mbox{}_{\mbox{\scriptsize $\!\;\!\R^3$}}} \!\!\!\!\;\!
D^\al u_{\nu_{\mbox{}_{1}}} D^\beta u_{\nu_{\mbox{}_{2}}}
D^\gamma u_{\nu_{\mbox{}_{3}}} \, dx $$
\mbox{} \vspace{-0.250cm} \\
with $\nu_1, \nu_2, \nu_3 \in \{\:\!1,2,3\;\!\}$
and
$|\,\al\,| = i, \,|\,\beta\,| = j, \,|\,\gamma\,| = k $. \\
\mbox{} \vspace{-0.750cm} \\

Let $n = 2$
and consider the terms appearing in $S_2$, \\
\mbox{} \vspace{-0.750cm} \\
$$  (D^\al u, D^\al (u\cdot \nabla u)), \qquad |\,\alpha\,| = 2\:\!.$$
\mbox{} \vspace{-0.350cm} \\
Thus $S_2$ is a sum of terms of the form \\
\mbox{} \vspace{-0.700cm} \\
$$(D^2 u, \,u\;\!D^3 u) \qand (D^2 u, \:\!Du\;\!D^2 u)\:\!.$$
\mbox{} \vspace{-0.450cm} \\
Using integration by parts,
a term $(D^2 u,\;\!Du \;\!D^2 u)$
can also be written as
a sum of terms $(D^2 u, \,u\;\!D^3 u)$.
Since, by (\ref{eq.2.M}), \\
\mbox{} \vspace{-0.750cm} \\
$$|\,(D^2 u,\,u \;\!D^3 u)\,| \,\leq\;\! M J_2(t) \:\! J_3(t), $$
we obtain that \\
\mbox{} \vspace{-1.150cm} \\
\beaa
\f{d}{dt}\,J_2^2(t) & \leq & C \:\! M J_2(t) \:\! J_3(t) \,-\, 2 \:\!J_3^2(t) \\
 & \leq & C^2 M^2 J_2^2(t)
\eeaa
for some constant $ C > 0 $.
Boundedness of $J_2(t) $ in $0 \leq t < T$ follows.
Similarly,
we get \\
\mbox{} \vspace{-0.950cm} \\
\beaa
\f{d}{dt}\,J_1^2(t) & \leq & C \:\! M J_1(t) \:\! J_2(t) \,-\, 2 \:\!J_2^2(t) \\
 & \leq & C^2 M^2 J_1^2(t),
\eeaa
so that the result is true
for $ J_{1}(t) $ as well.

Now,
let $n \geq 2$
and assume $J_n(t)$ to be bounded in $ 0 \leq t < T$.
We have, from (\ref{Navier.Stokes.eqs}), \\
\mbox{} \vspace{-0.570cm} \\
$$\f{1}{2}\;\f{d}{dt}\,J_{n+1}^2(t) \,\leq\,
S_{n+1}(t) \,-\;\! J^2_{n+2}(t), $$
\mbox{} \vspace{-0.130cm} \\
where $S_{n+1}(t)$ is a sum of terms \\
\mbox{} \vspace{-0.650cm} \\
$$T_{j}(t) \,=\,
(D^{n+2}u, \:\!D^j u\;\!D^{n+1-j}u\:\!),
\qquad 0 \leq j \leq n\:\!.$$
\mbox{} \vspace{-0.250cm} \\
There are three cases to consider: \\
\mbox{} \vspace{-0.250cm} \\
({\em i\/}) Let $0\leq j \leq n-2$.
We have,
by (\ref{eq.3.s2}), \\
\mbox{} \vspace{-1.050cm} \\
\beaa
|\,T_{j}(t)\,| & \leq &
C_n \;\! \|\;\!D^j u(\cdot,t) \,\ni \;\! J_{n+2}(t)\,(J_{n+1}(t) + J_0(t)\:\!) \\
 & \leq & C_n \;\!J_{n+2}(t) \:(J_{n+1}(t) + J_0(t)\:\!).
\eeaa
In the latter estimate we have used (\ref{eq.3.s2})
and the induction hypothesis. \\
\mbox{} \vspace{-0.250cm} \\
({\em ii\/}) Let $j = n-1$.
We have,
by (\ref{eq.3.s2}), \\
\mbox{} \vspace{-1.050cm} \\
\beaa
|\,T_{n-1}(t)\,| & \leq & C_n \,\|\;\!D^{n-1}u(\cdot,t)\,\ni
\;\! J_{n+2}(t) \,J_{2}(t) \\
 & \leq & C_n \,\|\;\!D^{n-1}u (\cdot,t)\,\ni \;\!J_{n+2}(t) \\
 & \leq & C_n \,(J_{n+1}(t) + J_0(t)\,\!) \,J_{n+2}(t).
\eeaa
In the second estimate we have used that
a bound for $J_2(t) $ is already shown. \\
\mbox{} \vspace{-0.250cm} \\
({\em iii\/}) Let $j = n$.
We have,
by (\ref{eq.3.s2}), \\
\mbox{} \vspace{-1.050cm} \\
\beaa
|\,T_{n}(t)\,| & \leq & C_n \;\!J_n(t)\, J_{n+2}(t) \, \|\;\!Du(\cdot,t)\,\ni \\
 & \leq & C_n \;\!J_{n+2}(t) \,(J_3(t) + J_0(t)\:\!) \\
 & \leq & C_n \;\!J_{n+2}(t) \,(J_{n+1}(t) + J_0(t)\,\!).
\eeaa
In the last estimate we have used that $n \geq 2$. \\
\mbox{} \vspace{-0.800cm} \\

These bounds prove that \\
\mbox{} \vspace{-0.650cm} \\
$$\f{d}{dt}\,J_{n+1}^2(t) \,\leq\,
C_n \;\!J_{n+2}(t)\:(J_{n+1}(t) + J_0(t)\:\!) \,-\, 2\:\!J^2_{n+2}(t), $$
\mbox{} \vspace{-0.200cm} \\
and boundedness of $J_{n+1}(t)$ in $0 \leq t < T$ follows.
\mbox{} \hfill $\Box$

%
%
\mbox{} \vspace{-0.650cm} \\

If $\;\!\|\,u(\cdot,t)\,\ni$
is bounded in some interval $0 \leq t < T$,
then all functions $J_n(t)$ are also bounded in $0 \leq t < T $
and, using (\ref{eq.3.s2}),
all space derivatives of $u(\cdot,t) $
are therefore bounded in maximum norm.
Estimates for the pressure and its
derivatives follow from the
Poisson equation
satisfied by $p(\cdot,t) $,\\
\mbox{} \vspace{-0.750cm} \\
\be
\label{Poisson.equation.pressure}
-\, \Delta \:\!p \,=\!
\sum_{\;i,\,j\,=\,1}^{3} \!\!D_i D_j(\:\!u_i\:\!u_j)\:\!.
\ee
\mbox{} \vspace{-0.150cm} \\
Time derivatives and mixed derivatives of $u$
can be expressed by space derivatives,
using the differential equation (\ref{Navier.Stokes.eqs}).
$\!$Hence, if
$\;\!\|\,u(\cdot,t)\,\ni$ is bounded in some interval $0 \leq t < T \!\;\!$,
then all derivatives of $u$ are bounded in the same interval,
and therefore the solution $(u,p)$
can be continued as a $C^\infty$ solution beyond $T\!\;\!$.
This proves Theorem \ref{t.1.max}.

\end{section}

\mbox{} \vspace{-0.920cm} \\
%
%
%
\begin{section}{Blow--up of
$\:\!\|\,\cD u(\cdot,t)\,\|_{\mbox{}_{\scriptstyle L^{q}}}\!\;\!$
for $ \f{3}{2} < q \leq 2 $}
\label{s3}
\mbox{} \vspace{-0.950cm} \\
\setcounter{equation}{0}

A physically important quantity is the vorticity,
$\;\!\omega(\cdot,t) = \nabla \!\,\!\times\!\:\! u(\cdot,t)$,
and the total enstrophy of the flow,
given by \\
\mbox{} \vspace{-0.950cm} \\
$$\int_{\mbox{}_{\mbox{\scriptsize $\!\;\!\R^3$}}} \!\!\!\:\!
|\,\omega(x,t)\,|^2\, dx\,\!. $$
\mbox{} \vspace{-0.300cm} \\
If \\
\mbox{} \vspace{-1.025cm} \\
$$\hat v(k) \,=\: (2\:\!\pi)^{-3/2} \!
\int_{\mbox{}_{\mbox{\scriptsize $\!\;\!\R^3$}}} \!\!\!\:\!
e^{-i\,k\;\!\cdot\;\!x} \;\!v(x)\, dx $$
\mbox{} \vspace{-0.100cm} \\
denotes the Fourier transform of a 3D field $v(x)$,
then we have \\
\mbox{} \vspace{-0.600cm} \\
$$ \hat\omega(k,t) \,=\:
i\;\!k \times \hat u(k,t),
\qquad
k \cdot \hat{u}(k,t) \,=\:0\:\!,$$
\mbox{} \vspace{-0.350cm} \\
and, therefore, \\
\mbox{} \vspace{-1.000cm} \\
$$|\:\hat\omega (k,t)\,| \,=\, |\,k\,|\:|\,\hat u(k,t)\,|\:\!.$$
\mbox{} \vspace{-0.200cm} \\
Using Parseval's relation,
one finds that \\
\mbox{} \vspace{-0.950cm} \\
\beaa
\int_{\mbox{}_{\mbox{\scriptsize $\!\;\!\R^3$}}} \!\!\!\:\!
 |\:\omega(x,t)\,|^2\, dx
 &=&
 \int_{\mbox{}_{\mbox{\scriptsize $\!\;\!\R^3$}}} \!\!\!\:\!
 |\:\hat \omega(k,t)\,|^2 \, dk \\
 &=&
 \sum_{i,\,j\,=\,1}^{3}
 \int_{\mbox{}_{\mbox{\scriptsize $\!\;\!\R^3$}}} \!\!\!\:\!
 k_i^2 \;\!|\,\hat u_j(k,t)\,|^2\, dk \\
 &=&
 \|\,\cD u(\cdot, t)\,\|^2.
\eeaa
\mbox{} \vspace{-0.300cm} \\
%
In this section,
we establish the blow--up of
$\;\!\|\,\cD u(\cdot, t)\,\|_{\mbox{}_{\scriptstyle L^{q}}}\!\;\!$,
for $\;\!\f{3}{2} < q \leq 2$,
with emphasis on $ q = 2 $.
The remaining case $ q > 2 \:\!$ is covered
in Section \ref{s4}.

%
%
\mbox{} \vspace{-0.900cm} \\
\begin{subsection}{Boundedness of $\;\!\|\,\cD u(\cdot,t)\,\|$
implies boundedness of $\;\!\|\,u(\cdot,t)\,\ni$}
\mbox{} \vspace{-0.250cm} \\
The basic result here is the following. \\
\mbox{} \vspace{-0.950cm} \\
%
%
%
%
\begin{theorem}
\label{t.3.1}
\label{theorem.Du.L2}
If \\
\mbox{} \vspace{-0.500cm} \\
\mbox{} \hspace{+4.000cm}
$ {\displaystyle
\sup_{0\,\leq\, t \,<\, T}
\!\;\!\|\,\cD u(\cdot,t) \,\|
\,=:\;\!
C_{\mbox{}_{\scriptstyle \!\;\!2}}
<\, \infty,
} $
\hfill $(3.1a)$ \\
\mbox{} \vspace{-0.250cm} \\
then \\
\mbox{} \vspace{-0.630cm} \\
\mbox{} \hspace{+4.750cm}
$ {\displaystyle
\label{eq.2.s3}
\sup_{0 \,\leq\, t \,<\, T} \!\;\!\|\,u(\cdot,t)\,\ni
\:\!\leq\;\! K
} $
\hfill $(3.1b)$ \\
\mbox{} \vspace{+0.050cm} \\
for some bound $ \:\!K \!\:\!$
that depends only on
$ \;\!C_{\mbox{}_{\scriptstyle \!\;\!2}}\!\;\!, \:\!T\!\:\!$
and
$ \;\!\|\, \hat{f}\,\|_{\mbox{}_{\scriptstyle L^{1}}} \!\;\!$,
where $ \hat{f} $ denotes the Fourier transform of
the initial state $f$.
In particular,
if $ \;\!T_{\!\;\!f} < \infty $,
then
$ {\displaystyle
\!
\sup_{0 \,\leq\, t \,<\, T_{\!\;\!f}} \!\!
\|\,\cD u(\cdot,t) \,\| \,=\, \infty
} $.
\end{theorem}
\setcounter{equation}{1}
%
%
\mbox{} \vspace{-0.470cm} \\
{\bf Proof:}
Taking the Fourier transform
of the Navier--Stokes equations,
we get\footnote{%
%
%
For more applications of Fourier transforms to (1.1)
along these lines, see Section 3.2
and \cite{Heywood1994, KreissHagstromLorenzZingano2002}.
}
%
%
\mbox{} \\
\mbox{} \vspace{-0.650cm} \\
$$\hat u_t +\;\! (\:\!u \!\;\!\cdot\!\;\! \nabla u\:\!)\hat{~} +\,
(\:\!\nabla p\:\!)\hat{~} =\, -\,|\,k\,|^2 \,\hat u, $$
\mbox{} \vspace{-0.250cm} \\
or,
setting
$\;\!Q(x,t) \;\!=\;\! - \;\!u \!\;\!\cdot\!\;\! \nabla u \;\!-\;\! \nabla p$, \\
\mbox{} \vspace{-0.650cm} \\
$$\hat u_t \;\!=\, -\,|\,k\,|^2 \,\hat u
\,+\, \hat Q(k,t), $$
\mbox{} \vspace{-0.450cm} \\
with
$ \hat{u}(\cdot,0) = \hat{f} $.
Since \\
\mbox{} \vspace{-1.000cm} \\
$$(\:\!u \!\;\!\cdot\!\;\! \nabla u\:\!)\hat{~} \;\!=\;\!
-\:\hat Q(k,t) \;\!-\;\! (\:\!\nabla p\:\!)\hat{~} $$
\mbox{} \vspace{-0.370cm} \\
is the orthogonal decomposition of the vector
$(\:\!u \!\;\!\cdot\!\;\! \nabla u\:\!)\hat{~}$
into a vector orthogonal to $k$ and \linebreak
a vector parallel to $k$,
it follows that
$|\,\hat Q(k,t)\,| \;\!\leq\;\! |\,(\;\!u \!\;\!\cdot\!\;\! \nabla u\;\!)\hat{~}\;\!|$.
One obtains,
for each $k$,

$$\hat u(k,t) \,=\, e^{-\;\!|\,k\,|^2 \;\!t} \hat f(k)
\,+ \int_0^t \!\;\!e^{-\;\!|\,k\,|^2\;\!(t-s)} \;\!\hat Q(k,s) \, ds, $$
and so, \\
\mbox{} \vspace{-0.950cm} \\
$$|\,\hat u(k,t)\,| \,\leq\,
e^{-\;\!|\,k\,|^2 \;\!t} \;\!|\, \hat f(k) \,|
\,+ \int_0^t \!\;\!e^{-\;\!|\,k\,|^2\;\!(t-s)} \,
|\,(\:\!u \!\;\!\cdot\!\;\! \nabla u\:\!)\hat{~}(k,s)\,|  \, ds. $$
\mbox{} \vspace{-0.200cm} \\
Integrating in $k \in \R^3$ one finds that \\
\mbox{} \vspace{-1.050cm} \\
\beaa
(\:\!2\:\!\pi\:\!)^{3/2} \: \|\,u(\cdot,t)\,\ni &\leq &
\int_{\mbox{}_{\mbox{\scriptsize $\!\;\!\R^3$}}} \!\!\!\;\!
 |\,\hat u(k,t)\,|\: dk \\
&\leq & \|\,\hat f \,\|_{\mbox{}_{\scriptstyle L^1}}
\;\!+\;\!
\int_0^t \!\;\!
\int_{\mbox{}_{\mbox{\scriptsize $\!\;\!\R^3$}}} \!\!\!\:\!
 e^{-\;\!|\,k\,|^2\;\!(t-s)} \,
 |\,(\:\!u \!\;\!\cdot\!\;\!\nabla u\:\!)\hat{~}(k,s)\,| \, dk \, ds.
\eeaa
We then apply the Cauchy--Schwarz inequality
to bound the inner integral
on the right-hand side by \\
\mbox{} \vspace{-1.050cm} \\
$$I_1^{1/2} \;\!I_2^{1/2}, $$
\mbox{} \vspace{-0.700cm} \\
where \\
\mbox{} \vspace{-1.350cm} \\
\beaa
\label{eq.3.ex}
 I_1 \,=\;\!
\int_{\mbox{}_{\mbox{\scriptsize $\!\;\!\R^3$}}} \!\!\!\:\!
 e^{-\;\!2\,|\,k\,|^2\;\!(t-s)}\, dk
 \;=\;
 C \cdot (\:\!t - s\:\!)^{-3/2}
\eeaa
\mbox{} \vspace{-0.450cm} \\
and \\
\mbox{} \vspace{-1.470cm} \\
\beaa
  I_2 &=&
  \int_{\mbox{}_{\mbox{\scriptsize $\!\;\!\R^3$}}} \!\!\!\:\!
  |\,(\:\!u \!\;\!\cdot\!\;\! \nabla u\:\!)(x,s)\,|^2 \, dx \\
  &\leq &
  \!\!\;\! C \cdot \|\,u(\cdot,s)\,\ni ^2 \,\|\,\cD u(\cdot,s)\,\|^2 \\
  &\leq &
  \!\!\;\! C \cdot C_{\scriptstyle \!\;\!2}^{\;\!2} \cdot
  \|\,u(\cdot,s)\,\ni^2,
\eeaa
\mbox{} \vspace{-0.500cm} \\
using Parseval's relation
and $(3.1a)$.
Thus
we have shown the estimate \\
\mbox{} \vspace{-0.550cm} \\
$$ (\:\!2\:\!\pi\:\!)^{3/2} \: \|\,u(\cdot, t)\,\ni
\,\leq\;
\|\,\hat f \,\|_{\mbox{}_{\scriptstyle L^1}}
+\, C \cdot C_{\mbox{}_{\scriptstyle \!\;\!2}} \!\:\!
\int_0^t (\:\!t - s\:\!)^{-3/4} \: \|\,u(\cdot,s)\,\ni \, ds
\qfor
0 \leq t < T $$
\mbox{} \vspace{-0.200cm} \\
for some constant $ C \!\:\!> 0 $.
By the singular Gronwall's lemma
given in
Lemma \ref{Lemma.singular.Gronwall.1}
next,
boundedness of $\;\!\|\,u(\cdot,t)\,\ni$
in the interval $0 \leq t < T$ follows,
as claimed.
\mbox{} \hfill $\Box$
%

%
\mbox{} \vspace{-0.050cm} \\
%
%
%
%
%
{\sc Remark.}
By the previous argument
and Lemma \ref{Lemma.singular.Gronwall.2} below,
we can see that
condition $(3.1a)$
also implies \\
\mbox{} \vspace{-1.020cm} \\
\be
\label{eq.2.s3}
\sup_{0 \,<\, t \,<\, T} \!\;\!t^{\,3/4} \, \|\,u(\cdot,t)\,\ni
\:\!\leq\;\! K_{\mbox{}_{\scriptstyle \!\;\!2}} \!\:\!(T)
\, \|\,f\,\|_{\mbox{}_{\scriptstyle L^{2}}}
\ee
\mbox{} \vspace{-0.070cm} \\
for some bound
$ \:\!K_{\mbox{}_{\scriptstyle \!\;\!2}} \!\;\!(T) > 0 \:\!$
that depends on
the values of
$ \;\!C_{\mbox{}_{\scriptstyle \!\;\!2}}\!\;\!, \;\!T\!\;\!$
only.

\newpage
\mbox{} \vspace{-1.300cm} \\

The following result is an important version
of Gronwall's lemma frequently used
for partial differential equations,
as in the proof of Theorem 3.1 above.

%
\mbox{} \vspace{-0.900cm} \\
%
%
%
%
%
\begin{lemma}
\label{Lemma.singular.Gronwall.1}
Let $ A \geq 0, \:\!B > 0, \;\!0 < \kappa < 1$.
Let $\;\!\phi \in C^{0}(\:\![\,0,T\;\![\;\!)$
satisfy \\
\mbox{} \vspace{-0.625cm} \\
\be
\label{eq.3.s3}
0 \,\leq\;\! \phi(t) \,\leq\:\!
A \;\!+\;\! B \! \int_0^t (t-s)^{-\kappa} \, \phi(s) \,ds
\qfor 0 \leq t < T.
\ee
\mbox{} \vspace{-0.200cm} \\
Then
$ \;\!\phi(t) \leq K\!\;\!(T) \:\!A \:\!$
for all $\;\!0 \leq t < T\!\;\!$,
where $K\!\;\!(T) > 0$
depends on $\:\!B\!\;\!, \:\!\kappa, T\!\;\!$
only.
\end{lemma}
%
%
\mbox{} \vspace{-0.300cm} \\
{\bf Proof:}
For convenience,
we provide a proof for
Lemma \ref{Lemma.singular.Gronwall.1}
that can be easily extended
to other useful similar statements
like
Lemma \ref{Lemma.singular.Gronwall.2}
below.
To this end,
we choose $\epsilon > 0$ with \\
\mbox{} \vspace{-1.150cm} \\
$$ \int_0^\epsilon \! \:\! \xi^{-\kappa} \, d\xi
\,=\, \f{\,\epsilon^{1-\kappa}}{1 - \kappa}
\,\leq\, \f{1}{2B}$$
\mbox{} \vspace{-0.150cm} \\
and,
given $\;\! t \in \:\![\,0, T\;\![\;\!$ arbitrary,
we take $ \;\!t_0 \!\:\!\in [\,0, t\;\!] $
such that
$ {\displaystyle
\;\! \phi(t_0)\;\!=
\max_{0 \,\leq\, s \,\leq\, t} \!\phi(s)
} $. \\
\mbox{} \vspace{-0.150cm} \\
{\sc Case I:} $ t_0 \geq\:\! \epsilon$.
Then we have \\
\mbox{} \vspace{-0.950cm} \\
\beaa
\phi(t_0) &\leq & A \,+\;\! B \! \int_0^{t_0\;\!-\,\epsilon} \hspace{-0.400cm}
(t_0 - s)^{- \kappa} \;\! \phi(s)\, ds
\;+\; B \!\int_{t_0\;\!-\,\epsilon}^{t_0} \hspace{-0.250cm}
(t_0 - s)^{- \kappa} \;\! \phi(s)\, ds \\
&\leq & A \,+\;\! B\;\! \epsilon^{- \kappa} \!\!\:\!
\int_0^{t_0} \!\! \phi(s) \,ds
\;+\; B\, \f{1}{2B}\; \phi(t_0)\:\!,
\eeaa
\mbox{} \vspace{-0.600cm} \\
so that \\
\mbox{} \vspace{-1.100cm} \\
$$\phi(t) \,\leq\, \phi(t_0) \,\leq\,
2\:\!A \,+\, 2\:\!B\;\!\epsilon^{- \kappa} \!\!\:\!
\int_0^{t} \!\;\!\phi(s)\, ds\:\!. $$
\mbox{} \vspace{-0.050cm} \\
{\sc Case II:} $0 \leq t_0 \!\:\!\leq \epsilon $.
We have \\
\mbox{} \vspace{-1.100cm} \\
\beaa
\phi(t_0) &\leq & A \,+\;\! B \! \int_0^{t_0} \!
(t_0 - s)^{- \kappa} \;\! \phi(s)\, ds \\
&\leq &
A \,+\;\! B \,\phi(t_0) \! \int_0^{\epsilon} \xi^{- \kappa} \, d\xi \\
&\leq &
A \,+\, \f{1}{2} \:\phi(t_0)\:\!,
\eeaa
\mbox{} \vspace{-0.600cm} \\
and so, \\
\mbox{} \vspace{-1.150cm} \\
$$\phi(t) \,\leq\, \phi(t_0) \,\leq\, 2 A. $$

\mbox{} \vspace{-0.750cm} \\

We have thus shown that \\
\mbox{} \vspace{-0.675cm} \\
$$\phi(t) \,\leq\, 2A \,+\, 2 \:\!B\;\! \epsilon^{- \kappa} \!\!\:\!
\int_0^{t} \!\;\! \phi(s) \,ds
\qfor 0 \leq t < T. $$
\mbox{} \vspace{-0.250cm} \\
This gives,
by standard Gronwall,
$ {\displaystyle
\;\!
\phi(t) \;\!\leq\;\! 2\:\!A\: \mbox{exp}\,
\bigl\{\;\!2\:\!B\;\!\epsilon^{- \kappa} \;\!T \;\!\bigr\}
\;\!
} $
for all $ 0 \leq t < T $.
\mbox{} \hfill $\Box$
%
%
\newpage
\mbox{} \vspace{-1.250cm} \\

In a similar way,
the following generalization of
Lemma \ref{Lemma.singular.Gronwall.1}
can be easily obtained. \\
\mbox{} \vspace{-0.950cm} \\
%
%
%
%
\begin{lemma}
\label{Lemma.singular.Gronwall.2}
Let $ A \geq 0 $,
$ B_{\mbox{}_{\scriptscriptstyle \!\;\!1}}\!\;\!, ...,
  B_{\mbox{}_{\scriptstyle \!\;\!n}} > 0$.
Let $\;\!\phi \in C^{0}(\:\![\,0,T\;\![\;\!)$
satisfy \\
\mbox{} \vspace{-0.675cm} \\
\be
\label{eq.4.s3}
0 \,\leq\;\! \phi(t) \,\leq\;\!
A \,+\;\! \sum_{j\,=\,1}^{n}
B_{\!\;\!j} \!\!\;\!
\int_0^t \!\:\!
s^{\!- \,\alpha_{\!\;\!j}} \;\!
(t-s)^{- \,\beta_{\!\;\!j}} \, \phi(s) \,ds
\qfor 0 \leq t < T,
\ee
\mbox{} \vspace{-0.275cm} \\
with
$ \:\!\alpha_{j} \!\:\!\geq 0 $,
$ \beta_{j} \!\;\!\geq 0 \:\!$
verifying
$ \;\!\alpha_{j} \!\;\!+\;\!\beta_{j} < 1 $
for all $ \:\!1 \leq j \leq n $.
Then $\,\phi(t) \leq \!\;\!K\!\;\!(T) \:\!A $
for all $\;\!0 \leq t < T \!\;\!$,
with the quantity $\:\!K\!\;\!(T) > 0 \,\!$
depending only on $\;\!T \!\,\!, \;\!n $
and $ \:\!B_{\!\;\!j}, \:\!\alpha_{j}, \:\!\beta_{j} $,
$ 1 \leq j \leq n $.
\end{lemma}
%
%
\mbox{} \vspace{-1.000cm} \\

\end{subsection}

%
%
\mbox{} \vspace{-1.100cm} \\
\begin{subsection}{Blow--up of
$\;\!\|\,\cD u(\cdot,t)\,\|_{\mbox{}_{\scriptstyle L^q}}$
for $ \f{3}{2} < q \leq 2 $}
\mbox{} \vspace{-0.200cm} \\
We now extend the proof of
Theorem \ref{theorem.Du.L2}
by using H\"older's inequality
instead of the Cauchy--Schwarz inequality
and
the Hausdorff--Young inequality
(see e.g.$\;$\cite{Grafakos2008}, p.$\;$104)
instead of Parseval's relation.

%
%
%
%
\begin{theorem}
\label{t.3.DLq}
\label{theorem.Du.Lq.(i)}
Let $\,\f{3}{2} < q \leq 2 $.
If \\
\mbox{} \vspace{-0.250cm} \\
\mbox{} \hspace{+4.650cm}
$ {\displaystyle
\sup_{0\,\leq\, t \,<\, T} \!\;\!
\|\,\cD u(\cdot,t) \,\|_{\mbox{}_{\scriptstyle L^{q}}}
=:\;\! C_{\mbox{}_{\scriptstyle \!\;\!q}} <\;\! \infty,
} $
\hfill $(3.5a)$ \\
\mbox{} \vspace{-0.250cm} \\
then \\
\mbox{} \vspace{-0.750cm} \\
\mbox{} \hspace{+5.625cm}
$ {\displaystyle
\sup_{0 \,\leq\, t \,<\, T} \!\;\!\|\,u(\cdot,t)\,\ni
\:\!\leq\;\! K
} $
\hfill $(3.5b)$ \\
\setcounter{equation}{5}
\mbox{} \vspace{-0.050cm} \\
for some bound $ \;\!K\!\;\!$
that depends only on
$ \;\!q, \;\!C_{\mbox{}_{\scriptstyle \!\;\!q}}\!\;\!, \,\!T\!\:\!$
and
$ \,\!\|\, \hat{f}\,\|_{\mbox{}_{\scriptstyle L^{1}}} \!\:\!$,
where $ \hat{f} $ denotes the Fourier transform of
the initial state $f$.
In particular,
if $ \;\!T_{\!\;\!f} < \infty $,
then \\
\mbox{} \vspace{-0.750cm} \\
\be
\label{eq.Dq.3halves}
\sup_{0 \,\leq\, t \,<\, T_{\!\;\!f}} \!\!
\|\,\cD u(\cdot,t) \,\||_{\mbox{}_{\scriptstyle L^{q}}}
\!\;\!=\;\! \infty.
\ee
\end{theorem}
%
%
\mbox{} \vspace{-0.250cm} \\
{\bf Proof:}
We use the same notation as in the proof
of Theorem \ref{theorem.Du.L2}.
Applying H\"older's inequality
to the integral \\
\mbox{} \vspace{-0.950cm} \\
$$I \;\!=
\int_{\mbox{}_{\mbox{\scriptsize $\!\;\!\R^3$}}} \!\!\!\:\!
e^{-\,|\,k\,|^2\;\!(t-s)} \,
|\,(\:\!u \!\;\!\cdot\!\;\! \nabla u\:\!)\hat{~}(k,s)\,|\: dk,$$
we obtain the bound \\
\mbox{} \vspace{-0.850cm} \\
$$I \,\leq\, I_1^{1/q} \:\!I_2^{1/q'}\!,
\quad\;\; \f{1}{q}+\f{1}{q'} = 1,$$
\mbox{} \vspace{-0.650cm} \\
with \\
\mbox{} \vspace{-0.850cm} \\
$$I_1 \,=
\int_{\mbox{}_{\mbox{\scriptsize $\!\;\!\R^3$}}} \!\!\!\:\!
e^{-\;\!q\,|\,k\,|^2\;\!(t-s)}\, dk
\;=\: \Bigl(\;\! \frac{\pi}{q}\,\Bigr)^{\!3/2} \:\!(t-s)^{-3/2} $$
and \\
\mbox{} \vspace{-1.300cm} \\
\beaa
I_2^{\,1/q'} &=&
\|\,(u \!\;\!\cdot\!\;\! \nabla u)\hat{~}(\cdot,s)\,\|_{\mbox{}_{\scriptstyle L^{q'}}} \\
&\leq &
3\: \|\,u\!\;\!\cdot\!\;\! \nabla u \,(\cdot,s)\,\|_{\mbox{}_{\scriptstyle L^q}} \\
&\leq &
C \:\|\,u(\cdot,s)\,\ni \,\|\, \cD u(\cdot,s)\,\|_{\mbox{}_{\scriptstyle L^q}}
\!\;\!,
\eeaa
\mbox{} \vspace{-0.550cm} \\
where
in the first estimate
we have used the Hausdorff--Young inequality,
since $\;\! q \leq 2 $.
We obtain that \\
\mbox{} \vspace{-0.650cm} \\
$$ (\:\!2\:\!\pi\:\!)^{3/2} \,
\|\,u(\cdot,t)\,\ni
\,\leq\;
\|\, \hat{f}\,\|_{\mbox{}_{\scriptstyle L^{1}}}
+\, C \! \int_0^t (t-s)^{-\;\!\kappa} \,
\|\,u(\cdot,s)\,\ni \, \|\,\cD u(\cdot,s)\,\|_{\mbox{}_{\scriptstyle L^q}}
\, ds\:\!, $$
\mbox{} \vspace{-0.240cm} \\
where
$ \;\! \kappa = \f{3}{2q} < 1 $,
in view that $\:\!q > \f{3}{2}$.
By Lemma \ref{Lemma.singular.Gronwall.1},
the result now readily follows.
\mbox{} \hfill $\Box$
%

%
\mbox{} \vspace{+0.090cm} \\
{\sc Remarks.}
({\em i\/})
By the argument above
and Lemma \ref{Lemma.singular.Gronwall.2},
we see that
$(3.5a)$
also gives \\
\mbox{} \vspace{-0.600cm} \\
\be
\label{eq.7.s3}
\sup_{0 \,<\, t \,<\, T} \!\;\!t^{\,3/4} \, \|\,u(\cdot,t)\,\ni
\:\!\leq\;\! K_{\mbox{}_{\scriptstyle \!\;\!q}} \!\;\!(T)
\, \|\,f\,\|_{\mbox{}_{\scriptstyle L^{2}}}
\ee
\mbox{} \vspace{-0.100cm} \\
for some coefficient
$ \,\!K_{\mbox{}_{\scriptstyle \!\;\!q}} \!\;\!(T) > 0 \;\!$
that depends on
the values of
$ \:\!q, \;\!C_{\mbox{}_{\scriptstyle \!\;\!q}}\!\;\!, \;\!T\!\;\!$
only. \\
\mbox{} \vspace{-0.350cm} \\
({\em ii\/})
The Navier--Stokes equations on the whole $\R^3\!$
enjoy the following scaling invariance:
If $\bigl(u(x,t),p(x,t)\bigr)$
solves the Navier--Stokes equations,
then, for every scaling parameter $\lam > 0$,
$\bigl(\:\!\lam \:\!u(\lam \:\!x,\lam^2 \,\!t),
\lam ^2 \:\!p(\lam \:\!x, \lam ^2 \:\!t)\:\!\bigr)$
solves the same equations.
The norms $\|\:u\:\|_{\mbox{}_{\scriptstyle L^3}}$
and \linebreak
$\|\,\cD u\,\|_{\mbox{}_{\scriptstyle L^{3/2}}}$,
which appear in the limiting values $q = 3$ in Theorem \ref{thm.1.3}
and $q = 3/2$ in Theorem \ref{t.3.DLq},
are also invariant under such $\lam$ scalings.
Better understanding of the scale invariant norms,
$\|\,u(\cdot,t)\,\|_{\mbox{}_{\scriptstyle L^3}}$
and
$\|\,\cD u(\cdot,t)\,\|_{\mbox{}_{\scriptstyle L^{3/2}}}\!\;\!$,
as $ t \to T_{\!\;\!f} $,
is likely to be important  for
further progress on the blow--up question. \\
\mbox{} \vspace{-0.400cm} \\
({\em iii\/})
Theorem \ref{t.3.DLq} can also be deduced
from the results of Section \ref{s4}
on $ \;\!\|\, u(\cdot,t) \,\|_{\mbox{}_{\scriptstyle L^{q}}} \!\;\!$,
$ q > 3 $,
using the Sobolev inequality (\ref{SNG.Du.Lq}).

\end{subsection}

%
%
%
\mbox{} \vspace{-0.930cm} \\
\begin{subsection}{A differential inequality for $\;\!\|\,\cD u(\cdot,t)\,\|$}
\mbox{} \vspace{-0.275cm} \\
Here
we derive
the estimate (\ref{lowerbound.Du.Lq})
for $\;\!\|\,\cD u(\cdot,t) \,\| \;\!$
from a nonlinear differential inequality
satisfied by $\:\!\|\,\cD u(\cdot,t) \,\| \:\!$
whether $ T_{\!\;\!f} \!\;\!$ is finite or not.
$\!\:\!$This method dates back to Leray \cite{Leray1934}. \\
\mbox{} \vspace{-1.000cm} \\
%
%
%
%
\begin{theorem}
\label{Differential.Inequality.Du.L2}
There is an absolute constant
$ {\displaystyle
\,0 < \!\;\!K
\!\:\!<
\mbox{\small $ {\displaystyle \frac{1}{\:\!32\:\!} }$}
} $
such that \\
\mbox{} \vspace{-0.700cm} \\
\be
\label{eq.8.s3}
\f{d}{dt}\: \|\,\cD u(\cdot,t)\,\|^2
\;\!\leq\, K \: \|\,\cD u(\cdot,t)\,\|^6 \!\;\!,
\quad \;\;
\forall \;\,
0 \leq t < T_{\!\;\!f}.
\ee
\end{theorem}
%
%
\mbox{} \vspace{-0.200cm} \\
{\bf Proof:}
We have,
using
$ \;\!\nabla \!\cdot u \;\!=\;\!0 $, \\
\mbox{} \vspace{-1.150cm} \\
\beaa
\f{1}{2}\: \f{d}{dt}\, \|\,\cD u(\cdot,t)\,\|^2
&=&
\sum_{j\,=\,1}^{3} \:\!(D_j u, D_j u_t) \\
&=&
\!\!\!- \sum_{j\,=\,1}^{3} \:\!
(D_j u, D_j (u \!\;\!\cdot\! \nabla u)\:\!)
\;-\;
\|\,\cD^2 u(\cdot,t)\,\|^2 \\
&=:&
S(t) \,-\, \|\,\cD^2 u(\cdot,t)\, \|^2 \!\;\!.\\
\eeaa
\mbox{} \vspace{-1.150cm} \\
Since
$(D_j u, u \,\!\cdot \nabla D_j u\:\!) \;\!=\;\! 0$,
because
$ \nabla \!\;\!\cdot u \;\!=\;\!0 $,
the nonlinear term $S(t)$ can be estimated by
$\;\!C \,\|\,\cD u(\cdot,t)\,\|_{\mbox{}_{\scriptstyle L^3}}^3 \!\;\!$,
for some constant $ C $.
Using the 3D Gagliardo--Nirenberg inequality \\
\mbox{} \vspace{-0.700cm} \\
$$\|\;v\;\|_{\mbox{}_{\scriptstyle L^3}} \leq\,
\Gamma \: \|\;v\;\|^{1/2} \: \|\:\cD v\;\|^{1/2}
\qquad
\forall \;\,
v \in H^{1}(\R^3)$$
\mbox{} \vspace{-0.300cm} \\
for $ v = D_j u $,
where $ \Gamma < 0.59 $,
%
%
and then
the Young's inequality
$ {\displaystyle
\;\!
a \:\!b \;\!\leq\:\!
1/4 \:a^{4} +\:\! 3/4 \;b^{4/3}
\!
} $,
one obtains (\ref{eq.8.s3}), with
$ K \!\;\!<\;\! 1/32 $,
as claimed.
\mbox{} \hfill $\Box$
%
\mbox{} \vspace{-0.100cm} \\

\noindent
{\sc Remark:}
A more involved derivation of (\ref{eq.8.s3})
in \cite{Heywood1994}, p.$\;\!$11-14,
gives that
$ {\displaystyle
\;\!0 < K \leq
\mbox{\small $ {\displaystyle \frac{1}{\;\!16 \;\!\pi^{2}} }$}
} $. \\

The next lemma shows
how nonlinear differential inequalities
such as (\ref{eq.8.s3}) above
can be used to derive lower bound estimates
in case of finite--time blow--up.

%
%
%
%
\begin{lemma}
\label{How.differential.inequalities.are.used}
Let $\,\mbox{\em w} \in \!\;\!C^{1}([\;\!0, T_{\!\;\!f}\,\![\;\!) $
be a positive function satisfying a differential inequality \\
\mbox{} \vspace{-0.100cm} \\
\mbox{} \hspace{+4.500cm}
$ {\displaystyle
\mbox{\em w}^{\prime}(t) \,\leq\;\! K \,
\mbox{\em w}(t)^{\mbox{}^{\scriptstyle \!\;\!\alpha}}
\qquad
\forall \;\,
0 \leq t < T_{\!\;\!f}
} $
\hfill $(3.9a)$ \\
\mbox{} \vspace{-0.120cm} \\
for some given
constants
$ \:\!K \!\:\!> 0 $, $ \alpha > 1 $.
If $\: T_{\!\;\!f} < \infty \;\!$
and
$ {\displaystyle
\!\!\!
\sup_{\;\;0 \,\leq\, t \,<\, T_{\!\;\!f}}
\!\!\!\!\mbox{\em w}(t)
\,=\, \infty
} $,
$\,$then
we have \\
\mbox{} \vspace{-0.175cm} \\
\mbox{} \hspace{+2.500cm}
$ {\displaystyle
\mbox{\em w}(t) \,\geq\,
\biggl(\;\!
\mbox{\small $ {\displaystyle
\frac{1}{K \!\cdot (\alpha - 1)\:\!} }$}
\,\biggr)^{\mbox{}^{\scriptstyle \!\!\!
\frac{\scriptstyle 1}{\scriptstyle \:\!\alpha - 1\:\!} }}
\hspace{-0.420cm} \cdot \hspace{+0.120cm}
\bigl(\;\! T_{\!\;\!f} -\;\!t
\;\!\bigr)^{\mbox{}^{\scriptstyle \!\!\!- \,
\frac{\scriptstyle 1}{\scriptstyle \:\!\alpha - 1\:\!} }}
\qquad
\forall \;\,
0 \leq t < T_{\!\;\!f}.
} $
\hfill $(3.9b)$ \\
\end{lemma}
\setcounter{equation}{9}
%
%
\mbox{} \vspace{-0.670cm} \\
{\bf Proof:}
Given $ \;\!t_0 \in [\,0, T_{\!\;\!f}\;\![\;\!$ arbitrary,
and setting
$ \;\!\mbox{\tt v} = \mbox{\tt v}(t) $
by
$ \;\!\mbox{\tt v}^{\prime}(t) =
K \,\!\mbox{\tt v}(t)^{\mbox{}^{\scriptstyle \!\;\!\alpha}} \!\!\:\!$,
$ \mbox{\tt v}(t_0) = \mbox{w}_0 $,
where
$ \,\!\mbox{w}_0 \!:= \mbox{w}(t_0) $,
we have
$\:\!\mbox{\tt v}(t) $
defined for all
$ {\displaystyle
t_0 \leq t < t_{\ast} \!\!\;\!:=\;\!
t_0 \!\;\!+ 1/\bigl(\:\!K\!\;\!(\alpha - 1) \:\!\mbox{w}_{0}^{\;\!\alpha -1})
} $.
Moreover,
one has
$\;\!\mbox{w}(t) \leq \mbox{\tt v}(t) $
for all $ \;\!t_0 \leq t < t_{\ast} $,
with
$ \;\!\mbox{\tt v}(t) $
given by \\
\mbox{} \vspace{-0.950cm} \\
$$
\mbox{\tt v}(t) \:=\;
\mbox{w}_{0} \cdot \Bigl(\;\! 1 \;\!-\:\! K \!\cdot\!\;\! (\alpha - 1)
\,\mbox{w}_{0}^{\,\alpha -1} \:\!\bigl(\:\!t - t_0 \bigr)
\:\!\Bigr)^{\mbox{}^{\scriptstyle \!\!\!-\,
\frac{\scriptstyle 1}{ \scriptstyle \alpha - 1\:\!} }}
\!\!,
\qquad
t_0 \leq t < t_{\ast},
$$
\mbox{} \vspace{-0.350cm} \\
so that,
in particular,
$ \;\!\mbox{\tt v}(t) \;\!\mbox{\small $\nearrow$}\:\! \infty \:\!$
as $\;\! t \;\!\mbox{\small $\nearrow$}\:\! t_{\ast} $.
This gives
$ \:\!t_{\ast} \leq T_{\!\;\!f} $,
which is $(3.9b)$ above.
\mbox{} \hfill $\Box$
%

From (\ref{eq.8.s3}),
we see that
$ \;\!\mbox{w}^{\prime}(t) \leq K \;\!\mbox{w}(t)^{3} $
for all $ 0 \leq t < T_{\!\;\!f} $,
where
$ \;\!\mbox{w}(t) := \|\, \cD u(\cdot,t) \,\|^{2} \!\:\!$,
$ 0 < \!\;\!K\!\;\! \leq 1/(16\;\!\pi^{2}) $.
Assuming $ T_{\!\;\!f} \!\;\!<\!\;\! \infty $,
we have
$\;\!\sup_{\;\!0 \,\leq\, t \,<\, T_{\!\;\!f}} \!\mbox{w}(t)
\;\!=\;\!\infty $
by Theorem~\ref{theorem.Du.L2}.
Therefore,
by Lemma \ref{How.differential.inequalities.are.used},
we get
the following result
(which dates back to Leray \cite{Leray1934}).

%
%

%
%
%
%
\begin{theorem}
\label{Blow.up.rate.Du.L2}
Assuming that $ \;\!T_{\!\;\!f} \!\;\!< \infty $,
we must then have \\
\mbox{} \vspace{-0.700cm} \\
\be
\label{eq.10.s3}
\|\,\cD u(\cdot,t)\,\| \;\geq\;
\f{c}{\;\!(\:\!T_{\!\;\!f} -\;\! t\:\!)^{1/4}}
\qquad
\forall \;\,
0 \leq t < T_{\!\;\!f},
\ee
\mbox{} \vspace{-0.200cm} \\
for some constant
$ {\displaystyle
\, c \:\!\geq \bigl\{\,2\;\!\pi\:\!\sqrt{\:\!2\;\!}\,
\bigr\}^{\mbox{}^{\scriptstyle \!1/2}}
\!>\:\!
2.98
\;\!
} $
$($independent of $\:\!f\!\;\!, \;\!u, \:\!T_{\!\;\!f})$.
\end{theorem}
%
%
\mbox{} \vspace{-0.550cm} \\
As with the other bounds for $u(\cdot,t)$
discussed in the text,
the optimal ($=$ largest, here) \linebreak
value of the absolute constant $\;\!c\;\!$
in (\ref{eq.10.s3}) above
is not known. \\
%
%
%

%
%

\end{subsection}

\end{section}

\mbox{} \vspace{-1.500cm} \\
%
%
%
\begin{section}{Blow--up of
$\:\!\|\,u(\cdot,t)\,\|_{\mbox{}_{\scriptstyle L^q}}$
for $ 3 < q \leq \infty $}
\label{s4}
\setcounter{equation}{0}
\mbox{} \vspace{-0.950cm} \\

Using the
Helmholtz projector $P_{\scriptscriptstyle \!\;\!H}$
(see e.g.$\;$\cite{FoiasManleyRosaTemam2004,%
Galdi1994, LopesNussenzveigZheng1999}),
one can write the incompres\-sible
Navier--Stokes equations
as \\
\mbox{} \vspace{-0.150cm} \\
\mbox{} \hspace{+2.250cm}
$ {\displaystyle
u_t \;\!=\: \Delta u \,-\,
P_{\scriptscriptstyle \!\;\!H}(u \cdot\!\;\!\nabla u),
\qquad
P_{\scriptscriptstyle \!\;\!H}(u \cdot\!\;\!\nabla u)
\:=\:
u \cdot \!\;\!\nabla u \,+\;\! \nabla \!\;\!p,
} $
\hfill $(4.1a)$ \\
\mbox{} \vspace{-0.150cm} \\
and,
if $e^{\Delta t}\!\;\!$ denotes the heat semigroup,
then one obtains,
by Duhamel's principle, \\
\mbox{} \vspace{-0.100cm} \\
\mbox{} \hspace{+3.150cm}
$ {\displaystyle
u(\cdot,t) \:=\:
e^{\Delta \:\!t} f(\cdot) \,-
\int_0^t \!\;\! e^{\Delta (t-s)}
P_{\scriptscriptstyle \!\;\!H}(u \cdot\!\;\!\nabla u)(\cdot,s)\, ds.
} $
\hfill $(4.1b)$ \\
\setcounter{equation}{1}
\mbox{} \vspace{-0.120cm} \\
It is not difficult to show that
the linear operators
$P_{\scriptscriptstyle \!\;\!H}$ and $e^{\Delta t}$
commute,
and these operators also commute with differentiation,
$D_{\!\;\!j} \!\;\!=\:\! \partial/\partial x_{\!\;\!j}$.
Using the Calderon--Zygmund theory \linebreak
of singular integrals
(see e.g.$\;$\cite{Stein1970}, Ch.$\;$2),
one shows the fundamental property
that the Helmholtz projector
is bounded in
$ L^{q} $ if $ 1 < q < \infty $.
That is,
for each $ 1 < q < \infty $
there exists some constant
$ C_{\!\;\!q} \!\;\!> 0 \;\! $
such that \\
\mbox{} \vspace{-0.750cm} \\
\be
\label{eq.2.s4}
\|\, P_{\scriptscriptstyle \!\;\!H} \:\!v \,
\|_{\mbox{}_{\scriptstyle L^{q}}}
\;\!\leq\,
C_{\!\;\!q} \,
\|\: v \:
\|_{\mbox{}_{\scriptstyle L^{q}}}
\qquad
\forall \;\,
v = (v_{\mbox{}_{1}}\!\:\!, v_{\mbox{}_{2}}\!\;\!, v_{\mbox{}_{3}})
\:\!\in L^{\mbox{}^{\scriptstyle q}}\!\:\!(\R^3)
\qquad
(1 < q < \infty).
\ee
\mbox{} \vspace{-0.300cm} \\
%
We will also need here
the following well known estimate
for solutions of the heat equation$\,\!$:
given any $ 1 \leq r \leq q \leq \infty $,
and any multi--index $ \alpha $,
we have,
for all $ \:\! t > 0 \:\!$: \\
\mbox{} \vspace{-0.600cm} \\
\be
\label{eq.3.s4}
\|\, D^{\alpha} \:\! e^{\Delta t} v \:
\|_{\mbox{}_{\scriptstyle L^{q}}}
\;\!\leq\,
C \:
\|\: v \:\|_{\mbox{}_{\scriptstyle L^{r}}}
\,
t^{\mbox{}^{\scriptstyle \!\;\!-\;\! \lambda \;\!-\;\! |\;\!\alpha\;\!|/2}}
\!,
\qquad
\lambda \,=\,
\f{3}{2} \,\biggl(\;\! \f{1}{r} \;\!-\;\! \f{1}{q} \;\!\biggr)
\ee
\mbox{} \vspace{-0.240cm} \\
for all $ v \in L^{r}(\R^3) $,
with $ C > 0 $ constant
depending only on the values of
$ \:\!q, \;\!r \:\!$ and $\;\! |\;\!\alpha\;\!| $.

%
%
\mbox{} \vspace{-0.900cm} \\
\begin{subsection}{Boundedness of
$\,\!\|\,u(\cdot,t)\,\|_{\mbox{}_{\scriptstyle L^q}}\!\,\!$
implies boundedness of $\,\!\|\,\cD u(\cdot,t)\,\|$}
\mbox{} \vspace{-0.850cm} \\
The basic result here is the following. \\
\mbox{} \vspace{-0.950cm} \\
%
%
%
%
\begin{theorem}
\label{t.4.1}
\label{theorem.4.Lq}
Let $\, 3 < q \leq \infty $.
If \\
\mbox{} \vspace{-0.195cm} \\
\mbox{} \hspace{+4.750cm}
$ {\displaystyle
\sup_{0\,\leq\, t \,<\, T}
\!\;\!\|\,u(\cdot,t) \,\|_{\mbox{}_{\scriptstyle L^{q}}}
=:\;\!
C_{\mbox{}_{\scriptstyle \!\;\!q}}
<\, \infty,
} $
\hfill $(4.4a)$ \\
\mbox{} \vspace{-0.250cm} \\
then \\
\mbox{} \vspace{-0.630cm} \\
\mbox{} \hspace{+4.450cm}
$ {\displaystyle
\sup_{0 \,\leq\, t \,<\, T} \!\;\!\|\,\cD u(\cdot,t)\,\|
\,\leq\;\! K_{\!\:\!q}(T) \, \|\, \cD f \,\|
} $,
\hfill $(4.4b)$ \\
\mbox{} \vspace{+0.050cm} \\
\setcounter{equation}{4}
where $ \;\!K_{\!\:\!q}(T) > 0 \:\!$
depends on the values of
$ \;\!q, C_{\mbox{}_{\scriptstyle \!\;\!q}}\!\;\!, \:\!T\!\:\!$
only.
In particular,
if $ \;\!T_{\!\;\!f} < \infty $,
then \\
\mbox{} \vspace{-0.700cm} \\
\be
\label{eq.5.s4}
\sup_{0 \,\leq\, t \,<\, T_{\!\;\!f}} \!\!
\|\, u(\cdot,t) \,\|_{\mbox{}_{\scriptstyle L^{q}}}
\:\!=\, \infty.
\ee
\end{theorem}
%
%

%
\mbox{} \vspace{-0.100cm} \\
{\bf Proof:}
(Note that
the result for $ q = \infty $
was already shown in Section \ref{s2},
and we provide another proof here.)
Given $ \;\!3 < q \leq \infty $,
let $\;\! \f{6}{5} < r \leq 2 \;\!$
be defined by \\
\mbox{} \vspace{-0.650cm} \\
$$\f{1}{q} \;\!+\;\! \f{1}{2} \;\!=\;\! \f{1}{r}\;\!. $$
\mbox{} \vspace{-0.250cm} \\
From $(4.1b)$,
we have,
for each $ 1 \leq j \leq 3 $: \\
\mbox{} \vspace{-0.600cm} \\
$$D_{\!\;\!j}\:\!u(\cdot,t) \,=\:
e^{\Delta \:\!t} D_{\!\;\!j}f
\;\!-
\int_0^t \!
D_{\!\;\!j}\;\! \bigl[\: e^{\Delta (t-s)}
P_{\scriptscriptstyle \!\;\!H}(u\cdot \nabla u)(\cdot,s)\,\bigr]
\, ds. $$
\mbox{} \vspace{-0.170cm} \\
Therefore,
with \\
\mbox{} \vspace{-1.000cm} \\
$$\kappa \,=\,
\f{3}{2} \;\!\Bigl(\f{1}{r} - \f{1}{2} \Bigr) +\;\! \f{1}{2}\;\!,$$
\mbox{} \vspace{-0.270cm} \\
the following estimates hold: \\
\mbox{} \vspace{-0.975cm} \\
\beaa
\|\,D_{\!\;\!j} u(\cdot,t)\,\|
&\leq &
\|\,e^{\Delta \:\!t} D_{\!\;\!j}f\,\|
\;+\;
C \!\:\! \int_0^t (t-s)^{-\kappa} \,
\|\,P_{\scriptscriptstyle \!\;\!H}(u \cdot\!\;\! \nabla u)(\cdot,s)
\,\|_{\mbox{}_{\scriptstyle L^r}} \, ds \\
&\leq &
\|\,D_{\!\;\!j} f\,\|
\;+\;
C \!\:\! \int_0^t (t-s)^{-\kappa} \,
\|\,u \cdot\!\;\! \nabla u\;\!(\cdot,s)
\,\|_{\mbox{}_{\scriptstyle L^r}} \, ds \\
&\leq &
\|\,D_{\!\;\!j}\!\;\!f\,\|
\;+\;
C \!\:\! \int_0^t (t-s)^{-\kappa} \,
\|\,u(\cdot,s) \,\|_{\mbox{}_{\scriptstyle L^q}} \;\!
\|\,\cD u(\cdot,s)\,\| \; ds
\eeaa
\mbox{} \vspace{-0.370cm} \\
In the first estimate,
we have applied (\ref{eq.3.s4});
the second estimate follows from (\ref{eq.2.s4}),
and the third estimate uses H{\"o}lder's inequality.
We thus have \\
\mbox{} \vspace{-1.050cm} \\
\beaa
\|\,\cD u(\cdot,t)\,\|
\;\,\leq\;
\sum_{j\,=\,1}^{3} \!\;\!\|\,D_{\!\;\!j}\!\;\!f\,\|
\;+\;
C \!\:\! \int_0^t (t-s)^{-\kappa} \,
\|\,u(\cdot,s) \,\|_{\mbox{}_{\scriptstyle L^q}} \;\!
\|\,\cD u(\cdot,s)\,\| \: ds\:\!.
\eeaa
\mbox{} \vspace{-0.350cm} \\
Let us note that \\
\mbox{} \vspace{-1.050cm} \\
$$\kappa \,=\, \f{3}{2\:\!q} + \f{1}{2} \,< 1\:\!, $$
\mbox{} \vspace{-0.250cm} \\
since $\;\!q > 3$.
Therefore,
recalling Lemma \ref{Lemma.singular.Gronwall.1},
we see that
$(4.4a)$ implies $(4.4b)$, as claimed.
By Theorem \ref{theorem.Du.L2},
this gives
(\ref{eq.5.s4})
if $T_{\!\;\!f} \!\;\!$ is finite,
and the proof is now complete.
\hfill $\Box$ \\
%

\end{subsection}

%
%
\mbox{} \vspace{-1.250cm} \\
\begin{subsection}{An Integral Inequality for
$\;\!\|\,u(\cdot,t)\,\|_{\mbox{}_{\scriptstyle L^q}}$,
$ \;\! 3 < q \leq \infty $}
\mbox{} \vspace{-0.825cm} \\

We now show a simple
nonlinear integral inequality
for the scalar function
$\;\!\|\,u(\cdot,t)\,\|_{\mbox{}_{\scriptstyle L^q}} \!\;\!$
that gives some local control
on the growth of
$\;\!\|\,u(\cdot,t)\,\|_{\mbox{}_{\scriptstyle L^q}}\!\;\!$.
This local control
together with Theorem \ref{theorem.4.Lq}
imply a lower bound for
$\;\!\|\,u(\cdot,t)\,\|_{\mbox{}_{\scriptstyle L^q}}\!\:\!$
if $ \;\!T_{\!\;\!f} \!\;\!< \infty $,
cf.$\;$Lemma \ref{How.integral.inequalities.are.used}
below.

%
%
%
%
\begin{theorem}
\label{t.4.2}
\label{Integral.inequality.u.Lq}
Let $ \,3 < q \leq \infty $
and set \\
%
\mbox{} \vspace{-0.075cm} \\
\mbox{} \hspace{+5.750cm}
$ {\displaystyle
\kappa \:=\:
\f{3}{2q} \;\!+\;\! \f{1}{2} \,<\, 1\:\!
} $.
\mbox{} \hfill $(4.6a)$ \\
\mbox{} \vspace{-0.050cm} \\
Then,
there is a constant
$ \;\!C_{\!\;\!q} \!\,\!>\!\;\! 0 $
$($depending only on $\;\!q\:\!)$
such that,
for any $ \;\!0 \leq t_0 < T_{\!\;\!f}\!\;\!$,
we have \\
\mbox{} \vspace{-0.250cm} \\
\mbox{} \hspace{+0.400cm}
$ {\displaystyle
\|\,u(\cdot,t)\,\|_{\mbox{}_{\scriptstyle L^q}}
\:\!\leq\;
\|\,u(\cdot,t_0)\,\|_{\mbox{}_{\scriptstyle L^q}}
+\,
C_{\!\;\!q} \!
\int_{\mbox{}_{\mbox{\footnotesize $\!\:\!t_0$}}}^{\mbox{\footnotesize $\;\!t$}}
\! (t-s)^{- \kappa} \;\!
\|\,u(\cdot,s)\,\|_{\mbox{}_{\scriptstyle L^q}}^2 \, ds\:\!,
\quad \;\;
\forall \;\,
t_0 \!\;\! \leq t < T_{\!\;\!f} \!\;\!
} $.
\mbox{} \hfill $(4.6b)$ \\
\setcounter{equation}{6}
\end{theorem}
%
%
\mbox{} \vspace{-0.600cm} \\
{\bf Proof:}
In the case $ \;\! 3 < q < \infty $,
we use
the following argument
(adapted
from \cite{Giga1986}).
Let
$ {\displaystyle
\;\!
r \:\!=\;\! \f{q}{2}
\;\!
} $
and note that \\
\mbox{} \vspace{-0.950cm} \\
$$
\kappa \;\!=\,
\f{3}{2}\;\!\Bigl(\f{1}{r} - \f{1}{q}\Bigr)
\:\!+\;\! \f{1}{2}\;\!.
$$
\mbox{} \vspace{-0.250cm} \\
We also have
$ {\displaystyle
\;\!\|\:u_i \;\!u_j\;\!\|_{\mbox{}_{\scriptstyle L^r}}
\!\:\!\leq\;\! \|\:u\:\|_{\mbox{}_{\scriptstyle L^q}}^2\!\;\!
} $,
since $\;\! 2 \;\!r =\:\! q\:\!$.
Using
$(4.1a)$
and Duhamel's principle,
we get \\
\mbox{} \vspace{-0.750cm} \\
$$u(\cdot,t) \:=\; e^{\Delta (t-t_0)} u(\cdot, t_0)
\,-\:\! \sum_{j\,=\,1}^{3}
\int_{\mbox{}_{\mbox{\footnotesize $\!\:\!t_0$}}}^{\mbox{\footnotesize $\;\!t$}}
\!\! e^{\Delta (t-s)} \:\!
P_{\scriptscriptstyle \!\;\!H} \,\!\bigl[\;\!
D_j (\:\!u_j\;\!u\:\!)\:\!(\cdot,s) \,\bigr] \: ds,
\quad \;\;
t_0 \!\;\!\leq t < T_{\!\;\!f} \!\;\!,$$
\mbox{} \vspace{-0.150cm} \\
which gives,
by (\ref{eq.2.s4}) and (\ref{eq.3.s4})
above, \\
\mbox{} \vspace{-1.200cm} \\
\beaa
\|\,u(\cdot,t)\,\|_{\mbox{}_{\scriptstyle L^q}}
&\!\!\;\!\leq &
\|\,u(\cdot,t_0)\,\|_{\mbox{}_{\scriptstyle L^q}}
\,+\;
C_{\!\;\!q}\!\;\! \sum_{j\,=\,1}^{3}
\int_{\mbox{}_{\mbox{\footnotesize $\!\:\!t_0$}}}^{\mbox{\footnotesize $\;\!t$}}
\! (t-s)^{- \kappa} \,
\|\,P_{\scriptscriptstyle \!\;\!H} \:\!(u_j \;\!u)\:\!(\cdot,s)
\,\|_{\mbox{}_{\scriptstyle L^r}}
\, ds \\
&\!\!\;\!\leq &
\|\,u(\cdot,t_0)\,\|_{\mbox{}_{\scriptstyle L^q}}
\,+\;
C_{\!\;\!q} \!\;\!\sum_{j\,=\,1}^{3}
\int_{\mbox{}_{\mbox{\footnotesize $\!\:\!t_0$}}}^{\mbox{\footnotesize $\;\!t$}}
\! (t-s)^{- \kappa} \,
\|\,u_j \:\!u\;\!(\cdot,s) \,\|_{\mbox{}_{\scriptstyle L^r}}
\, ds \\
&\!\!\;\!\leq &
\|\,u(\cdot,t_0)\,\|_{\mbox{}_{\scriptstyle L^q}}
\,+\;
C_{\!\;\!q} \!
\int_{\mbox{}_{\mbox{\footnotesize $\!\:\!t_0$}}}^{\mbox{\footnotesize $\;\!t$}}
\! (t-s)^{- \kappa} \,
\|\, u(\cdot,s)\,\|_{\mbox{}_{\scriptstyle L^q}}^2
\, ds
\eeaa
\mbox{} \vspace{-0.250cm} \\
for all $\;\! t_0 \leq t < T_{\!\;\!f} \!\;\!$.
This shows the result
if $ \;\!3 < q < \infty $,
as claimed.
The proof in the
case $ \;\!q = \infty $
is due to Leray
and is developed in
Chapters 2 and 3
of \cite{Leray1934}.
\mbox{} \hfill $\Box$ \\
%
%
\mbox{} \vspace{-0.550cm} \\

The following lemma
shows how (4.6) is used
to yield
the fundamental lower bound (1.12)
for
$ \;\!\|\, u(\cdot,t) \,\|_{\mbox{}_{\scriptstyle L^{q}}} \!\;\!$,
$ \;\! 3 < q \leq \infty $,
in case of
finite--time blow--up.

%
%
%
%
\begin{lemma}
\label{How.integral.inequalities.are.used}
Let $\,\mbox{\em w} \in \!\;\!C^{0}([\;\!0, T_{\!\;\!f}[\;\!) $
be some positive function
such that we have,
for certain
$\;\! B > 0 $, $ \alpha > 1 $, $ \kappa < 1 $
constant, \\
\mbox{} \vspace{-0.175cm} \\
\mbox{} \hspace{+2.500cm}
$ {\displaystyle
\mbox{\em w}(t)
\;\leq\;
\mbox{\em w}(t_0)
\:+\,
B \!
\int_{\mbox{}_{\mbox{\footnotesize $\!\:\!t_0$}}}^{\mbox{\footnotesize $\;\!t$}}
\! (t-s)^{- \kappa} \:
\mbox{\em w}(s)^{\alpha}
\:ds
\qquad
\forall \;\,
t_0 \leq t < T_{\!\;\!f}
} $
\hfill $(4.7a)$ \\
\mbox{} \vspace{-0.080cm} \\
for each
$ \;\!0 \leq t_0 \!\;\! < T_{\!\;\!f} $.
If,
in addition,
$\: T_{\!\;\!f} < \infty \;\!$
and
$ {\displaystyle
\!\!\!
\sup_{\;\;0 \,\leq\, t \,<\, T_{\!\;\!f}}
\!\!\!\!\mbox{\em w}(t)
\,=\, \infty
} $,
$\,$then
it follows that \\
%
\mbox{} \vspace{-0.175cm} \\
\mbox{} \hspace{+2.250cm}
$ {\displaystyle
\mbox{\em w}(t)
\:>\;
\bigl(\:\! \alpha - 1 \bigr)
\cdot
\biggl(\;\!
\mbox{\small $ {\displaystyle
\frac{1 - \kappa}{\;\!B \, \alpha^{\mbox{}^{\scriptstyle \alpha}} } }$}
\;\!\biggr)^{\mbox{}^{\scriptstyle \!\!\!
\frac{\scriptstyle 1}{\scriptstyle \:\!\alpha - 1\:\!} }}
\hspace{-0.470cm} \cdot \hspace{+0.120cm}
\bigl(\;\! T_{\!\;\!f} -\;\!t
\;\!\bigr)^{\mbox{}^{\scriptstyle \!\!\!- \,
\frac{\scriptstyle 1 - \kappa}{\scriptstyle \:\!\alpha - 1\:\!} }}
\qquad
\forall \;\,
0 \leq t < T_{\!\;\!f}.
} $
\hfill $(4.7b)$ \\
\end{lemma}
\setcounter{equation}{7}
%
%
\mbox{} \vspace{-0.700cm} \\
{\bf Proof:}
Let $ \lambda > 1 $.
Given $ \;\!t_0 \in [\,0, T_{\!\;\!f}\;\![\;\!$ arbitrary,
by $(4.7a)$
we must have
$ \;\!{\tt w}(t) < \lambda \,{\tt w}(t_0) \;\!$
if $ \;\!t > t_0 \:\!$ is close to $ t_0 $.
In fact,
setting
$ \;\!\tau_{\ast} \!\:\!> 0 \;\!$
by \\
\mbox{} \vspace{-1.000cm} \\
$$
\tau_{\ast} \;\!:=\;\;\!
\min \,\Biggl\{\;\! T_{\!\;\!f}, \;
t_0 \:\!+\;\!
\biggl[ \frac{\;\!(1 - \kappa)\,(\lambda - 1)\;\!}
 {\lambda^{\mbox{}^{\scriptstyle \!\;\!\alpha}} B \:
 {\tt w}(t_0)^{\mbox{}^{\scriptstyle \alpha - 1}}}
\biggr]^{\mbox{}^{\scriptstyle \!\!\:\!
\frac{\scriptstyle 1}{\scriptstyle 1 \,-\,\kappa} }}
\Biggr\},
$$
\mbox{} \vspace{-0.200cm} \\
we have
$ {\displaystyle
\;\!
{\tt w}(t) < \lambda \,{\tt w}(t_0)
\;\!
} $
for all
$ \:\!t_0 \leq t < \tau_{\ast} $.
Because,
if this were false,
we could then find
$\;\!t_{1} \!\;\!\in \:]\,t_0, \:\!\tau_{\ast}\:\![\;\!$
such that
$ \;\! {\tt w}(t) < \lambda \, {\tt w}(t_0) \;\!$
for all
$ \:\!t_0 \leq t < t_{1} $,
while
$ \;\!{\tt w}(t_1) = \lambda \, {\tt w}(t_0) $.
This would give,
by $(4.7a)$
and the choice of
$ \tau_{\ast} \!\;\!$ above, \\
\mbox{} \vspace{-1.000cm} \\
\beaa
\lambda \, {\tt w}(t_0)
\;=\;\;\! {\tt w}(t_1)
&\!\!\;\!\leq &
\!{\tt w}(t_0) \,+\, B \!
\int_{\mbox{}_{\mbox{\footnotesize $\!\:\!t_0$}}}^{\mbox{\footnotesize $\;\!t_{1}$}}
\! (t_{1} - s)^{- \kappa} \:
{\tt w}(s)^{\alpha}
\:ds \\
&\!\!\;\!< &
\!{\tt w}(t_0) \,+\, B \!
\int_{\mbox{}_{\mbox{\footnotesize $\!\:\!t_0$}}}^{\mbox{\footnotesize $\;\!t_{1}$}}
\! (t_{1} - s)^{- \kappa} \:
\lambda^{\!\;\!\alpha} \,{\tt w}(t_0)^{\alpha}
\:ds \\
&\!\!\;\!< &
{\tt w}(t_0) \,+\, B \, \lambda^{\!\;\!\alpha} \,
{\tt w}(t_0)^{\mbox{}^{\scriptstyle \!\;\!\alpha}} \;\!
\frac{\,(\tau_{\ast} - t_0)^{1 - \kappa}}{1 - \kappa}
\;\:\leq\;\;\!
\lambda \,{\tt w}(t_0),
\eeaa
which could not be.
Hence,
we have
$ {\displaystyle
\;\!
{\tt w}(t) < \lambda \,{\tt w}(t_0)
\;\!
} $
for all
$ \:\!t_0 \leq t < \tau_{\ast} $,
as claimed,
and
in particular
$ \;\!{\tt w}\;\! $
is bounded
on $\,[\,t_0, \;\!\tau_{\ast}\:\![ $.
Since,
by assumption,
$ {\tt w} \;\!$ is unbounded
in $\,[\,t_0, \;\!T_{\!\;\!f}\:\![ $,
we must have
$ \;\!T_{\!\;\!f} > \tau_{\ast} $,
that is, \\
\mbox{} \vspace{-0.950cm} \\
$$
{\tt w}(t_0)
\;>\;
c(\lambda) \cdot
\bigl(\;\! T_{\!\;\!f} -\:\! t_0\:\!
\bigr)^{\mbox{}^{\scriptstyle \!\!\!-\,
\frac{\scriptstyle \:\!1 \:\!-\;\! \kappa\;\!}
     {\scriptstyle \alpha \;\!-\;\! 1} }}
\!\!\;\!,
\qquad
c(\lambda) \;=\;
\Bigl(\;\!
\mbox{\small $ {\displaystyle \frac{1 - \kappa}{B} }$}
\,\Bigr)^{\mbox{}^{\scriptstyle \!\!
\frac{\scriptstyle 1}{\scriptstyle \alpha \;\!-\;\! 1} }}
\hspace{-0.480cm} \cdot \hspace{+0.130cm}
\Bigl(\;\!
\mbox{\small $ {\displaystyle
\frac{\lambda - 1}{\lambda^{\mbox{}^{\scriptstyle \!\;\!\alpha}}}
} $}
\;\!\Bigr)^{\mbox{}^{\scriptstyle \!\!
\frac{\scriptstyle 1}{\scriptstyle \alpha \;\!-\;\! 1} }}
$$
\mbox{} \vspace{-0.250cm} \\
for
$ \;\!t_0 \!\;\!\in\!\;\! [\,0, T_{\!\;\!f} [\;\!$
arbitrary.
$\!\;\!$The largest value of $ c(\lambda) $
is obtained
by choosing
$ \;\!\lambda = \alpha/(\alpha - 1) $,
which yields the estimate $(4.7b)$.
\mbox{} \hfill $\Box$ \\
%
%
\mbox{} \vspace{-0.725cm} \\

From (4.6)
and Lemma \ref{How.integral.inequalities.are.used},
we get
for $\;\!3 < q \leq \infty \:\!$
the lower bound estimate \\
\mbox{} \vspace{-0.250cm} \\
\mbox{} \hspace{+3.000cm}
$ {\displaystyle
\|\,u(\cdot,t)\,\|_{\mbox{}_{\scriptstyle L^q}}
\;\!\geq\;
c_q \!\:\!\cdot
\bigl(\:\!T_{\!\;\!f} -\:\! t \:\!
\bigr)^{\mbox{}^{\scriptstyle \!\!\!-\,
\frac{\scriptstyle \;\!q \;\!-\;\!3\;\!}{\scriptstyle 2\;\!q} }}
\quad \;\;\;\,
\forall \;\,
0 \leq t < T_{\!\;\!f}
} $
\mbox{} \hfill $(4.8a)$ \\
\mbox{} \vspace{-0.100cm} \\
if
$\;\!T_{\!\;\!f}\!\;\!< \infty $,
where \\
\mbox{} \vspace{-0.150cm} \\
\mbox{} \hspace{+2.100cm}
$ {\displaystyle
c_q \;\!=\;
\frac{q - 3}{\;\!8\;\!q\,C_{\!\;\!q}}
} $
\mbox{} \hspace{+0.150cm}
if $\;3 < q < \infty $,
\mbox{} \hspace{+0.750cm}
$ {\displaystyle
c_{\infty} \;\!=\;
\frac{1}{\;8\;\!C_{\!\;\!\infty}}
} $
\mbox{} \hspace{+0.150cm}
if $\;q = \infty $,
\mbox{} \hfill $(4.8b)$ \\
\setcounter{equation}{8}
\mbox{} \vspace{+0.025cm} \\
with $ C_{\!\;\!q} \!\;\!> 0$
given in $(4.6b)$ above.
This proves Theorem \ref{thm.1.3}
of Section \ref{s1}
for $\:\!3 < q \leq \infty$.
(Another proof
for $ 3 \leq q < \infty $
is given in
Subsection 4.3.)
%
%
Using the Gagliardo inequality \\
\mbox{} \vspace{-0.550cm} \\
\be
\label{eq.9.s4}
\|\: {\sf u} \:\ni
\,\leq\,
K(q) \:
\|\: {\sf u} \:
\|_{\mbox{}_{\scriptstyle L^{2}}}
  ^{\mbox{}^{\scriptstyle 1 - \theta}}
\|\: \nabla \!\;\!{\sf u} \:
\|_{\mbox{}_{\scriptstyle L^{q}}}
  ^{\mbox{}^{\scriptstyle \:\!\theta}}
\!\;\!,
\quad \;\;\,
\theta \,=\,
\mbox{\normalsize $ {\displaystyle
\frac{3 \:\!q}{\;\! 5 \:\!q - 6\;\!} }$}
\quad \;\;\;\,
(\:\! 3 < q \leq \infty\:\!),
\ee
\mbox{} \vspace{-0.185cm} \\
which holds for arbitrary
$ \:\!{\sf u} \in L^{2}(\R^3) \cap W^{1,\mbox{\footnotesize $q$}}(\R^3) $,
we obtain,
from (2.1) and (4.8), \\
\mbox{} \vspace{-0.570cm} \\
\be
\label{eq.10.s4}
\|\, \cD u(\cdot,t) \,
\|_{\mbox{}_{\scriptstyle L^{q}}}
\;\!\geq\:
\hat{c}_{q} \,
\|\, f \,
\|_{\mbox{}_{\scriptstyle L^{2}}}
  ^{\scriptstyle \!\;\! -\,
    \frac{\scriptstyle \;\!2\:\!q \;\!-\;\! 6\;\!}
         {\scriptstyle 3 \:\!q} }
\mbox{\large (} T_{\!\;\!f} \!\;\!-\:\! t \:\!
\mbox{\large )}^{\scriptstyle \!\;\! -\,
                 \frac{\scriptstyle \;\!5\:\!q \;\!-\;\! 6\;\!}
                      {\scriptstyle 6 \:\!q} }
\quad \;\;
\forall \;\,
0 \leq t < T_{\!\;\!f}
\quad \;\;\;
\mbox{($\:\!$if $\,T_{\!\;\!} \!\;\!< \infty\:\!$)}
\ee
\mbox{} \vspace{-0.200cm} \\
for each
$ \;\!3 < q \leq \infty $,
and some constant
$ \hat{c}_{q} \!\:\!> 0 \:\!$
that depends only on $\:\!q$.
For
$ \:\!q = 3 $,
we can \linebreak
similarly obtain \\
\mbox{} \vspace{-0.670cm} \\
\be
\label{eq.11.s4}
\|\, \cD u(\cdot,t) \,
\|_{\mbox{}_{\scriptstyle L^{3}}}
\;\!\geq\:
c({\epsilon}) \:
\|\, f \,
\|_{\mbox{}_{\scriptstyle L^{2}}}
  ^{\mbox{}^{\scriptstyle \!\;\! -\,4\;\!\mbox{\footnotesize $\epsilon$}}}
\mbox{\large (} T_{\!\;\!f} \!\;\!-\:\! t \:\!
\mbox{\large )}^{\scriptstyle \!\;\! -\,
                 \frac{\scriptstyle 1}{\scriptstyle \;\!2\;\!}
                 \;\!+\;\! \mbox{\footnotesize $\epsilon$} }
\quad \;\;
\forall \;\,
0 \leq t < T_{\!\;\!f}
\quad \;\;\;
\mbox{($\:\!$if $\,T_{\!\;\!} \!\;\!< \infty\:\!$)}
\ee
\mbox{} \vspace{-0.250cm} \\
for each
$ \:\!0 < \epsilon \leq 1/2 $,
and some constant
$\:\!c(\epsilon) > 0 \:\!$
depending only on $\:\!\epsilon$,
using (2.1), (4.8)
and the 3D inequalities \\
\mbox{} \vspace{-0.600cm} \\
\be
\label{eq.12.s4}
\|\: {\sf u} \:\|_{\mbox{}_{\scriptstyle L^{\scriptstyle r}}}
\leq\,
K(r) \:
\|\: {\sf u} \:
\|_{\mbox{}_{\scriptstyle L^{2}}}
  ^{\mbox{}^{\scriptstyle 2/\mbox{\footnotesize $r$}}}
\|\: \nabla \!\;\!{\sf u} \:
\|_{\mbox{}_{\scriptstyle L^{3}}}
  ^{\mbox{}^{\scriptstyle \:\!1 \;\!-\;\! 2/\mbox{\footnotesize $r$}}}
\quad \;\;\;\,
(\:\! 3 \leq r < \infty\:\!).
\ee
\mbox{} \vspace{-0.200cm} \\
This completes the proof of (\ref{limit.Du.Lq}).
($\:\!$For Theorem \ref{theorem.Du.Lq},
see also (\ref{limit.u.L3}),
(\ref{SNG.Du.Lq}) and (\ref{eq.16.s4}).)
\mbox{} \vspace{-0.670cm} \\

\end{subsection}

%
%
\mbox{} \vspace{-1.250cm} \\
\begin{subsection}{A Differential Inequality for
$\;\!\|\,u(\cdot,t)\,\|_{\mbox{}_{\scriptstyle L^q}}$,
$ \;\! 3 < q < \infty $}
\mbox{} \vspace{-0.700cm} \\

We recall the fundamental estimate for the pressure
$ p(\cdot,t) $
obtained by the Calderon--Zygmund theory
applied to the Poisson equation
(\ref{Poisson.equation.pressure}), \\
\mbox{} \vspace{-0.650cm} \\
\be
\label{pressure.estimate}
\|\, p(\cdot,t) \,
\|_{\mbox{}_{\scriptstyle L^{\scriptstyle r}}}
\leq\,
C_{r} \,
\|\, u(\cdot,t) \,\|_{\mbox{}_{\scriptstyle L^{2\:\!{\scriptstyle r}}}}
                    ^{\mbox{}^{\scriptstyle \:\!2}}
\quad \;\;\;\,
\forall \;\, 0 \leq t < T_{\!\;\!f}
\quad \;\;\;\;\;\,
\mbox{$(\:\!1 < r < \infty \:\!)$},
\!\!\!
\ee
\mbox{} \vspace{-0.200cm} \\
see e.g.$\;$\cite{Galdi1994, LopesNussenzveigZheng1999}.
The basic result in this subsection is
the following differential inequality.

%
%
%
%
\begin{theorem}
\label{thm.4.3}
\label{theorem.differential.inequality.Du.Lq}
Let $\, T_{\!\;\!f}\!\;\!\leq \infty \:\!$
and $\;\! 3 < q < \infty $.
Then
there exists an absolute constant
$\:\!K_{q} $
$($depending only on $\:\!q\:\!)$
such that \\
\mbox{} \vspace{-0.760cm} \\
\be
\label{differential.inequality.Du.Lq}
\label{eq.14.s4}
\f{d}{d\:\!t\,} \:
\|\, u(\cdot,t) \,\|_{\mbox{}_{\scriptstyle L^{q}}}
                    ^{\mbox{}^{\scriptstyle \;\!q}}
\,\leq\:
K_{\!\;\!q}
\!\:\!\cdot\!\:\!
\Bigl(\,
\|\, u(\cdot,t) \,\|_{\mbox{}_{\scriptstyle L^{q}}}
                    ^{\mbox{}^{\scriptstyle \;\!q}}
\Bigr)^{\!\!\;\!\frac{\scriptstyle \;\!q \;\!-\;\!1\;\!}{q\;\!-\;\!3}}
\quad \;\;\;\,
\forall \;\,
0 \leq t < T_{\!\;\!f}.
\ee
\end{theorem}
%
%
\mbox{} \vspace{-0.200cm} \\
{\bf Proof:}
Given $\:\!\delta > 0 $,
let $\;\!L^{\prime}_{\delta}(\cdot) \;\!$
be a regularized sign function
(see e.g.$\;$\cite{KreissLorenz1989}, p.$\;$136),
and let \linebreak
$ \Phi_{\delta}({\tt u}) \!\:\!:= L_{\delta}({\tt u})^{\mbox{}^{\scriptstyle q}} \!$.
Multiplying the equation for $ u_{i}(\cdot,t)$
by
$ \Phi^{\prime}_{\delta}(u_{i}(\cdot,t)) $,
integrating on $ \R^{3}$
and letting
$ \:\!\delta \rightarrow 0 $,
we get \\
\mbox{} \vspace{-0.100cm} \\
\mbox{} \hspace{+2.000cm}
$ {\displaystyle
\f{d}{d\:\!t\,} \:
\|\, u_{i}(\cdot,t) \,\|_{\mbox{}_{\scriptstyle L^{q}}}
                        ^{\mbox{}^{\scriptstyle \:\!q}}
\:\!+\:
q \!\;\!\cdot\!\;\! (q - 1) \!
\int_{\mbox{}_{\mbox{\scriptsize $\!\;\!\R^3$}}} \!\!\!
|\, u_{i}(x,t) \,|^{\mbox{}^{\scriptstyle \:\!q - 2}} \;\!
|\, \nabla u_{i}\,|^{\mbox{}^{\scriptstyle \:\!2}}
\: dx
\;\leq
} $ \\
\mbox{} \vspace{+0.050cm} \\
\mbox{} \hspace{+5.200cm}
$ {\displaystyle
\leq \;
q \!\;\!\cdot\!\;\! (q - 1) \!
\int_{\mbox{}_{\mbox{\scriptsize $\!\;\!\R^3$}}} \!\!\!
|\, p(x,t) \,| \;
|\, u_{i}(x,t) \,|^{\mbox{}^{\scriptstyle \:\!q - 2}} \;\!
|\, \nabla u_{i}\,|
\: dx
} $ \\
\mbox{} \vspace{+0.120cm} \\
for all
$ 1 \leq i \leq 3 $,
$ \;\!0 \leq t < T_{\!\;\!f}$.
Using (\ref{pressure.estimate})
and H\"older's inequality,
we have \\
\mbox{} \vspace{-0.020cm} \\
\mbox{} \hspace{+2.000cm}
$ {\displaystyle
\int_{\mbox{}_{\mbox{\scriptsize $\!\;\!\R^3$}}} \!\!\!
|\, p(x,t) \,| \;
|\, u_{i}(x,t) \,|^{\mbox{}^{\scriptstyle \:\!q - 2}} \;\!
|\, \nabla u_{i}(x,t)\,|
\: dx
\;\leq
} $ \\
\mbox{} \vspace{-0.150cm} \\
\mbox{} \hfill
$ {\displaystyle
\leq\;
C(q) \:
\|\, u(\cdot,t) \,\|_{\mbox{}_{\scriptstyle L^{q+2}}}
                    ^{\mbox{}^{\scriptstyle \:\!2}}
\:\!
\|\, u_{i}(\cdot,t) \,
\|_{\mbox{}_{\scriptstyle L^{q+2}}}
  ^{\mbox{}^{\scriptstyle \frac{\scriptstyle q \;\!-\;\! 2}
                               {\scriptstyle 2} }}
\:\!
\biggl(
\int_{\mbox{}_{\mbox{\scriptsize $\!\;\!\R^3$}}} \!\!\!\:\!
|\, u_{i}(x,t) \,|^{\mbox{}^{\scriptstyle \:\!q - 2}} \;\!
|\, \nabla u_{i}\,|^{\mbox{}^{\scriptstyle \:\!2}}
\: dx
\:\!
\biggr)^{\!\!\;\!\frac{\scriptstyle 1}{\scriptstyle 2}}
} $ \\
\mbox{} \vspace{+0.050cm} \\
for each $\:\!i$,
and some constant $ \:\!C(q) > 0 \:\!$
that depends on the value of $\:\!q\:\!$ only.
In terms of
$ {\displaystyle
v(x,t) = \bigl(
v_{\mbox{}_{1}}(x,t), \;\! v_{\mbox{}_{2}}(x,t), \;\!
v_{\mbox{}_{3}}(x,t) \:\!\bigr)
} $
given by \\
\mbox{} \vspace{-0.150cm} \\
\mbox{} \hspace{+4.250cm}
$ {\displaystyle
v_{i}(x,t) \;\!:=\;
\bigl|\, u_{i}(x,t)
\,\bigr|^{\mbox{}^{\scriptstyle \!\;\!
\frac{\scriptstyle q}{\scriptstyle 2} }}
\!\!\;\!,
\qquad
1 \leq i \leq 3,
} $
\mbox{} \hfill $(4.15a)$ \\
\mbox{} \vspace{-0.050cm} \\
we therefore have \\
\mbox{} \vspace{-0.150cm} \\
\mbox{} \hspace{+1.000cm}
$ {\displaystyle
\f{d}{d\:\!t\,} \:
\|\, v_{i}(\cdot,t) \,\|_{\mbox{}_{\scriptstyle L^{2}}}
                        ^{\mbox{}^{\scriptstyle \:\!2}}
+\:
4 \,\Bigl(\;\! 1 - \frac{1}{q} \;\!\Bigr) \,
\|\, \nabla v_{i}(\cdot,t) \,
\|_{\mbox{}_{\scriptstyle L^{2}}}^{\mbox{}^{\scriptstyle \:\!2}}
\;\!\leq
} $ \\
\mbox{} \vspace{-0.200cm} \\
\mbox{} \hfill
$ {\displaystyle
\leq\;
2 \;\! q \,
\Bigl(\;\! 1 - \frac{1}{q} \;\!\Bigr) \,
C(q) \:
\|\, v(\cdot,t) \,
\|_{\mbox{}_{\scriptstyle L^{\beta}}}
  ^{\mbox{}^{\scriptstyle \:\!\frac{\scriptstyle 4}{\scriptstyle q} }}
\:\!
\|\, v_{i}(\cdot,t) \,
\|_{\mbox{}_{\scriptstyle L^{\beta}}}
  ^{\mbox{}^{\scriptstyle \!\:\!\frac{\scriptstyle q \;\!-\;\! 2}{\scriptstyle q} }}
\:\!
\|\, \nabla v_{i}(\cdot,t) \,
\|_{\mbox{}_{\scriptstyle L^{2}}}
} $ \\
\mbox{} \vspace{-0.000cm} \\
for each $\:\!i$,
where
$ {\displaystyle
\;\!
\mbox{\small $\beta$} \;\!=\;\!
\mbox{\small 2} + \mbox{\small 4}/q
} $.
Using the inequality \\
\mbox{} \vspace{-0.720cm} \\
$$
\|\; {\sf v} \;\|_{\mbox{}_{\scriptstyle L^{\beta}}}
\;\!\leq\;
K\!\;\!(\mbox{\small $\beta$}) \:
\|\; {\sf v} \;
\|_{\mbox{}_{\scriptstyle L^{2}}}
  ^{\mbox{}^{\scriptstyle \frac{\scriptstyle q \;\!-\;\! 1}
                               {\scriptstyle q \;\!+\;\! 2} }}
\:\!
\|\; \nabla {\sf v} \;
\|_{\mbox{}_{\scriptstyle L^{2}}}
  ^{\mbox{}^{\scriptstyle \frac{\scriptstyle 3}
                               {\scriptstyle q \;\!+\;\! 2} }}
\qquad
\forall \;\,
{\sf v} \in H^{1}(\R^3),
$$
\mbox{} \vspace{-0.230cm} \\
where
the constant
$ K\!\;\!(\mbox{\small $\beta$}) > 0 \;\!$
depends only on $ \mbox{\small $\beta$} $,
and
summing on $ i = 1, 2, 3 $,
we obtain \\
\mbox{} \vspace{-0.150cm} \\
\mbox{} \hspace{+0.350cm}
$ {\displaystyle
\f{d}{d\:\!t\,} \:
\|\, v(\cdot,t) \,\|_{\mbox{}_{\scriptstyle L^{2}}}
                    ^{\mbox{}^{\scriptstyle \:\!2}}
+\:
4 \,\Bigl(\;\! 1 - \frac{1}{q} \;\!\Bigr) \,
\|\, \cD v(\cdot,t) \,
\|_{\mbox{}_{\scriptstyle L^{2}}}^{\mbox{}^{\scriptstyle \:\!2}}
\,
\leq\;
C_{\!\;\!q} \:
\|\, v(\cdot,t) \,
\|_{\mbox{}_{\scriptstyle L^{2}}}
  ^{\mbox{}^{\scriptstyle \frac{\scriptstyle q \;\!-\;\!1}{\scriptstyle q} }}
\|\, \cD v(\cdot,t) \,
\|_{\mbox{}_{\scriptstyle L^{2}}}
  ^{\mbox{}^{\scriptstyle \frac{\scriptstyle q \;\!+\;\!3}{\scriptstyle q} }}
} $
\mbox{} \hfill $(4.15b)$ \\
\setcounter{equation}{15}
\mbox{} \vspace{+0.025cm} \\
for all $ \;\!0 \leq t < T_{\!\;\!f} $,
and some constant
$ C_{\!\;\!q} \!\:\!> 0 $
that depends on $\:\!q\:\!$ only.
This gives \\
\mbox{} \vspace{-0.700cm} \\
$$
\f{d}{d\:\!t\,} \:
\|\, v(\cdot,t) \,\|_{\mbox{}_{\scriptstyle L^{2}}}
                    ^{\mbox{}^{\scriptstyle \:\!2}}
\,\leq\:
K_{\!\;\!q} \!\:\!\cdot\!\:\!
\Bigl(\,
\|\, v(\cdot,t) \,\|_{\mbox{}_{\scriptstyle L^{2}}}
                    ^{\mbox{}^{\scriptstyle \:\!2}}
\Bigr)^{\!\!\;\!\frac{\scriptstyle \;\!q \;\!-\;\!1\;\!}{q\;\!-\;\!3}}
\quad \;\;\;\,
\forall \;\,
0 \leq t < T_{\!\;\!f}
$$
\mbox{} \vspace{-0.250cm} \\
for some constant
$ K_{\!\;\!q} \!\;\!> 0 \:\!$
depending on $\:\!q\:\!$ only,
which is equivalent to
(\ref{differential.inequality.Du.Lq}).
\mbox{} \hfill $\Box$ \\
%
%
\mbox{} \vspace{-0.650cm} \\

It follows from
the previous proof
that the estimate (4.15$b$)
is valid more generally
for any $\:\! 2 < q < \infty $.
Taking
$\:\!q = 3 $,
it gives
that
$ {\displaystyle
\;\!
d/dt \,
\|\, v(\cdot,t) \,
\|_{L^{2}}^{\:\!2}
\!\;\!<\:\! 0
\:\!
} $
if
$ {\displaystyle
\;\!
\|\, v(\cdot,t) \,
\|_{\mbox{}_{\scriptstyle L^{2}}}
\!\:\!
} $
is appropriately small;
since
$ {\displaystyle
\;\!
\|\, v(\cdot,t) \,
\|_{\mbox{}_{\scriptstyle L^{2}}}^{\:\!2}
\!\:\!=\;\!
\|\, u(\cdot,t) \,
\|_{\mbox{}_{\scriptstyle L^{3}}}^{\:\!3}
\!\:\!
} $
in this case,
cf.$\;$(4.15$a$),
we con\-clude that
$ {\displaystyle
\|\, u(\cdot,t) \,
\|_{\mbox{}_{\scriptstyle L^{3}}}
\!\:\!
} $
is monotonically decreasing
in $ t $
when
$ {\displaystyle
\:\!
\|\, u(\cdot,0) \,
\|_{\mbox{}_{\scriptstyle L^{3}}}
\!\:\!
} $
is sufficiently small,
i.e., \\
\mbox{} \vspace{-1.000cm} \\
\be
\label{eq.16.s4}
\|\, u(\cdot,0) \,\|_{\mbox{}_{\scriptstyle L^{3}}}
\:\!<\;
\eta_{\mbox{}_{\scriptstyle 3}}
\;\;\;
\Longrightarrow
\;\;\;
T_{\!\;\!f} =\;\! \infty
\ee
\mbox{} \vspace{-0.200cm} \\
for some
absolute value
$ \:\!\eta_{\mbox{}_{\scriptstyle 3}} \!\:\!> 0 $.
This shows (1.12),
(\ref{global.existence.Lq})
for $ \;\! q = 3 $,
and also,
using (\ref{SNG.Du.Lq}), \linebreak
the bound
(\ref{lowerbound.Du.Lq})
for $ \:\!q = 3/2 $,
thus completing the proof
of Theorems \ref{thm.1.2}
and \ref{thm.1.3}
above.
It also implies,
by~(\ref{limit.u.L3})
and
Gronwall's lemma,
that
$ {\displaystyle
\;\!
\|\, \cD v(\cdot,t) \,
\|_{L^{2}}^{\:\!2}
\!\;\!
} $
cannot be integrable
on $ \;\![\,0, \;\!T_{\!\;\!f}\;\!] \;\!$
if $\;\!T_{\!\;\!f} \!\;\!< \infty $,
or,
in terms of $ \:\!u(\cdot,t) $,
that
we have \\
\mbox{} \vspace{-0.625cm} \\
\be
\label{eq.17.s4}
\sum_{i\,=\,1}^{3} \,
\int_0^{T_{\!\;\!f}} \!\!
\int_{\mbox{}_{\mbox{\scriptsize $\!\;\!\R^3$}}} \!\!\!\:\!
|\, u_{i}(x,t) \,| \;\:\!
|\, \nabla u_{i}(x,t) \,|^{\:\!2}
\; dx \, dt
\;\;\!=\;\;\!
\infty
\qquad
\mbox{$(\:\!$if $\;\!T_{\!\;\!f} \!\;\!< \infty\:\!)$}.
\ee
\mbox{} \vspace{-0.000cm} \\
On the other hand,
taking $ \:\!q > 3 $
in (4.15),
we reobtain the fundamental estimate (1.12),
in view of
Lemma \ref{How.differential.inequalities.are.used}.
Another important consequence
is the following.
Using that \\
\mbox{} \vspace{-0.570cm} \\
\be
\label{eq.18.s4}
\mbox{} \;\;
\|\; {\sf v} \;\|_{\mbox{}_{\scriptstyle L^{2}}}
\;\!\leq\;
K\!\;\!(q)\,
\|\; {\sf v} \;\|_{\mbox{}_{\scriptstyle L^{4/q}}}
                 ^{\mbox{}^{\scriptstyle \:\! 1 \;\!-\;\! \delta}}
\;\!
\|\; \nabla {\sf v} \;
\|_{\mbox{}_{\scriptstyle L^{2}}}^{\mbox{}^{\scriptstyle \;\! \delta}}
\!\:\!,
\qquad
\delta \,=\,
\frac{\;\!3 \:\!q - 6\;\!}{3\:\!q - 2}
\qquad
\,
(\:\!2 \leq q < \infty\:\!),
\ee
\mbox{} \vspace{-0.200cm} \\
we have,
for $ \:\! q > 3 \:\!$: \\
\mbox{} \vspace{-0.220cm} \\
\mbox{} \hspace{+1.750cm}
$ {\displaystyle
\|\, v(\cdot,t) \,
\|_{\mbox{}_{\scriptstyle L^{2}}}
  ^{\mbox{}^{\scriptstyle
    \frac{\scriptstyle \;\!q \;\!-\;\!1\;\!}{\scriptstyle q} }}
\:\!
\|\, \cD v(\cdot,t) \,
\|_{\mbox{}_{\scriptstyle L^{2}}}
  ^{\mbox{}^{\scriptstyle
    \frac{\scriptstyle \;\!q \;\!+\;\!3\;\!}{\scriptstyle q} }}
\;\!=
} $ \\
\mbox{} \vspace{-0.170cm} \\
\mbox{} \hspace{+4.250cm}
$ {\displaystyle
=\;
\|\, v(\cdot,t) \,
\|_{\mbox{}_{\scriptstyle L^{2}}}
  ^{\mbox{}^{\scriptstyle
    \frac{\scriptstyle 2}{\scriptstyle \;\!3 \:\!q \;\!-\;\!6\;\!} }}
\:\!
\|\, v(\cdot,t) \,
\|_{\mbox{}_{\scriptstyle L^{2}}}
  ^{\mbox{}^{\scriptstyle
    \frac{\scriptstyle \;\!3\:\!q \;\!-\;\!2\;\!}
         {\scriptstyle \;\!3 \:\!q \;\!-\;\!6\;\!}
    \;\!
    \bigl( 1 \;\!-\;\! \frac{\scriptstyle \;\!3\;\!}{\scriptstyle q} \bigr) }}
\|\, \cD v(\cdot,t) \,
\|_{\mbox{}_{\scriptstyle L^{2}}}
  ^{\mbox{}^{\scriptstyle
    \frac{\scriptstyle \;\!q \;\!+\;\!3\;\!}{\scriptstyle q} }}
} $ \\
\mbox{} \vspace{-0.170cm} \\
\mbox{} \hspace{+4.250cm}
$ {\displaystyle
\leq\;
C(q) \;
\|\, v(\cdot,t) \,
\|_{\mbox{}_{\scriptstyle L^{4/q}}}
  ^{\mbox{}^{\scriptstyle
    \frac{\scriptstyle \;\!4\;\!}{\scriptstyle q} \;\!
    \frac{\scriptstyle \;\!q \;\!-\;\! 3 \;\!}
         {\scriptstyle \;\!3 \:\!q \;\!-\;\!6\;\!} }}
\|\, v(\cdot,t) \,
\|_{\mbox{}_{\scriptstyle L^{2}}}
  ^{\mbox{}^{\scriptstyle
    \frac{\scriptstyle 2}{\scriptstyle \;\!3 \:\!q \;\!-\;\!6\;\!} }}
\|\, \cD v(\cdot,t) \,
\|_{\mbox{}_{\scriptstyle L^{2}}}
  ^{\mbox{}^{\scriptstyle \:\!2 }}
} $ \\
\mbox{} \vspace{-0.130cm} \\
\mbox{} \hspace{+4.250cm}
$ {\displaystyle
=\;
C(q) \;
\|\, u(\cdot,t) \,
\|_{\mbox{}_{\scriptstyle L^{2}}}
  ^{\mbox{}^{\scriptstyle
    \frac{\scriptstyle \;\!2 \:\!q \;\!-\;\! 6 \;\!}
         {\scriptstyle \;\!3 \:\!q \;\!-\;\! 6\;\!} }}
\|\, u(\cdot,t) \,
\|_{\mbox{}_{\scriptstyle L^{q}}}
  ^{\mbox{}^{\scriptstyle
    \frac{\scriptstyle q}{\scriptstyle \;\!3 \:\!q \;\!-\;\!6\;\!} }}
\|\, \cD v(\cdot,t) \,
\|_{\mbox{}_{\scriptstyle L^{2}}}
  ^{\mbox{}^{\scriptstyle \:\!2 }}
} $ \\
\mbox{} \vspace{-0.000cm} \\
by (4.15$a$) and (\ref{eq.18.s4})
above.
As
$ {\displaystyle
\:\!
\|\, u(\cdot,t) \,\|_{\mbox{}_{\scriptstyle L^{2}}}
\!\;\!
} $
never increases,
this gives,
because of (4.15$b$),
that
$ {\displaystyle
\:\!
\|\, v(\cdot,t) \,\|_{\mbox{}_{\scriptstyle L^{2}}}
\:\!
\bigl(\:\!=\;\!
\|\, u(\cdot,t) \,\|_{\mbox{}_{\scriptstyle L^{q}}}
\!\;\!
\bigr)
} $
is monotonically decreasing
in time
whenever we have \\
\mbox{} \vspace{-0.750cm} \\
$$
\|\, u(\cdot,0) \,
\|_{\mbox{}_{\scriptstyle L^{2}}}
  ^{\mbox{}^{\scriptstyle
    \frac{\scriptstyle \;\!2 \:\!q \;\!-\;\! 6 \;\!}
         {\scriptstyle \;\!3 \:\!q \;\!-\;\! 6\;\!} }}
\|\, u(\cdot,0) \,
\|_{\mbox{}_{\scriptstyle L^{q}}}
  ^{\mbox{}^{\scriptstyle
    \frac{\scriptstyle q}{\scriptstyle \;\!3 \:\!q \;\!-\;\!6\;\!} }}
<\;
\eta_{\mbox{}_{\scriptstyle q}}
$$
\mbox{} \vspace{-0.270cm} \\
for some value
$ \:\!\eta_{\mbox{}_{\scriptstyle q}} \!\:\!> 0 \:\!$
appropriately small
(depending only on $\:\!q\:\!$).
Together with (\ref{eq.16.s4}),
this shows (\ref{global.existence.Lq}),
Theorem \ref{thm.1.4},
for $ 3 \leq q < \infty $.
$\!$The proof for $ q = \infty \:\!$
is given in~\cite{Leray1934}. \\
\mbox{} \vspace{-0.950cm} \\

We finish this Section
with a few last remarks.
Using (4.8)
and the 3D inequality \\
\mbox{} \vspace{-0.625cm} \\
\be
\label{SNG.D2u.Lq}
\label{eq.19.s4}
\|\: \mbox{\sf u} \:\|_{\mbox{}_{\scriptstyle L^{r(q)}}}
\,\!\leq\,
K_{\!\:\!q} \:
\|\, {\cal D}^{2} \mbox{\sf u}
\:\|_{\mbox{}_{\scriptstyle L^{q}}}
\!\;\!,
\quad
\;\;
r(q) \;\!=\;\! \frac{3 \:\!q}{\;\!3 - 2 \:\!q\,}
\quad
\;\;\;\;
\biggl(\;\!
1 \leq q < \mbox{\normalsize $ {\displaystyle \frac{\;\!3\;\!}{2} }$}
\;\!\biggr)
\!\;\!,
\ee
\mbox{} \vspace{-0.200cm} \\
where
$ K_{\!\:\!q} \!\;\!> 0 \;\!$
depends only on $\:\!q $,
we obtain
the lower bound estimate \\
\mbox{} \vspace{-0.825cm} \\
\be
\label{eq.20.s4}
\|\, \cD^2 u(\cdot,t) \,
\|_{\mbox{}_{\scriptstyle L^{q}}}
\;\!\geq\:
c_{\mbox{}_{\scriptstyle q}} \;\!
\mbox{\large (} T_{\!\;\!f} \!\;\!-\:\! t \:\!
\mbox{\large )}^{\mbox{}^{\scriptstyle \!\!\!\;\! -\,
                 \frac{\scriptstyle 3}{\scriptstyle \;\!2\;\!}
                 \;\!
                 \frac{\scriptstyle \;\! q \;\!-\;\!1\;\!}{\scriptstyle q} }}
\quad \;
\forall \;\,
0 \leq t < T_{\!\;\!f}
\qquad
\mbox{($\;\!$if $\,T_{\!\;\!} \!\;\!< \infty\:\!$)}
\ee
\mbox{} \vspace{-0.270cm} \\
for each
$\:\!1 \leq q < 3/2 $,
and some constant
$ c_{\mbox{}_{\scriptstyle q}} \!\:\!> 0 \:\!$
that depends only on $\:\!q$.
The estimate (\ref{eq.20.s4})
has been recently shown
in \cite{RobinsonSadowskiSilva2012}
to hold for $ q = 2 $ as well,
but its validity for arbitrary $ q \geq 3/2 $
seems to be still open.
%
%
The general fact that
the norms
$ {\displaystyle
\|\, \cD^{n} u(\cdot,t) \,\|_{\mbox{}_{\scriptstyle L^{q}}}
\!\;\!
} $,
$ 1 \leq q \leq \infty $,
$ n \geq 2 $,
do all blow up
as $ \;\!t \mbox{\small $\nearrow$}\;\! T_{\!\;\!f} $
in case $\;\!T_{\!\;\!f} \!\;\!< \infty \;\!$
is a direct consequence
of (\ref{limit.u.L3}), (4.8)
and the family of
3D Gagliardo inequalities
given by \\
\mbox{} \vspace{-0.770cm} \\
$$
\|\; {\sf u} \;\|_{\mbox{}_{\scriptstyle L^{q}}}
\leq\,
K\!\;\!(q,r) \:
\|\; {\sf u} \;
\|_{\mbox{}_{\scriptstyle L^{2}}}
  ^{\mbox{}^{\scriptstyle \:\! 1 - \theta}}
\:\!
\|\; \cD^n {\sf u} \;
\|_{\mbox{}_{\scriptstyle L^{r}}}
  ^{\mbox{}^{\scriptstyle \:\! \theta}}
\!\;\!,
\quad
%
%
\theta \,=\, \frac{\;\!\mbox{\small $ 1/2 \;\!-\;\! 1/q $}\;\!}
                  {\;\!\mbox{\small $ 1/2 \;\!+\;\! n/3 \;\!-\;\! 1/r $}\;\!},
\quad
r \geq
\max \, \Bigl\{\;\! \mbox{\small $1$}, \,
\frac{\;\!\mbox{\small $3 \:\!q$}\;\!}
     {\;\!\mbox{\small $n \:\!q \;\!+\;\! 3$}\;\!}
\,\Bigr\}
$$
%
%
for $\:\! n \geq \mbox{\normalsize $2$} $,
$ \:\!\mbox{\normalsize $3$} \leq q \leq \mbox{\normalsize $\infty$} \:\!$
arbitrary,
provided that
$ {\displaystyle
(n, q, r) \ne
(\mbox{\normalsize $2$}, \mbox{\normalsize $\infty$}, \mbox{\normalsize $3/2$})
} $,
$ {\displaystyle
(n, q, r) \ne
(\mbox{\normalsize $3$}, \mbox{\normalsize $\infty$}, \mbox{\normalsize $1$} )
} $.

\end{subsection}

\end{section}

\mbox{} \vspace{-0.950cm} \\
%
%
%
\begin{section}{Comparison of blow--up functions}
\label{s5}
\mbox{} \vspace{-1.030cm} \\
\setcounter{equation}{0}

Let $ \;\!3 \leq q < r \leq \infty \:\!$
and assume that
$ T_{\!\;\!f} < \infty $.
Theorem \ref{thm.1.3}
yields the lower bounds \\
\mbox{} \vspace{-0.630cm} \\
$$
\|\,u(\cdot,t)\,\|_{\mbox{}_{\scriptstyle L^q}}
\geq\:
c_q \;\!
(T_{\!\;\!f} -\;\!t)^{- \;\!\kappa(q)}
\!\;\!,
\quad \;\;
\|\,u(\cdot,t)\,\|_{\mbox{}_{\scriptstyle L^r}}
\geq\;
c_r \;\!
(T_{\!\;\!f} -\;\!t)^{- \;\!\kappa(r)}
$$
\mbox{} \vspace{-0.230cm} \\
with positive constants
$\,\!c_q, \:\!c_r$
and \\
\mbox{} \vspace{-0.690cm} \\
$$
0 \:\!\leq\:\! \kappa(q)
\:\!=\,
\f{q-3}{2\:\!q}
\;\!<\:\!
\kappa(r)
\:\!=\,
\f{r-3}{2r}\;\!.
$$
\mbox{} \vspace{-0.250cm} \\
Thus,
the lower bound for the $L^r\!\;\!$ norm
blows up faster than the
lower bound for the $L^q\!\;\!$ norm.
This suggests that \\
\mbox{} \vspace{-0.590cm} \\
\be
\label{eq.1.s5}
\f{\;\!\|\,u(\cdot,t)\,\|_{\mbox{}_{\scriptstyle L^r}}}
  {\;\!\|\,u(\cdot,t)\,\|_{\mbox{}_{\scriptstyle L^q}}}
\,\to \infty
\qas
t\:\! \to T_{\!\;\!f}
\qif
3 \leq q < r \leq \infty.
\ee
\mbox{} \vspace{-0.150cm} \\
A precise result can be obtained
by using boundedness of the $L^2\!\;\!$ norm
and interpolation:
defining
$\;\!0 < \lam < 1\;\!$
by \\
\mbox{} \vspace{-0.900cm} \\
$$
\f{1}{q}
\:=\:
\f{\lam}{2} \,+\, \f{1-\lam}{r}\;\!
$$
\mbox{} \vspace{-0.170cm} \\
and
recalling the interpolation estimate
$ {\displaystyle
\;\!
\|\,u(\cdot,t)\,\|_{\mbox{}_{\scriptstyle L^q}}
\leq\;\!
\|\,u(\cdot,t)\,\|^\lam \,
\|\,u(\cdot,t)\,\|_{\mbox{}_{\scriptstyle L^r}}^{1-\lam}
\!\:\!
} $,
$\;\!$one gets \\
\mbox{} \vspace{-0.550cm} \\
\be
\label{eq.5.lb}
\label{eq.2.s5}
\|\,u(\cdot,t)\,\|_{\mbox{}_{\scriptstyle L^r}}^\lam
\leq\,
\|\,f\,\|^{\lambda} \,
\f{\;\!\|\,u(\cdot,t)\,\|_{\mbox{}_{\scriptstyle L^r}}}
  {\;\!\|\,u(\cdot,t)\,\|_{\mbox{}_{\scriptstyle L^q}}}
\;\!,
\qquad
\lambda \,=\, \frac{\;\!\mbox{\small $ 1/q \;\!-\;\! 1/r $}\;\!}
                  {\;\!\mbox{\small $ 1/2 \;\!-\;\! 1/r $}\;\!}.
\ee
\mbox{} \vspace{-0.050cm} \\
Using the lower bound on the blow--up of
$\;\!\|\,u(\cdot,t)\,\|_{\mbox{}_{\scriptstyle L^r}}\!\;\!$
provided by Theorem \ref{thm.1.3},
one then obtains (\ref{eq.1.s5})
with
an algebraic lower bound,
as described next.

%
%
%
%
\begin{theorem}
Let $\;\!3 \leq q < r \leq \infty $,
and assume that $ \;\!T_{\!\;\!f} < \infty $.
Then
there is a constant $\:\!c(f) = c(f;q,r) > 0$,
depending on $\;\!q, \;\!r $
and the initial state $f\!\:\!$,
such that \\
\mbox{} \vspace{-0.100cm} \\
\mbox{} \hspace{+3.000cm}
$ {\displaystyle
\f{\,\|\,u(\cdot,t)\,\|_{\mbox{}_{\scriptstyle L^r}}}
    {\,\|\,u(\cdot,t)\,\|_{\mbox{}_{\scriptstyle L^q}}}
\,\geq\;
c(f) \cdot
\bigl(\:\!T_{\!\;\!f} -\;\! t\;\!
\bigr)^{\mbox{}^{\scriptstyle \!\!-\,\gamma}}
\qquad
\forall \;\,
0 \leq t < T_{\!\;\!f},
} $
\mbox{} \hfill $(5.3a)$ \\
\mbox{} \vspace{-0.000cm} \\
where \\
\mbox{} \vspace{-0.425cm} \\
\mbox{} \hspace{+5.500cm}
$ {\displaystyle
\gamma
\,=\:
\frac{\;\!r - 3\;\!}{\;\!r - 2\;\!} \cdot
\frac{\;\!r - q\;\!}{q \:\!r}
\;\!>\;\! 0\:\!.
} $
\mbox{} \hfill $(5.3b)$ \\
\setcounter{equation}{3}
\end{theorem}
%
%

In a similar way,
using the 3D inequality \\
\mbox{} \vspace{-0.600cm} \\
\be
\label{eq.4.s5}
\|\; {\sf u} \;\|_{\mbox{}_{\scriptstyle L^{q}}}
\leq\,
K\!\;\!(q) \:
\|\; {\sf u} \;\|^{\scriptstyle \:\!1 - \theta}
\;\!
\|\; \cD {\sf u} \;\|^{\scriptstyle \:\!\theta}
\!,
\qquad
\theta \,=\,
\frac{\;\!3\;\!}{2} \cdot
\frac{\;\!q - 2\;\!}{q}
\qquad
(\:\!2 \leq q \leq 6 \:\!)
\ee
\mbox{} \vspace{-0.200cm} \\
and assuming
$ \:\!T_{\!\;\!f} \!\;\!< \infty $,
one obtains \\
\mbox{} \vspace{-0.550cm} \\
\be
\f{\;\!\|\,\cD u(\cdot,t)\,\|\;\!}
    {\;\;\,\|\,u(\cdot,t)\,\|_{\mbox{}_{\scriptstyle L^q}}}
\,\geq\;
\tilde{c}(f) \cdot
\bigl(\:\!T_{\!\;\!f} -\;\! t \;\!
\bigr)^{\mbox{}^{\scriptstyle \!\! -\,\tilde{\gamma}}}
\!,
\qquad
\tilde{\gamma} \,=\,
\frac{\;\!6 - q\;\!}{8 \:\! q}
\;\!\geq\;\! 0
\ee
\mbox{} \vspace{-0.100cm} \\
for all
$\:\!0 \leq t < T_{\!\;\!f} \!\;\!$,
and every
$ \:\!2 \leq q \leq 6 $,
where the constant
$\:\!\tilde{c}(f) \:\!$
depends on
$\:\!q\:\!$
and
$ \:\!\|\,f\,\| $. \linebreak
\mbox{} \vspace{-0.850cm} \\
%
%

In the remaining part of this section
we compare the growth of
$ {\displaystyle
\;\!
\|\,u(\cdot,t)\,\|_{\mbox{}_{\scriptstyle L^q}}^{\:\!q}
\:\!
\|\,u(\cdot,t)\,\|_{\mbox{}_{\scriptstyle \infty}}^{\:\!2}
} $,
for $3 < q < \infty $,
with the growth of
$ {\displaystyle
\;\!
\|\,u(\cdot,t)\,\|_{\mbox{}_{\scriptstyle L^{3}}}
\:\!
\|\,u(\cdot,t)\,\|_{\mbox{}_{\scriptstyle \infty}}^{\:\!q}
} $,
as $\:\! t \mbox{\small $\nearrow$}\;\! T_{\!\;\!f}$.
Setting $ r = \infty \:\!$
and $\:\! \lam = 2/q \:\!$
in (\ref{eq.5.lb}) above,
we have \\
\mbox{} \vspace{-0.700cm} \\
\be
\label{eq.r.infty}
\label{eq.6.s5}
\|\,u(\cdot,t)\,\ni^{\mbox{}^{\scriptstyle 2/q}}
\leq\:
\|\,f\,\|^{\mbox{}^{\scriptstyle 2/q}} \;\!
\f{\;\!\|\,u(\cdot,t)\,\ni}
{\;\|\,u(\cdot,t)\,\|_{\mbox{}_{\scriptstyle L^q}}}
\qquad
\forall \;\,
0 \leq t < T_{\!\;\!f}.
\ee
\mbox{} \vspace{-0.150cm} \\
Together with (\ref{limit.u.L3})
and (5.3),
one obtains the following theorem.

%
%
%
%
\begin{theorem}
Let $ \;\!3 < q < \infty $,
and assume that $\;\! T_{\!\;\!f} < \infty $.
Then
we have \\
\mbox{} \vspace{-0.670cm} \\
\be
\label{eq.7.s5}
\lim_{\mbox{\footnotesize $t$} \:\!\mbox{\scriptsize $\nearrow$}\: T_{\!\;\!f}}
\f{\;\!\|\,u(\cdot,t)\,\ni^{\mbox{}^{\scriptstyle q}}}
  {\;\!\|\,u(\cdot,t)\,\|_{\mbox{}_{\scriptstyle L^q}}^{\mbox{}^{\scriptstyle \:\!q}}}
\cdot
\f{\;\!\|\,u(\cdot,t)\,\|_{\mbox{}_{\scriptstyle L^3}}}
  {\;\!\|\,u(\cdot,t)\,\ni^{\mbox{}^{\scriptstyle 2}}}
\;=\: \infty\:\!.
\ee
\end{theorem}
%
%
\mbox{} \vspace{-0.150cm} \\
{\sc Remark:}
Since the proof of (\ref{limit.u.L3})
in \cite{Seregin2012} is very involved,
we give here a direct elementary argument
for the weaker statement \\
\mbox{} \vspace{-0.720cm} \\
\be
\label{eq.5.sup}
\label{eq.8.s5}
\sup_{0 \,\leq\, \mbox{\footnotesize $t$} \,<\, T_{\!\;\!f}}
\!\:\!
\f{\;\!\|\,u(\cdot,t)\,\ni^{\mbox{}^{\scriptstyle q}}}
  {\;\!\|\,u(\cdot,t)\,\|_{\mbox{}_{\scriptstyle L^q}}
                         ^{\mbox{}^{\scriptstyle \:\!q}}}
\cdot
\f{\;\!\|\,u(\cdot,t)\,\|_{\mbox{}_{\scriptstyle L^3}}}
  {\;\!\|\,u(\cdot,t)\,\ni^{\mbox{}^{\scriptstyle 2}}}
\;=\: \infty\:\!.
\ee
\mbox{} \vspace{-0.100cm} \\
Denoting by
$\langle u, v \rangle \!\;\!= \sum_i u_i\;\!v_i$
the Euclidean inner product in $\R^3\!\:\!$,
we have \\
\mbox{} \vspace{-1.050cm} \\
\beaa
\f{1}{q}\: \f{d}{dt}\;\!
\|\,u(\cdot,t)\,\|_{\mbox{}_{\scriptstyle L^q}}
                  ^{{\scriptstyle \;\!q}}
\!\!\!
&=&
\!\!
\int_{\mbox{}_{\mbox{\scriptsize $\!\;\!\R^3$}}} \!\!\!
|\,u(x,t)\,|^{q-2}\;\!
\langle \;\!u(x,t), \;\!u_t(x,t) \;\!\rangle \, dx \\
&=&
\!\!\!\!- \!\:\!
\int_{\mbox{}_{\mbox{\scriptsize $\!\;\!\R^3$}}} \!\!\!
|\,u\,|^{q-2}\;\!\langle u, \nabla p \rangle \;\! dx
\;\!- \!\;\!
\int_{\mbox{}_{\mbox{\scriptsize $\!\;\!\R^3$}}} \!\!\!
|\,u\,|^{q-2} \langle \:\!u, u \!\;\!\cdot\!\:\! \nabla u \:\!\rangle \;\! dx
\;\!+ \!\:\!
\int_{\mbox{}_{\mbox{\scriptsize $\!\;\!\R^3$}}} \!\!\!
|\,u\,|^{q-2} \langle \:\!u, \Delta u \:\!\rangle \;\! dx \\
&=:&
T_p \:\!+\;\! T_c \:\!+\;\! T_v
\eeaa
Using integration by parts,
one obtains that $\:\!T_c \!\;\!=\:\!0$.
Also, \\
\mbox{} \vspace{-1.050cm} \\
\beaa
T_v &=&
\!\sum_{i,\;\!j\,=\,1}^{3}
\int_{\mbox{}_{\mbox{\scriptsize $\!\;\!\R^3$}}} \!\!\!\:\!
|\,u\,|^{q-2} \,u_i \;\!D_j^2 \:\!u_i \, dx \\
&=&
\!\!-\!\;\! \sum_{i,\;\!j\,=\,1}^{3}
\int_{\mbox{}_{\mbox{\scriptsize $\!\;\!\R^3$}}} \!\!\!\:\!
|\,u\,|^{q-2} \, |\;\!\cD u\,|^2  \, dx
\;-\; (q-2) \sum_{j\,=\,1}^{3}
\int_{\mbox{}_{\mbox{\scriptsize $\!\;\!\R^3$}}} \!\!\!\:\!
|\,u\,|^{q-4} \;\!\langle \:\!u, D_j u \:\!\rangle^2 \, dx \\
&\leq &
0\:\!.
\eeaa
\mbox{} \vspace{-0.500cm} \\
The pressure term is \\
\mbox{} \vspace{-1.230cm} \\
\beaa
T_p &=&
\!- \sum_{j\,=\,1}^{3}
\int_{\mbox{}_{\mbox{\scriptsize $\!\;\!\R^3$}}} \!\!\!\:\!
|\,u\,|^{q-2} \, u_j \;\! D_j \:\!p \: dx \\
&=&
(q - 2) \sum_{j\,=\,1}^{3}
\int_{\mbox{}_{\mbox{\scriptsize $\!\;\!\R^3$}}} \!\!\!\:\!
|\,u\,|^{q-4} \;\!\langle \:\!u, D_j \:\!u \:\!\rangle \,
u_j \,p \: dx \\
&\leq &
C_{\!\;\!q} \:\|\,u(\cdot,t)\,\ni ^{q-2} \,
\|\,\cD u(\cdot,t) \,\| \: \|\,p(\cdot,t)\,\|\:\!.
\eeaa
\mbox{} \vspace{-0.320cm} \\
For $\;\!\|\,p(\cdot,t)\,\|$,
we use the bound
(from Fourier transform,
plus Parseval's relation) \\
\mbox{} \vspace{-1.070cm} \\
\beaa
\|\,p(\cdot,t)\,\|
&\leq & \!\!\!\!
\sum_{\mbox{}\;\;i,\;\!j\,=\,1}^{3} \!\!
\|\,u_i \;\!u_j \,(\cdot,t)\,\| \\
&\leq &
C \; \|\,u(\cdot,t)\,\|_{\mbox{}_{\scriptstyle L^4}}^2 \\
&\leq &
C \;\|\,u(\cdot,t)\,\|_{\mbox{}_{\scriptstyle L^3}} \,
\|\,\cD u(\cdot,t)\,\|\,\!.
\eeaa
\mbox{} \vspace{-0.500cm} \\
Thus,
we have shown the estimate \\
\mbox{} \vspace{-0.650cm} \\
\be
\label{eq.5.qest}
\f{d}{dt}\, \|\,u(\cdot,t)\,\|_{\mbox{}_{\scriptstyle L^q}}^q
\leq\,
C_{\!\;\!q} \,\|\,u(\cdot,t) \,\ni^{q-2} \,
\|\,u(\cdot,t)\,\|_{\mbox{}_{\scriptstyle L^3}} \,
\|\,\cD u(\cdot,t)\,\|^2.
\ee
\mbox{} \vspace{-0.050cm} \\
Setting\\
\mbox{} \vspace{-1.250cm} \\
\be
\label{eq.5.psi}
h(t) \,:=\;
\f{\;\!\|\,u(\cdot,t)\,\ni^{\mbox{}^{\scriptstyle \:\!q}}}
  {\:\|\,u(\cdot,t)\,\|_{\mbox{}_{\scriptstyle L^q}}
                       ^{\mbox{}^{\scriptstyle \:\!q}}}
\cdot
\f{\;\!\|\,u(\cdot,t)\,\|_{\mbox{}_{\scriptstyle L^3}}}
  {\;\!\|\,u(\cdot,t)\,\ni^{\mbox{}^{\scriptstyle 2}}}\;\!,
\ee
\mbox{} \vspace{-0.020cm} \\
we have that
if $h(t)$ were bounded by some quantity
$h_{\sf max}$
in the interval $ 0 \leq t < T_{\!\;\!f} \!\;\!$,
then
the estimate (\ref{eq.5.qest})
would give \\
\mbox{} \vspace{-0.620cm} \\
$$\f{d}{dt}\,
\|\,u(\cdot,t)\,\|_{\mbox{}_{\scriptstyle L^q}}
                  ^{\mbox{}^{\scriptstyle \:\!q}}
\;\!\leq\:
C_{\!\;\!q} \, h_{\sf max} \,
\|\,\cD u(\cdot,t) \,\|^2 \:
\|\,u(\cdot,t) \,\|_{\mbox{}_{\scriptstyle L^q}}
                   ^{\mbox{}^{\scriptstyle \:\!q}}
\!\;\!.$$
\mbox{} \vspace{-0.190cm} \\
Since
$\int_0^{T_{\!\;\!f}} \!\;\! \|\,\cD u(\cdot,t)\,\|^2\, ds\;\!$
is finite by Theorem \ref{t.2.energy},
this would give
(by Gronwall's lemma)
boundedness of
$\;\!\|\,u(\cdot,t)\,\|_{\mbox{}_{\scriptstyle L^q}}$
in $ 0 \leq t < T_{\!\;\!f}$.
This contradiction proves (\ref{eq.5.sup}).
\mbox{} \hfill $\Box$
%

\end{section}

%
%
\begin{section}{The Beale--Kato--Majda blow--up condition}
\mbox{} \vspace{-1.030cm} \\
\setcounter{equation}{0}

In this section,
we recall a few basic facts
on the flow vorticity
$ \;\!\omega(\cdot,t) := \nabla \!\times\!\;\!u(\cdot,t) $,
which satisfies the related equation \\
\mbox{} \vspace{-0.650cm} \\
\be
\label{eq.1.s6}
\label{eq.vorticity}
\omega_t \,+\:
u(\cdot,t) \!\;\!\cdot\!\;\! \nabla \, \omega(\cdot,t)
\;=\;
\Delta \;\!\omega(\cdot,t)
\,+\,
\omega(\cdot,t) \!\;\!\cdot\!\;\! \nabla \;\! u(\cdot,t).
\ee
\mbox{} \vspace{-0.220cm} \\
From our definition
of the $L^{2}$-norm
$ \;\!\|\:\cdot\:\| $,
it is readily seen that
$ {\displaystyle
\;\!
\|\, \omega(\cdot,t) \,\| =
\|\, \cD u(\cdot,t) \,\|
} $,
and, more generally, \\
\mbox{} \vspace{-0.770cm} \\
\be
\label{eq.2.s6}
\|\, \cD^{\ell + 1} \:\! u(\cdot,t) \,\|
\;=\;
\|\, \cD^{\ell} \:\! \omega(\cdot,t) \,\|
\qquad
\forall \;\, \ell \geq 0,
\ee
\mbox{} \vspace{-0.250cm} \\
so that we have,
in case $ \:\!T_{\!\;\!f} \!\;\!< \infty $,
that
$ {\displaystyle
\;\!\|\, \cD^{\ell} \:\! \omega(\cdot,t) \,\|
\;\!\rightarrow\;\!
\infty
\;\!
} $
as
$ \;\!t \:\!\mbox{\footnotesize $\nearrow$}\,\:\! T_{\!\;\!f} $
for all $ \:\!\ell \geq 0 $.
Similar considerations
are obtained from \\
\mbox{} \vspace{-0.670cm} \\
\be
\label{eq.3.s6}
\|\, \cD^{\ell + 1} \:\! u(\cdot,t) \,\|_{\mbox{}_{\scriptstyle L^{q}}}
\;\!\leq\:
K\!\;\!(\ell,q) \:
\|\, \cD^{\ell} \:\! \omega(\cdot,t) \,\|_{\mbox{}_{\scriptstyle L^{q}}}
\!\!
\qquad
\forall \;\, \ell \geq 0,
\;\:
1 < q < \infty,
\ee
\mbox{} \vspace{-0.250cm} \\
which follows from the
Calderon--Zygmund theory
of singular operators,
see e.g.$\;$\cite{Galdi1994, LopesNussenzveigZheng1999}.
Another important property
of $ \omega(\cdot,t) $
is that it stays bounded in $ L^{1} \!\:\!$,
as observed in \cite{Constantin1990}. \\
\mbox{} \vspace{-0.950cm} \\
%

%
%
%
%
\begin{theorem}
$(i)$
Let
$ {\displaystyle
\,
\omega =
\bigl(\;\!
\omega_{\mbox{}_{1}}\!\;\!, \;\!\omega_{\mbox{}_{2}}\!\;\!, \;\!\omega_{\mbox{}_{3}}
\bigr)
\:\!
} $
be the vorticity.
If
$\,\omega_{\mbox{}_{\scriptstyle i}}\!\;\!(\cdot,0) \in L^{1}(\R^3) \:\!$
for some~$\;\!i$,
then
$\;\!\omega_{\mbox{}_{\scriptstyle i}}\!\;\!(\cdot,t) \:\!$
remains in $\:\! L^{1}(\R^3) \:\!$
for $\,t > 0 $,
with \\
\mbox{} \vspace{-0.700cm} \\
\be
\label{eq.4.s6}
\|\, \omega_{\mbox{}_{\scriptstyle \!\;\!i}}\!\;\!(\cdot,t) \,\|_{\mbox{}_{\scriptstyle L^{1}}}
\;\!\leq\;\;\!
\|\, \omega_{\mbox{}_{\scriptstyle \!\;\!i}}\!\;\!(\cdot,0) \,\|_{\mbox{}_{\scriptstyle L^{1}}}
+\,
\mbox{\small $ {\displaystyle \frac{\,1\,}{2} }$}
\,
\|\, u(\cdot,0) \,\|^{2}
\qquad
\forall \;\,
0 \leq t < T_{\!\;\!f}.
\ee
\mbox{} \vspace{-0.200cm} \\
$(ii)$
If
$\,\omega(\cdot,0) \in L^{1}(\R^3) $,
then \\
\mbox{} \vspace{-0.720cm} \\
\be
\label{eq.5.s6}
\|\, \omega(\cdot,t) \,\|_{\mbox{}_{\scriptstyle L^{1}}}
\;\!\leq\;\;\!
\|\, \omega(\cdot,0) \,\|_{\mbox{}_{\scriptstyle L^{1}}}
+\,
\mbox{\small $ {\displaystyle \frac{\sqrt{\:\!3\,}\,}{2} }$}
\,
\|\, u(\cdot,0) \,\|^{2}
\qquad
\forall \;\,
0 \leq t < T_{\!\;\!f}.
\ee
\end{theorem}
%
%

%
\mbox{} \vspace{-0.200cm} \\
{\bf Proof:}
Again,
we use regularized sign functions
$ \,\!L^{\prime}_{\delta}(\:\!\cdot\:\!) \:\!$
as defined in \cite{KreissLorenz1989}, p.$\:$136,
where $ \delta > 0 $ is arbitrary.
Multiplying the $i$-th component of
equation (6.1) above
by
$\:\!L^{\prime}_{\delta}(\omega_{\mbox{}_{\scriptstyle \!\;\!i}}(\cdot,t)) $
and
integrating on $\;\!\R^{3} \!\times\!\;\! [\,0, \;\!t\;\!] $,
we get,
letting $ \delta \,\mbox{\footnotesize $\searrow$}\, 0 $, \\
\mbox{} \vspace{-0.800cm} \\
$$
\|\, \omega_{\mbox{}_{\scriptstyle \!\;\!i}}\!\;\!(\cdot,t) \,\|_{\mbox{}_{\scriptstyle L^{1}}}
\;\!\leq\;\;\!
\|\, \omega_{\mbox{}_{\scriptstyle \!\;\!i}}\!\;\!(\cdot,0) \,\|_{\mbox{}_{\scriptstyle L^{1}}}
+\,
\sum_{j\,=\,1}^{3}
\int_{\mbox{}_{\scriptstyle 0}}^{\mbox{\footnotesize $\:\!t$}} \!\!\;\!
\int_{\mbox{}_{\mbox{\scriptsize $\!\;\!\R^3$}}} \!\!\!
|\, \omega_{\mbox{}_{\scriptstyle \!\;\!j}}(x,t) \,| \;
|\, D_{\mbox{}_{\scriptstyle \!j}}
u_{\mbox{}_{\scriptstyle \!\;\!i}}(x,t) \,| \;
dx \, d\tau
$$
\mbox{} \vspace{-0.200cm} \\
for all $\;\!0 \leq t < T_{\!\;\!f} $,
which gives
(\ref{eq.4.s6})
by the Cauchy--Schwarz inequality and
(\ref{integral.Du.L2}), (\ref{eq.2.s6}).
Moreover,
summing on $\;\!1 \leq i \leq 3 \;\!$
and applying again the Cauchy--Schwarz inequality
to estimate the integral term,
one obtains (\ref{eq.5.s6}),
using (\ref{integral.Du.L2}) and (\ref{eq.2.s6})
once more.
\mbox{} \hfill $\Box$
%
%

\newpage

\mbox{} \vspace{-0.750cm} \\

We now turn to the
$L^{2}$ norm of $ \omega(\cdot,t) $,
which will quickly lead us to the
following blow--up result,
originally obtained for the
Euler's equations in \cite{BealeKatoMajda1984}.

%
%
%
%
\begin{theorem}
$($Beale--Kato--Majda\/$)$.
If $\,T_{\!\;\!f} \!\;\!< \infty $,
then
$ {\displaystyle
\:\!
\int_{\mbox{}_{\scriptstyle 0}}^{\;\!T_{\mbox{}_{\!\,\!f}}}
\!\:\!
\|\, \omega(\cdot,t) \,\|_{\mbox{}_{\scriptstyle \:\!\infty}} \,dt
\:=\: \infty
} $.
\end{theorem}
%
%

%
\mbox{} \vspace{-0.400cm} \\
{\bf Proof:}
Multiplying the $i$-th component of
equation (6.1)
by
$\;\!\omega_{\mbox{}_{\scriptstyle \!\;\!i}}(\cdot,t) $
and
integrating on $\;\!\R^{3} \!\times\!\;\! [\,0, \;\!t\;\!] $,
we get,
summing on $ \:\!1 \leq i \leq 3 $, \\
\mbox{} \vspace{-0.050cm} \\
\mbox{} \hspace{+0.750cm}
$ {\displaystyle
\mbox{\normalsize $ {\displaystyle \frac{\,1\,}{2} }$} \:
\mbox{\normalsize $ {\displaystyle \frac{d}{\;\!d \:\!t\;\!} }$} \,
\|\, \omega(\cdot,t) \,\|^{2}
\;\!+\;
\|\, \cD \omega(\cdot,t) \,\|^{2}
\;=
\sum_{\mbox{}\;i, \,j \,=\, 1}^{3}
\int_{\mbox{}_{\mbox{\scriptsize $\!\;\!\R^3$}}} \!\!\!\:\!
\omega_{\mbox{}_{\scriptstyle i}}(x,t) \:
\omega_{\mbox{}_{\scriptstyle \!\:\!j}}(x,t) \,
D_{\scriptstyle \!\:\!j} \:\!u_{\mbox{}_{\scriptstyle i}}(x,t)
\, dx
} $ \\
\mbox{} \vspace{+0.050cm} \\
\mbox{} \hspace{+6.700cm}
$ {\displaystyle
\leq\;\,
\|\, \omega(\cdot,t) \,\|_{\mbox{}_{\mbox{}_{\scriptstyle \infty}}}
\!\!
\sum_{\mbox{}\;i, \,j \,=\, 1}^{3}
\int_{\mbox{}_{\mbox{\scriptsize $\!\;\!\R^3$}}} \!\!\!\:\!
|\, \omega_{\mbox{}_{\scriptstyle \!\:\!j}}(x,t) \,| \;
|\;\!D_{\scriptstyle \!\:\!j} \:\! u_{\mbox{}_{\scriptstyle i}}(x,t) \,|
\, dx
} $ \\
\mbox{} \vspace{+0.150cm} \\
\mbox{} \hspace{+6.700cm}
$ {\displaystyle
\leq\;
\sqrt{\:\!3\:} \;
\|\, \omega(\cdot,t) \,\|_{\mbox{}_{\mbox{}_{\scriptstyle \infty}}}
\;\!
\|\, \omega(\cdot,t) \,\|\;\;\!
\|\, \cD u(\cdot,t) \,\|
} $ \\
\mbox{} \vspace{+0.000cm} \\
\mbox{} \hspace{+6.700cm}
$ {\displaystyle
=\;
\sqrt{\:\!3\:} \;
\|\, \omega(\cdot,t) \,\|_{\mbox{}_{\mbox{}_{\scriptstyle \infty}}}
\;\!
\|\, \omega(\cdot,t) \,\|^{2}
} $ \\
\mbox{} \vspace{+0.050cm} \\
for all
$\;\! 0 \leq t < T_{\!\;\!f} $,
using the Cauchy-Schwarz inequality
and
that
$ {\displaystyle
\;\!
\|\, \cD u(\cdot,t) \,\| \;\!=\;\!
\|\, \omega(\cdot,t) \,\|
} $,
see (\ref{eq.2.s6}) above.
By the standard Gronwall lemma,
this gives \\
\mbox{} \vspace{-0.525cm} \\
$$
\|\, \omega(\cdot,t) \,\|
\;\;\!\leq\;\,
\|\, \omega(\cdot,0) \,\|
\:\!\cdot\,
\mbox{exp}\,
\biggl\{\:\! \sqrt{\:\!3\:} \!\:\!
\int_{\mbox{}_{\scriptstyle 0}}^{\mbox{\footnotesize $\;\!t$}} \!\!\;\!
\|\, \omega(\cdot,\tau) \,\|_{\mbox{}_{\mbox{}_{\scriptstyle \infty}}}
\:\!
d\tau
\,\biggr\}
\qquad
\forall \;\,
0 \leq t < T_{\!\;\!f}.
$$
\mbox{} \vspace{-0.075cm} \\
Therefore,
if we had
$ {\displaystyle
\:\!
\int_{\mbox{}_{\scriptstyle 0}}^{\;\!T_{\mbox{}_{\!\,\!f}}} \!
\|\, \omega(\cdot,t) \,\|_{\mbox{}_{\scriptstyle \:\!\infty}} \,dt
\;\!< \infty
} $,
we would have
$ {\displaystyle
\;\!
\|\, \omega(\cdot,t) \,\|
\;\!
} $
bounded in
$ \;\![\,0, \;\!T_{\!\;\!f}\;\![ $. \linebreak
\mbox{} \vspace{-0.450cm} \\
That is,
$ \|\, \cD u(\cdot,t) \,\| $
would be
bounded in
$ \;\![\,0, \;\!T_{\!\;\!f}\;\![ $,
contradicting
Leray's estimate
(\ref{eq.10.s3}).
\mbox{} \hfill $\Box$ \\
%
%

\end{section}

\mbox{} \vspace{-1.400cm} \\
%

%
%

\begin{section}{Appendix: some auxiliary results}

\setcounter{equation}{0}
%

%
%
\mbox{} \vspace{-1.500cm} \\
\begin{subsection}{Heat equation estimates}
\label{s.7.1}
\mbox{} \vspace{-0.320cm} \\
Consider the Cauchy problem for the heat equation, \\
\mbox{} \vspace{-0.650cm} \\
\be
u_t \;\!=\:\Delta u,
\quad u\;\!= f \qat t = 0, \;\; x \in \R^N\!.
\ee
\mbox{} \vspace{-0.250cm} \\
We assume that
$f \in L^r(\R^N\!\;\!)$
for some  $1 \leq r \leq \infty$.
The solution is \\
\mbox{} \vspace{-0.250cm} \\
\mbox{} \hspace{+4.800cm}
$ {\displaystyle
\label{eq.ap.conv}
u(x,t) \,=\;\!
\int_{\mbox{}_{\mbox{\scriptsize $\!\;\!\R^N$}}} \!\!\!\!\!\:\!
\Phi(x-y,t) \:\!f(y)\, dy
} $,
\mbox{} \hfill $(7.2a)$ \\
\mbox{} \vspace{-0.150cm} \\
where $\Phi$ is the heat kernel, \\
\mbox{} \vspace{-0.200cm} \\
\mbox{} \hspace{+5.000cm}
$ {\displaystyle
\Phi(x,t) \,=\:
(\:\!4\pi t\:\!)^{-N/2} \;\!e^{-\,|\,x\,|^2/4t}\!\;\!.
} $
\mbox{} \hfill $(7.2b)$ \\
\setcounter{equation}{2}
\mbox{} \vspace{-0.250cm} \\
We have \\
\mbox{} \vspace{-1.000cm} \\
\be
\label{eq.ap.heat}
0 < \Phi(x,t) \leq (\:\!4\pi t\:\!)^{-N/2}\!\;\!,
\quad \;\;
\int_{\mbox{}_{\mbox{\scriptsize $\!\;\!\R^N$}}} \!\!\!\!\!\:\!
\Phi(x,t)\,dx \,=\;\! 1 \qfor t > 0\:\!.
\ee
\mbox{} \vspace{-0.200cm} \\
It is elementary to show that
the function $u(x,t)$
defined by (\ref{eq.ap.conv})
is a $C^\infty\!\;\!$ function for $t > 0$.
Furthermore, for $t > 0$,
all derivatives of $u(x,t)$
can be obtained by differentiating
the convolution integral (7.2$a$)
under the integral sign.
In particular,
we have, for all space derivatives
$D^\alpha u$, \\
\mbox{} \vspace{-0.800cm} \\
\be
D^\alpha u(x,t) \,=
\int_{\mbox{}_{\mbox{\scriptsize $\!\;\!\R^N$}}} \!\!\!\!\!
D^{\alpha}_{\!x} \:\!\Phi(x-y,t) \;\! f(y)\: dy.
\ee
\mbox{} \vspace{-0.150cm} \\
We show the following solution estimates:

%
%
%
%
\begin{theorem}
Let $1 \leq r \leq q \leq \infty$.
Then we have \\
\mbox{} \vspace{-0.140cm} \\
\mbox{} \hspace{+4.000cm}
$ {\displaystyle
\label{ap.heat.1}
\|\,u(\cdot,t)\,\|_{\mbox{}_{\scriptstyle L^q}} \leq\:
(\:\!4\pi t\:\!)^{-\lambda} \, \|\,f\,\|_{\mbox{}_{\scriptstyle L^r}}
\quad \;\,
\forall \;\, t > 0
} $,
\mbox{} \hfill $(7.5a)$ \\
\mbox{} \vspace{-0.130cm} \\
and,
for every space derivative $D^\alpha\!\,\!$, \\
\mbox{} \vspace{-0.150cm} \\
\mbox{} \hspace{+3.500cm}
$ {\displaystyle
\label{ap.heat.2}
\|\;\!D^\alpha u(\cdot,t)\,\|_{\mbox{}_{\scriptstyle L^q}}
\leq\:
C \,t^{-\lambda - |\,\alpha\,|/2} \,
\|\,f\,\|_{\mbox{}_{\scriptstyle L^r}}
\quad \;\,
\forall \;\, t > 0,
} $
\mbox{} \hfill $(7.5b)$ \\
\mbox{} \vspace{-0.250cm} \\
with \\
\mbox{} \vspace{-0.550cm} \\
\mbox{} \hspace{+6.100cm}
$ {\displaystyle
\lambda \,=\:
\f{N}{2} \;\!\Bigl(\:\!\f{1}{r} - \f{1}{q} \:\!\Bigr).
} $
\mbox{} \hfill $(7.5c)$ \\
\mbox{} \vspace{+0.050cm} \\
The constant $C$ depends on $\;\!r, q, \alpha$ and $\;\!N\!\:\!$,
but is independent of $\;\!t$ and $f\!\;\!$.
\end{theorem}
%
%
\setcounter{equation}{5}

To prove (7.5$a$),
we first show that the $L^r$--norm of
$\;\!u(\cdot,t)$ cannot grow in time:

%
%
%
%
\begin{lemma}
\label{l.ap.3}
For $\,1 \leq r \leq \infty$,
we have \\
\mbox{} \vspace{-0.770cm} \\
\be
\|\,u(\cdot,t)\,\|_{\mbox{}_{\scriptstyle L^r}}
\leq\; \|\,f\,\|_{\mbox{}_{\scriptstyle L^r}}
\quad \;\,
\forall \;\,t \geq 0\:\!.
\ee
\end{lemma}
%
%
\mbox{} \vspace{-0.350cm} \\
{\bf Proof:}
For $r = \infty$ and for $r = 1$
the estimate follows directly from
$ {\displaystyle
\int_{\mbox{}_{\mbox{\scriptsize $\!\;\!\R^N$}}} \!\!\!\!\!\:\!
|\, \Phi(x,t)\,| \, dx \,=\, 1
} $. \\
\mbox{} \vspace{-0.575cm} \\
For $\;\!1 < r < \infty$,
define $r'$ by \\
\mbox{} \vspace{-0.950cm} \\
$$\f{1}{r'} \,+\, \f{1}{r} \,=\,1. $$
\mbox{} \vspace{-0.400cm} \\
Using H\"older's inequality, \\
\mbox{} \vspace{-1.050cm} \\
\beaa
|\,u(x,t)\,| &\leq &
\int_{\mbox{}_{\mbox{\scriptsize $\!\;\!\R^N$}}} \!\!\!\!\!\:\!
\Phi(x-y,t)^{1/r'} \;\!
\Phi(x-y,t)^{1/r} \;\! |\,f(y)\,|\: dy \\
&\leq &
\!\!\Bigl(
\int_{\mbox{}_{\mbox{\scriptsize $\!\;\!\R^N$}}} \!\!\!\!\!\:\!
\Phi(x-y,t)\,dy \;\!\Bigr )^{\!1/r'} \!\;\!
\Bigl(
\int_{\mbox{}_{\mbox{\scriptsize $\!\;\!\R^N$}}} \!\!\!\!\!\:\!
\Phi(x-y,t) \,|\,f(y)\,|^r\,dy \;\!\Bigr)^{\!1/r}\!,
\eeaa
\mbox{} \vspace{-0.450cm} \\
so that
we have \\
\mbox{} \vspace{-1.000cm} \\
$$ |\,u(x,t)\,|^r  \;\!\leq\:\!
\int_{\mbox{}_{\mbox{\scriptsize $\!\;\!\R^N$}}} \!\!\!\!\!\:\!
\Phi(x-y,t) \, |\,f(y)\,|^r \, dy. $$
\mbox{} \vspace{-0.200cm} \\
Integration in $x$ proves the lemma.
\mbox{} \hfill $\Box$ \\
%
%
\mbox{} \vspace{-0.650cm} \\

In the next lemma,
we estimate the maximum norm of $u(\cdot,t)$ by
$\|\,f\,\|_{\mbox{}_{\scriptstyle L^r}}\!\;\!$.

%
%
%
%
\begin{lemma}
\label{l.ap.4}
For $\;\!1 \leq r < \infty$,
we have \\
\mbox{} \vspace{-0.750cm} \\
\be
\|\,u(\cdot,t)\,\ni \:\!\leq\;\!
(\:\!4\pi t\:\!)^{-N/(2\:\!r)}\;\!
\|\,f\,\|_{\mbox{}_{\scriptstyle L^r}}.
\ee
\end{lemma}
%
%
{\bf Proof:}
For $r = 1$,
the estimate follows
from the bound
$\Phi(x,t) \!\;\!\leq\!\;\! (4\pi t)^{-N/2}\!\;\!$.
If \mbox{$\;\!1 < r < \infty$}, \linebreak
define $r'$ by $\f{1}{r'} + \f{1}{r} = 1\;\!$
and obtain \\
\mbox{} \vspace{-1.100cm} \\
\beaa
|\,u(x,t)\,| &\leq &
\!\!
\int_{\mbox{}_{\mbox{\scriptsize $\!\;\!\R^N$}}} \!\!\!\!\!\:\!
\Phi(x-y,t)^{1/r'} \Phi(x-y,t)^{1/r} \;\!|\,f(y)\,|\: dy \\
&\leq &
\!
(\:\!4\pi t\:\!)^{-N/(2\:\!)r} \!
\int_{\mbox{}_{\mbox{\scriptsize $\!\;\!\R^N$}}} \!\!\!\!\!\:\!
\Phi(x-y,t)^{1/r'} \:\!|\,f(y)\,|\:dy
\;\,\leq\;\,
(\:\!4\pi t\:\!)^{-N/2r} \;\! \|\,f\,\|_{\mbox{}_{\scriptstyle L^r}}
\!\;\!.
\eeaa
In the last estimate we have used H\"older's inequality,
and the lemma is proved.
\mbox{} \hfill $\Box$ \\
%
%
\mbox{} \vspace{-0.750cm} \\

Using the bounds of the two previous lemmas,
the estimate (7.5$a$)
follows by a simple  argument:
for $\;\!1 \leq r \leq q < \infty$,
we have \\
\mbox{} \vspace{-1.100cm} \\
\beaa
\|\,u(\cdot,t)\,\|_{\mbox{}_{\scriptstyle L^q}}^q
&=&
\int_{\mbox{}_{\mbox{\scriptsize $\!\;\!\R^N$}}} \!\!\!\!\!\:\!
|\,u(x,t)\,|^q \: dx \\
&=&
\int_{\mbox{}_{\mbox{\scriptsize $\!\;\!\R^N$}}} \!\!\!\!\!\:\!
|\,u(x,t)\,|^{q-r} \;\!|\,u(x,t)\,|^r  \;\! dx \\
&\leq &
\|\,u(\cdot,t)\,\ni^{q-r} \;\!\|\,u(\cdot,t)\,\|_{\mbox{}_{\scriptstyle L^r}}^r \\
&\leq &
(\:\!4\pi t\:\!)^{-\;\!(q-r)N/(2\:\!r)} \;\!
\|\,f\,\|_{\mbox{}_{\scriptstyle L^r}}^q
\!\;\!.
\eeaa
\mbox{} \vspace{-0.475cm} \\
Taking the $q$--th root, we obtain (7.5$a$).
\mbox{} \hfill $\Box$ \\
%
\mbox{} \vspace{-0.200cm} \\
%
%
{\bf Estimates of Derivatives:}
Consider the function \\
\mbox{} \vspace{-0.650cm} \\
$$g(z) \,=\: e^{-\,|\,z\,|^2\!\;\!/4}\!\;\!,
\quad |\,z\,|^2 =\, z_1^2 + \ldots + z_N^2. $$
\mbox{} \vspace{-0.300cm} \\
It is easy to see that \\
\mbox{} \vspace{-1.050cm} \\
$$D^\alpha_{\!z} \:\!g(z) =\;\! p_\alpha(z) \;\!g(z), $$
\mbox{} \vspace{-0.300cm} \\
where $p_\alpha(z)$ is a polynomial.
Therefore,
$D^\alpha _{\!z} \:\!g(z)\:\!$
is a bounded function,
and \\
\mbox{} \vspace{-0.600cm} \\
$$ \int_{\mbox{}_{\mbox{\scriptsize $\!\;\!\R^N$}}} \!\!\!\!\!\:\!
 |\;\!D^\alpha_{\!z} \;\!g(z)\,| \: dz\, =: \;\!
c(\alpha,N) < \infty. $$
Since \\
\mbox{} \vspace{-1.000cm} \\
$$\Phi(x,t) \,=\, (\:\!4\pi t\:\!)^{-\:\!N/2}\;\!
g(x\;\!t^{-1/2}), $$
we have \\
\mbox{} \vspace{-1.000cm} \\
$$D^\alpha_{\!x} \;\!\Phi(x,t) \,=\,
(\:\!4\pi t\:\!)^{-N/2}\, t^{-|\,\alpha\,|/2} \;\!
(D^\alpha _{\!z} \:\!g)(x \;\!t^{-1/2}). $$
\mbox{} \vspace{-0.250cm} \\
Integrating over $x\in \R^N\!\;\!$
and
using the substitution \\
\mbox{} \vspace{-0.650cm} \\
$$ x \;\!t^{-1/2} \,=\, z,
\quad  dx \,=\: t^{N/2} \;\!dz, $$
\mbox{} \vspace{-0.650cm} \\
one obtains that \\
\mbox{} \vspace{-1.100cm} \\
\beaa
\int_{\mbox{}_{\mbox{\scriptsize $\!\;\!\R^N$}}} \!\!\!\!\!\:\!
|\,D^\alpha_{\!x} \;\!\Phi(x,t)\,|\: dx
&=&
(\:\!4\pi t\:\!)^{-\;\!N/2} \,t^{(N - |\,\alpha\,|\;\!)/2}
\;\!c(\alpha,N) \\
&=&
\tilde{c}(\alpha ,N) \, t^{-|\,\alpha\,|/2}\!\;\!.
\eeaa
One can now repeat the arguments given above
for the case $\alpha = 0$
and derive estimates of
$\;\!\|\;\!D^\alpha u(\cdot,t)\,\|_{\mbox{}_{\scriptstyle L^q}}$
in terms of
$\;\!\|\,f\,\|_{\mbox{}_{\scriptstyle L^r}}\!\;\!$.
Instead of \\
\mbox{} \vspace{-0.600cm} \\
$$ \int_{\mbox{}_{\mbox{\scriptsize $\!\;\!\R^N$}}} \!\!\!\!\!\:\!
\Phi(x,t)\, dx \,=\, 1, $$
\mbox{} \vspace{-0.400cm} \\
one uses the equation \\
\mbox{} \vspace{-0.750cm} \\
$$\int_{\mbox{}_{\mbox{\scriptsize $\!\;\!\R^N$}}} \!\!\!\!\!\:\!
|\,D^\alpha_{\!x} \;\!\Phi(x,t)\,|\, dx \,=\,
c\,t^{-\,|\,\alpha\,|/2}\!\;\!, $$
\mbox{} \vspace{-0.050cm} \\
with
$ c > 0 $
constant
depending on $ \alpha $, \mbox{\small $N$}
only.
In this way,
the estimate (7.5$b$) follows.
\mbox{} \hfill $\Box$ \\
\mbox{} \\
\mbox{} \vspace{-0.150cm} \\
%
%
%
%
{\sc Remark:}
Consider a convolution integral \\
\mbox{} \vspace{-0.600cm} \\
$$ v(x) \;=\,
\int_{\mbox{}_{\mbox{\scriptsize $\!\;\!\R^N$}}} \!\!\!\!\!
k(x-y) \, f(y) \: dy,  $$
\mbox{} \vspace{-0.170cm} \\
where $ k \in L^1(\R^N\!\:\!), f \in L^r(\R^N\!\:\!)$.
Arguing as in the proof of Lemma \ref{l.ap.3},
one obtains that \\
\mbox{} \vspace{-0.700cm} \\
$$\|\;v\;\|_{\mbox{}_{\scriptstyle L^r}}
\:\!\leq\;
\|\;k\;\|_{\mbox{}_{\scriptstyle L^1}} \:\!
\|\,f\,\|_{\mbox{}_{\scriptstyle L^r}}\!\;\!, $$
\mbox{} \vspace{-0.330cm} \\
which is Young's inequality.
As before,
define $r'$ by
$\,\f{1}{r} + \f{1}{r'} = 1$.
If one argues as in the proof of Lemma \ref{l.ap.4},
one obtains the bound \\
\mbox{} \vspace{-0.700cm} \\
$$\|\;v\;\ni
\;\!\leq\;
\|\,k\,\|_{\mbox{}_{\scriptstyle L^{r'}}} \;\!
\|\,f\,\|_{\mbox{}_{\scriptstyle L^r}}\!\;\!. $$
\mbox{} \vspace{-0.250cm} \\
Therefore,
the estimates (7.5$a$) and (7.5$b$)
can also be proved
by computing the $L^q$--norms of
$\;\!\Phi(\cdot,t)$ and $D^\alpha \;\!\Phi(\cdot,t)$.

\end{subsection}

%
%
\mbox{} \vspace{-0.900cm} \\
\begin{subsection}{The Helmholtz projector $P_{\scriptscriptstyle \!\;\!H}$}
\label{s.7.2}
\mbox{} \vspace{-0.870cm} \\

Let $v \!\;\!: \R^3 \to \R^3$ denote a given vector field.
If $v$ satisfies some restrictions,
there is a unique decomposition of $\;\!v\:\!$
into a divergence--free field, $w$,
and a gradient field, $\nabla \phi$, \\
\mbox{} \vspace{-0.600cm} \\
\be
\label{eq.ap.helm}
v \,=\,
w \,-\:\! \nabla \phi,
\qquad
\nabla \!\;\!\cdot\;\! w \,=\,0\:\!.
\ee
\mbox{} \vspace{-0.300cm} \\
Then,
in suitable function spaces,
the assignment
$v \to w =: \!\:\!P_{\scriptscriptstyle \!\;\!H}v \;\!$
defines a projection oper\-ator $P_{\scriptscriptstyle \!\;\!H}$,
called the Helmholtz projector.
An important and nontrivial result
is the boundedness of $P_{\scriptscriptstyle \!\;\!H}$
in $L^q$ for every $q$ with $1 < q < \infty$.
We will not prove this result here,
but only discuss how it is related to the
Calderon--Zygmund theory of singular integrals. \\
\mbox{} \vspace{-0.950cm} \\

To begin with, let $ v \in C^\infty_0 $, i.e.,
$v$ is  smooth and compactly supported.
First assume that a decomposition (\ref{eq.ap.helm})
holds with $C^\infty$ functions $w$ and $\phi$.
Then,
taking the divergence of the equation
$v \;\!=\;\! w\;\! - \nabla \phi$,
one obtains \\
\mbox{} \vspace{-0.050cm} \\
\mbox{} \hspace{+6.150cm}
$ {\displaystyle
-\,\Delta \:\!\phi \:=\: \nabla \!\:\!\cdot\:\! v\:\!.
} $
\mbox{} \hfill $(7.9a)$ \\
\mbox{} \vspace{-0.070cm} \\
Assuming that
$\:\!\phi(x) \to 0\;\!$ as $\;\!|\,x\,| \to \infty$,
we have \\
\mbox{} \vspace{-0.050cm} \\
\mbox{} \hspace{+4.250cm}
$ {\displaystyle
\label{eq.ap.helm2}
\phi(x) \:=\; \f{1}{4\:\!\pi} \!\;\!
\int_{\mbox{}_{\mbox{\scriptsize $\!\;\!\R^3$}}} \!\!\!\:\!
|\,x - y\,|^{-1} \;\!\nabla \!\cdot v(y) \, dy\:\!.
} $
\mbox{} \hfill $(7.9b)$ \\
\setcounter{equation}{9}
\mbox{} \vspace{-0.120cm} \\
One can now reverse the process,
defining $\phi$ by (7.9$b$)
and setting \\
\mbox{} \vspace{-0.600cm} \\
\be
w \:=\:
v \,+\, \nabla \phi
\;=:\,
P_{\scriptscriptstyle \!\;\!H} v\:\!.
\ee
\mbox{} \vspace{-0.300cm} \\
Then the Helmholtz decomposition (\ref{eq.ap.helm}) is obtained,
leading to the following definition. \\
\mbox{} \vspace{-0.100cm} \\
%
%
%
%
%
\noindent
{\bf Definition 7.1:}
{\it If $\;\!v \in C_0^\infty(\R^3\!\;\!)$,
then define
$\:\!w \;\!=\;\! P_{\scriptscriptstyle \!\;\!H}v\,\!$
by} \\
\mbox{} \vspace{-0.150cm} \\
\mbox{} \hspace{+3.000cm}
$ {\displaystyle
w_j(x)
\;=\;
v_j(x) \,+\, \f{1}{4\pi} \sum_{i\,=\,1}^{3}
\int_{\mbox{}_{\mbox{\scriptsize $\!\;\!\R^3$}}} \!\!\!\!\:\!
D_{\!\;\!x_j}|\,x - y\,|^{-1}
\;\! D_i \;\! v_i(y)\: dy
} $
\mbox{} \hfill $(7.11a)$ \\
\mbox{} \vspace{+0.100cm} \\
\mbox{} \hspace{+4.250cm}
$ {\displaystyle
\label{eq.ap.w}
=\;
v_j(x) \,+\, \f{1}{4\pi} \sum_{i\,=\,1}^{3}
\int_{\mbox{}_{\mbox{\scriptsize $\!\;\!\R^3$}}} \!\!\!\:\!
|\,x - y\, |^{-1} \;\!
D_j D_i \;\!v_i(y)\: dy,
} $
\mbox{} \hfill $(7.11b)$ \\
\setcounter{equation}{11}
\mbox{} \vspace{-0.050cm} \\
{\it for each $1 \leq j \leq 3$.} \\
%
%
\mbox{} \vspace{-0.100cm} \\
The following result about
$P_{\scriptscriptstyle \!\;\!H}\!$
is easy to prove: \\
\mbox{} \vspace{-1.000cm} \\
%
%
%
%
\begin{lemma}
\label{l.ap.6}
Let $\;\!v\in C_0^\infty(\R^3\!\;\!)$.\\
$(i)$ Define $\;\!w \:\!= P_{\scriptscriptstyle \!\;\!H}v$
by $\:\!(7.11b)$
and define $\phi$ by $\:\!(7.9b)$.
Then
$\;\!w \in C^\infty(\R^3\!\;\!), \,\phi \in C^\infty(\R^3\!\;\!)\:\!$
and $\:\!(\ref{eq.ap.helm})$ holds.
Furthermore,
we have \\
\mbox{} \vspace{-0.700cm} \\
\be
\label{eq.ap.decay}
|\,w(x)\,| \:+\: |\,\phi(x)\,| \,\to\, 0 \qas |\,x\,|\to \infty\:\!.
\ee
\mbox{} \vspace{-0.300cm} \\
$(ii)$ Conversely,
if $\;\!w, \,\phi \in C_0^{\infty}(\R^3\!\;\!) $
satisfy $\:\!(\ref{eq.ap.helm})$ and $\:\!(\ref{eq.ap.decay})$,
then $\:\!w\:\!$ and $\:\!\phi\:\!$
agree with the functions
defined in $\:\!(7.11b)$ and $\:\!(7.9b)$
above.
\end{lemma}
%
%
\mbox{} \vspace{-1.100cm} \\

So far we have assumed $v \in C^\infty_0\!\:\!$.
$\!\;\!$To obtain estimates of
$P_{\scriptscriptstyle \!\;\!H}v\:\!$
in terms of $\;\!v$,
which allow \linebreak
to extend the operator
$\!\;\!P_{\scriptscriptstyle \!\;\!H}\!\:\!$
to less regular functions,
the theory of singular integrals can be applied.
To discuss this,
let \\
\mbox{} \vspace{-0.700cm} \\
$$K_{ij}(z)\:\!:= \;\!D_i\;\!D_j \Bigl(\,|\,z\,|^{-1}\:\!\Bigr)
\;\qfor z \in \R^3\!,
\quad z \not= 0\:\!. $$
\mbox{} \vspace{-0.400cm} \\
Then we have \\
\mbox{} \vspace{-0.750cm} \\
$$K_{ij}(z) \,=\, 3 \;\!z_i \;\!z_j\;\!|\,z\,|^{-5} \qfor i \not= j
\qand K_{jj}(z) \,=\, (\:\!3\;\! z_j^2 - |\,z\,|^2)\:|\,z\,|^{-5}.$$
\mbox{} \vspace{-0.300cm} \\
Thus,
the kernels $K_{ij}(z)$ are homogeneous of degree $-3$,
i.e., \\
\mbox{} \vspace{-0.700cm} \\
$$K_{ij}(z) \,=\, |\,z\,|^{-3} \;\!K_{ij}(z^0),
\;\quad z^0 = z/|\,z\,|\:\!. $$
\mbox{} \vspace{-0.300cm} \\
If one integrates by parts in (7.9$b$)
and formally differentiates under the integral sign,
then one obtains \\
\mbox{} \vspace{-0.350cm} \\
\mbox{} \hspace{+3.850cm}
$ {\displaystyle
\label{eq.ap.helm4}
D_j \:\!\phi(x) \:=\: \f{1}{4\:\!\pi} \:\!
\sum_{i\,=\,1}^{3}
\int_{\mbox{}_{\mbox{\scriptsize $\!\;\!\R^3$}}} \!\!\!\!\;\!
K_{ij}(x-y) \, v_i(y)\, dy\:\!.
} $
\mbox{} \hfill $(7.13a)$ \\
\mbox{} \vspace{+0.175cm} \\
However,
since $|\;\!K_{ij}(z)\,| \sim |\,z\,|^{-3}$,
the integrals in (7.13$a$)
do not exist as Lebesgue integrals,
but must be interpreted as {\it principle values},
that is, \\
\mbox{} \vspace{-0.150cm} \\
\mbox{} \hspace{+3.500cm}
$ {\displaystyle
\label{eq.ap.helm5}
D_j \:\!\phi(x) \:=\: \f{1}{4\:\!\pi}\:\!
\sum_{i\,=\,1}^{3}
\lim_{\eps \,\to\, 0}
\int_{\mbox{}_{\scriptstyle \!\;\!|\,x - y\,| \,\geq\, \eps}}
\hspace{-1.100cm} K_{ij}(x-y) \, v_i(y)\, dy\:\!.
} $
\mbox{} \hfill $(7.13b)$ \\
\setcounter{equation}{13}
\mbox{} \vspace{+0.150cm} \\
Together with the equation $\;\!w = v + \nabla \phi$,
one is led to the following definition of
$P_{\scriptscriptstyle \!\;\!H}v$. \\
\mbox{} \vspace{-0.100cm} \\
%
%
%
%
{\bf Definition 7.2:}
If $\:\!v \in C_0^\infty(\R^3\!\;\!)$,
then define
$\;\!w \;\!=\:\! P_{\scriptscriptstyle \!\;\!H}v \:\!$
by \\
\mbox{} \vspace{-0.750cm} \\
\be
w_j(x) \:=\: v_j(x) \,+\,
\f{1}{4\:\!\pi}\:\!
\sum_{i\,=\,1}^{3}
\lim_{\eps \,\to\, 0}
\int_{\mbox{}_{\scriptstyle \!\;\!|\,x - y\,| \,\geq\, \eps}}
\hspace{-1.100cm} K_{ij}(x-y) \, v_i(y)\: dy,
\quad \;\;
1 \leq j \leq 3\:\!.
\ee
%
%
\mbox{} \vspace{-0.400cm} \\

\noindent
{\sc Remark:}
It is not difficult to show that Definitions 7.1 and 7.2
are equivalent for functions $\:\!v \in C_0^\infty \!\:\!$.
However,
Definition 7.2 has the advantage that, formally,
the function $\:\!v\:\!$ is not required to be differentiable.
Also, it has the form \\
\mbox{} \vspace{-0.630cm} \\
$$w \,=\, P_{\scriptscriptstyle \!\;\!H}v
\,=\: v \,+\;\! \nabla \phi $$
\mbox{} \vspace{-0.350cm} \\
with $\nabla \phi\:\!$
determined by the singular integral
in (7.13$b$),
to which the
Calderon--Zygmund theory of singular integrals
can be applied.
Using the estimate \\
\mbox{} \vspace{-0.750cm} \\
$$|\,K_{ij}(x-y)\,| \;\leq\;
 \f{C_1}{\,1\;\!+\;\!|\;\!x\;\!|^3} \qfor
 \;|\,y\,| \;\!\leq\;\! C_0 $$
\mbox{} \vspace{-0.230cm} \\
and large $|\;\!x\;\!|$,
it is easy to see
that \\
\mbox{} \vspace{-0.700cm} \\
$$|\,w(x)\,| \;\leq\;
 C\,\bigl(\,1 \:\!+\;\! |\;\!x\;\!|^3\;\!\bigr)^{-1}
 \qwith
 C \:\!=\, C(v), $$
\mbox{} \vspace{-0.300cm} \\
so that $\;\!w \in L^q(\R^3\!\;\!)\;\!$
if $1 < q \leq \infty $.
Here, as before,
we have assumed that
$\;\!w = P_{\scriptscriptstyle \!\;\!H}v$ and
$v \in C_0^\infty(\R^3\!\;\!)$.
Furthermore, we have $w \in C^\infty(\R^{3}\!\;\!) $
since the integrals in (7.13$b$)
can be differentiated under the integral sign
if all derivatives are moved through integration by parts
from the kernels to $v_i(y)$.
Thus,
we have shown the following:
if $1 < q \leq \infty$,
then
the Helmholtz projector $P_{\scriptscriptstyle \!\;\!H}$
maps $C_0^\infty(\R^{3}\!\;\!)$ into
$C^\infty(\R^3\!\;\!) \cap L^q(\R^3\!\;\!)$.

We state the following important theorem,
which can be obtained from the Calderon--Zygmund theory
(see e.g.$\;$\cite{Galdi1994}, p.$\;$109).

%
%
%
%
\begin{theorem}
\label{t.6.PH}
For $\;\!1 < q < \infty$,
there is a constant $\:\!C_{\!\;\!q}\!\;\!> 0\:\!$
with \\
\mbox{} \vspace{-0.700cm} \\
\be
\|\,P_{\scriptscriptstyle \!\;\!H} \:\!v\,\|_{\mbox{}_{\scriptstyle L^q}}
\;\!\leq\:
C_{\!\;\!q} \: \|\;v\;\|_{\mbox{}_{\scriptstyle L^q}}\!\;\!,
\quad \;\; v \in C_0^\infty\!\:\!.
\ee
\mbox{} \vspace{-0.330cm} \\
Since $\;\!C_0^\infty\!$ is dense in $L^q\!$,
the operator $\:\!P_{\scriptscriptstyle \!\;\!H}\!\:\!$
can be extended uniquely
as a bounded linear operator from $\:\!L^q\!\;\!$ into itself.
\end{theorem}
%
%

\end{subsection}

\end{section}

%
%
%
\mbox{} \vspace{-1.050cm} \\
%
%
%


%
%
\mbox{} \\
\mbox{} \\
\mbox{} \\
\mbox{} \\
{\sc Jens Lorenz} \\
Department of Mathematics \\
University of New Mexico \\
Albuquerque, NM 87131, USA \\
E-mail: {\tt lorenz@math.unm.edu} \\
\mbox{} \\
\mbox{} \vspace{-0.350cm} \\
{\sc Paulo R. Zingano} \\
Departamento de Matematica \\
Universidade Federal do Rio Grande do Sul \\
Porto Alegre, RS 91509, Brasil \\
E-mail: {\tt paulo.zingano@ufrgs.br} \\

%
%

\end{document}